\documentclass[12pt]{article}

\usepackage[dvips]{graphicx}
\usepackage{latexsym,amsmath,amsfonts,amscd}
\usepackage{epsfig}
\usepackage{changebar}
\usepackage{multirow}
\usepackage{verbatim}
\usepackage{graphicx}
\usepackage{amsmath}
\usepackage{amssymb}
\usepackage{color}

\newcommand{\bk}{\mathbf{k}}
\newcommand{\bx}{\mathbf{x}}
\newcommand{\ph}{\varphi}
\newcommand{\eps}{\varepsilon}
\newcommand{\bE}{\mathbf{E}}
\newcommand{\ti}{\mathrm{t}}
\newcommand{\nV}{\Psi}
\newcommand{\kbt}{k_B T_L}
\newcommand{\ak}{\alpha_K}
\newcommand{\hw}{\hbar \omega_{p}}
\newcommand{\mass}{m^{*}}
\newcommand{\qe}{\gamma}

\newfont{\iams}{msbm9}
\newcommand{\rind}{\mbox{\iams \symbol{'122}}}
\newcommand{\Itre}{\int_{\scriptstyle \rind^{3}}}
\newcommand{\Ipm}{\int_0^\pi \! \! d\ph' \int_{-1}^1 \! d\mu' \:}
\newcommand{\Iwmp}{\int_{0}^{+ \infty} \! \! dw \int_{-1}^1 \! d\mu
                   \int_0^\pi \! \! d\ph \:}
\newcommand{\Iwm}{\int_{0}^{+ \infty} \! \! dw \int_{-1}^1 \! d\mu
                    \:}
\newcommand{\argf}{( \mathbf{\cdot} )}
\newcommand{\ot}{\frac{1}{2}}
\newcommand{\dm}{\displaystyle}

\newcommand{\ix}{\int_{x_{i-\frac{1}{2}}}^{x_{i+\frac{1}{2}}}}
\newcommand{\iy}{\int_{y_{j-\frac{1}{2}}}^{y_{j+\frac{1}{2}}}}
\newcommand{\iw}{\int_{w_{k-\frac{1}{2}}}^{w_{k+\frac{1}{2}}}}
\newcommand{\imu}{\int_{\mu_{m-\frac{1}{2}}}^{\mu_{m+\frac{1}{2}}}}
\newcommand{\iphi}{\int_{\ph_{n-\frac{1}{2}}}^{\ph_{n+\frac{1}{2}}}}
\newcommand{\kpol}{\ell}
\begin{document}

\baselineskip=2pc

\begin{center}
{\bf A discontinuous Galerkin  solver for Boltzmann Poisson systems
in nano devices\footnote{Support from the Institute of Computational
Engineering and Sciences and the University of Texas Austin is
gratefully acknowledged.}}
\end{center}

\vspace{.10in}

\centerline{Yingda Cheng\footnote{E-mail: ycheng@math.utexas.edu.},
Irene M. Gamba\footnote{E-mail: gamba@math.utexas.edu. Research
supported by  NSF-0807712  and NSF-FRG-0757450.}}

\smallskip

\centerline{Department of Mathematics and ICES, University of Texas,
Austin, TX 78712}

\bigskip

\centerline{Armando Majorana\footnote{E-mail: majorana@dmi.unict.it.
Research supported by Italian PRIN 2006: Kinetic and continuum
models for particle transport in gases and semiconductors:
analytical and computational aspects.}}

\smallskip

\centerline{Dipartimento di Matematica e Informatica, Universit\`a
di Catania, Catania, Italy}

\bigskip

\centerline{\em and}

\bigskip

\centerline{Chi-Wang Shu\footnote{E-mail: shu@dam.brown.edu.
Research supported by NSF grant DMS-0809086 and DOE grant
DE-FG02-08ER25863.}}

\smallskip

\centerline{Division of Applied Mathematics, Brown University,
Providence, RI 02912}

\vspace{.4in}

\centerline{\bf Abstract}

\vspace{.16in}

In this paper, we present results of a discontinuous Galerkin (DG)
scheme applied to deterministic computations of the transients for
the Boltzmann-Poisson system describing  electron transport in
semiconductor devices. The collisional term models optical-phonon
interactions which become dominant under strong energetic conditions
corresponding to nano-scale active regions under applied bias. The
proposed numerical technique is a finite element method using
discontinuous piecewise polynomials as basis functions on
unstructured meshes.  It is applied to simulate hot electron
transport in bulk silicon, in a silicon $n^+$-$n$-$n^+$ diode and in
a double gated 12nm MOSFET. Additionally, the obtained results are
compared to those of a high order WENO scheme simulation and DSMC
(Discrete Simulation Monte Carlo) solvers.

\vfill


\newpage

\section{Introduction}

The evolution of the electron distribution function $f(\ti, \bx,
\bk)$ in semiconductors in dependence of time $\ti$, position $\bx$
and electron wave vector $\bk$ is governed by the Boltzmann
transport equation (BTE) \cite{lundstrom, Mark90, ferr91}
\begin{equation}
\label{bte} \frac{\partial f}{\partial \ti} + \frac{1}{\hbar}
\nabla_{\bk} \, \eps \cdot \nabla_{\bx} f -
 \frac{q}{\hbar} \bE \cdot \nabla_{\bk} f = Q(f) \, ,
\end{equation}
where $\hbar$ is the reduced Planck constant, and $q$ denotes the
positive elementary charge. The function $\eps(\bk)$ is the energy
of the considered crystal conduction band measured from the band
minimum; according to the Kane dispersion relation, $\eps$ is the
positive root of
\begin{equation}
\eps ( 1 + \alpha \eps) = \frac{\hbar^{2} k^{2}}{2 \mass} \, ,
\end{equation}
where $\alpha$ is the non-parabolicity factor and $\mass$ the
effective electron mass. The electric field $\bE$ is related to the
doping density $N_{D}$ and the electron density $n$, which equals
the zero-order moment of the electron distribution function $f$, by
the Poisson equation
\begin{equation}
\nabla_{\bx} \left[ \epsilon_{r}(\bx) \, \nabla_{\bx} V \right]
 = \frac{q}{\epsilon_{0}} \left[ n(\ti,\bx) - N_{D}(\bx) \right] ,
\quad \bE = - \nabla_{\bx} V \, , \label{eqPoiss}
\end{equation}
where $\epsilon_{0}$ is the dielectric constant of the vacuum,
$\epsilon_{r}(\bx)$ labels the relative dielectric function
depending on the material and $V$ the electrostatic potential. The
collision operator $Q(f)$ takes into account acoustic deformation
potential and optical intervalley scattering \cite{ziman,tomi}. For
low electron densities, it reads
\begin{equation}
Q(f)(\ti, \bx, \bk) = \Itre \left[ S(\bk',\bk) f(\ti, \bx, \bk') -
S(\bk,\bk') f(\ti, \bx, \bk) \right] d \bk' \label{ope_coll}
\end{equation}
with the scattering kernel
\begin{eqnarray}
  S(\bk,\bk') & = & (n_{q} + 1) \, K \, \delta(\eps(\bk') - \eps(\bk) + \hw)
  \nonumber \\
  & & \mbox{} + n_{q} \, K \, \delta(\eps(\bk') - \eps(\bk) - \hw)
   +  K_{0} \, \delta(\eps(\bk') - \eps(\bk)) \,
  \label{kernel}
\end{eqnarray}
and $K$ and $K_{0}$ being constant for silicon. The symbol $\delta$
indicates the usual Dirac distribution and $\omega_{p}$ is the
constant phonon frequency. Moreover,
$$
  n_{q} = \left[ \exp \left( \frac{\hw}{\kbt} \right) - 1 \right]^{-1}
$$
is the occupation number of phonons, $k_B$ is the Boltzmann constant
and $T_L$ is the constant lattice temperature.

Semiclassical description of electron flow in semiconductors thus is
an equation in six dimensions (plus time if the device is not in
steady state) for a truly 3-D device, and four dimensions for a 1-D
device. This  heavy computational cost explains why the BP system is
traditionally simulated by the Direct Simulation Monte Carlo (DSMC)
methods \cite{jaco89}. In recent years, deterministic solvers to the
BP system were proposed \cite{fate93, MP,
carr02,carr021,cgms03,carr03,cgms06, gm07}. These methods provide
accurate results which, in general, agree well with those obtained
from Monte Carlo (DSMC) simulations, often at a fractional
computational time.  Moreover, they can resolve transient details
for the \emph{pdf}, which are difficult to compute with DSMC
simulators.
The methods proposed in \cite{cgms03,cgms06} used weighted
essentially non-oscillatory (WENO) finite difference schemes to
solve the Boltzmann-Poisson system.  The advantage of the WENO
scheme is that it is relatively simple to code and very stable even
on coarse meshes for solutions containing sharp gradient regions.  A
disadvantage of the WENO finite difference method is that it
requires smooth meshes to achieve high order accuracy, hence it is
not very flexible for adaptive meshes.

On the other hand, motivated by the easy {\it hp-}adaptivity and
simple communication pattern of the discontinuous Galerkin (DG)
 methods, researchers have
worked on developing the DG method for solving the Boltzmann
equation and its macroscopic models \cite{Chen_95_JCP_Q,
Chen_95_VLSI, LS1,LS2, GP}. The type of DG method that we will
discuss here is a class of finite element methods originally devised
to solve hyperbolic conservation laws containing only first order
spatial derivatives, e.g. \cite{cs1,cs2,cs3,cs4,cs5}. Using
completely discontinuous polynomial space for both the test and
trial functions in the spatial variables and coupled with explicit
and nonlinearly stable high order Runge-Kutta time discretization,
the method has the advantage of flexibility for arbitrarily
unstructured meshes, with a compact stencil, and with the ability to
easily accommodate arbitrary {\it hp-}adaptivity. For more details
about DG scheme for convection dominated problems, we refer to the
review paper \cite{dgsurvey}. The DG method was later generalized to
the local DG (LDG) method to solve the convection diffusion equation
\cite{cs6} and elliptic equations \cite{abcm}. It is $L^2$ stable
and locally conservative, which makes it particularly suitable to
treat the Poisson equation.   In our previous work
\cite{chgms-sispad07,Cheng_08_JCE_BP}, we proposed the first DG
solver for (\ref{bte}) and showed some preliminary numerical
calculations for one- and two-dimensional devices.  In this paper,
we will carefully formulate the DG-LDG scheme for the
Boltzmann-Poisson system and perform extensive numerical studies to
validate our calculation.

This paper is organized as follows: in Section 2, we review the
change of variables in \cite{MP,carr03}. In Section 3, we study the
DG-BTE solver for 1D diodes. Section 4 is devoted to the discussion
of the 2D double gate MOSFET DG solver. Conclusions and final
remarks are presented in Section 5. Some technical implementations
of the DG solver are collected in the Appendix.

\section{Change of variables}
In this section, we review the change of variables proposed in
\cite{MP,cgms03}.

For the numerical treatment of the system (\ref{bte}),
(\ref{eqPoiss}), it is convenient to introduce suitable
dimensionless quantities and variables.
We assume $T_L = 300 \, K$. Typical values for length, time and
voltage are $\ell_* = 10^{-6} \, m$, $t_* = 10^{-12} \, s$ and $V_*
= 1 \, \mbox{Volt}$, respectively.  Thus, we define the
dimensionless variables
$$
  (x,y,z) = \frac{\bx}{\ell_*}  \, , \quad
  t = \frac{\ti}{t_*} \, , \quad
  \nV = \frac{V}{V_*} \, , \quad
  (E_x, E_y, E_z) = \frac{\bE}{E_*} \,
$$
with $E_* = 0.1 \, V_* \, \ell_{*}^{-1}$ and
$$
\dm
  E_{x} = - c_{v} \frac{\partial \nV}{\partial x} \, , \quad
  E_{y} = - c_{v} \frac{\partial \nV}{\partial y} \, , \quad
  c_v = \frac{V_*}{\ell_* E_*} \, .
$$
 In correspondence to \cite{MP} and \cite{carr03}, we perform a coordinate
transformation for $\bk$ according to
\begin{equation}\label{ktransform}
\bk = \frac{\sqrt{2 m^* \kbt}}{\hbar} \sqrt{w(1+\ak w)} \left(\mu,
\sqrt{1-\mu^2} \cos \ph, \sqrt{1-\mu^2} \sin \ph \right),
\end{equation}
where the new independent variables are the dimensionless energy
$\dm w = \frac{\eps}{\kbt}$, the cosine of the polar angle $\mu$ and
the azimuth angle $\ph$ with $\ak = \kbt \alpha$.
The main advantage of the generalized spherical coordinates
(\ref{ktransform}) is the easy treatment of the Dirac distribution
in the kernel (\ref{kernel}) of the collision term. In fact, this
procedure enables us to transform the integral operator
(\ref{ope_coll}) with the nonregular kernel $S$  into an
integral-difference operator, as shown in the following.

We are interested in studying two-dimensional problems in real space
but, of course, in the whole three-dimensional $\bk$-space.
Therefore, it is useful to consider the new unknown function $\Phi$
related to the electron distribution function via
$$
\Phi(t, x, y, w, \mu, \ph) = s(w) f(\ti, \bx, \bk) \, ,
$$
where
\begin{equation}
 s(w) = \sqrt{w(1+\ak w)}(1+2\ak w),
 \label{sw}
\end{equation}
is proportional to the Jacobian of the change of variables
(\ref{ktransform}) and, apart from a dimensional constant factor, to
the density of states.
This allows us to write the free streaming operator of the
dimensionless Boltzmann equation in a conservative form, which is
appropriate for applying standard numerical schemes used for
hyperbolic partial differential  equations. Due to the symmetry of
the problem and of the collision operator, we have
\begin{equation}
 \Phi(t, x, y, w, \mu, 2 \pi - \ph) = \Phi(t, x, y, w, \mu, \ph) \, .
 \label{symm}
\end{equation}
Straightforward but cumbersome calculations end in the following
transport e\-qua\-tion for $\Phi$:
\begin{equation}
\frac{\partial\Phi}{\partial t} + \frac{\partial}{\partial x} (g_1
\Phi) + \frac{\partial}{\partial y} (g_2 \Phi) +
\frac{\partial}{\partial w} (g_3 \Phi) + \frac{\partial}{\partial
\mu} (g_4 \Phi) + \frac{\partial}{\partial \ph} (g_5 \Phi) = C(\Phi)
\, . \label{eqPhi}
\end{equation}
The functions $g_i$ $(i=1,2,..,5)$ in the advection terms depend on
the independent variables $w$, $\mu$, $\ph$  as well as on time and
position
 via the electric field.
They are given by
\begin{eqnarray*}
g_1 \argf & = & c_x \frac{\mu \sqrt{w(1+\ak w)}}{1+2 \ak w} \, ,
\\
g_2 \argf & = & c_x \frac{\sqrt{1-\mu^2} \sqrt{w(1+\ak w)}
\cos\ph}{1+2 \ak w} \, ,
\\
g_3 \argf & = & \mbox{} - 2 c_k \frac{\sqrt{w(1+\ak w)}}{1+2 \ak w}
\left[ \mu \, E_x(t,x,y) + \sqrt{1-\mu^2} \cos\ph \, E_y(t,x,y)
\right]  ,
\\
g_4 \argf & = & \mbox{} - c_k \frac{\sqrt{1-\mu^2}}{\sqrt{w(1+\ak
w)}}
 \left[ \sqrt{1-\mu^2} \, E_x(t,x,y) - \mu \cos\ph \, E_y(t,x,y) \right] \, ,
\\
g_5 \argf & = & c_k \frac{\sin\ph}{\sqrt{w(1+\ak w)} \sqrt{1-\mu^2}}
\, E_y(t,x,y)
\end{eqnarray*}
with
$$
  \dm
  c_x = \frac{t_*}{\ell_*} \sqrt{\frac{2 \, \kbt}{\mass}}
  \mbox{ and} \quad
  c_k = \frac{t_* q E_*}{\sqrt{2 \mass \kbt}} \, .
$$
The right hand side of (\ref{eqPhi}) is the integral-difference
operator
\begin{eqnarray*}
&& C(\Phi)(t,x,y,w,\mu,\ph) =  s(w) \left\{ c_{0} \Ipm
\Phi(t,x,y,w,\mu ',\ph')
   \right. \\
&& \left. +  \Ipm [ c_{+} \Phi(t,x,y,w + \qe,\mu ',\ph')
   + c_{-} \Phi(t,x,y,w - \qe,\mu ',\ph') ] \right\}
\\[5pt]
&& \mbox{} - 2 \pi [c_0 s(w) +  c_+ s(w - \qe) + c_- s(w + \qe)]
   \Phi(t,x,y,w,\mu,\ph)  \, ,
\end{eqnarray*}
where
$$
\dm
  (c_0, c_+, c_-) = \frac{2 \mass \, t_*}{\hbar^3} \sqrt{2 \, \mass \, \kbt}
                    \left(K_0 , (n_{q} + 1) K , n_{q} K \right) ,
  \quad
   \qe = \frac{\hw}{\kbt} \,
$$
are dimensionless parameters. We remark that the $\delta$
distributions in the kernel $S$ have been eliminated which leads to
the shifted arguments of $\Phi$. The para\-me\-ter $\qe$ represents
the jump constant corresponding to the quantum of energy $\hw$. We
have also taken into account (\ref{symm}) in the integration with
respect to $\ph'$. Since the energy variable $\omega$ is not
negative, we must consider null $\Phi$ and the function $s$, if the
argument $\omega - \gamma$ is negative.

In terms of the new variables the  electron density becomes
$$
  \dm n(t_* t, \ell_* x, \ell_* y)
   = \Itre f(t_* t, \ell_* x, \ell_* y, \bk) \: d \bk
   = \left( \frac{\sqrt{2 \,\mass \kbt }}{\hbar} \right)^{\! \! 3}
   \rho(t,x,y) \, ,
$$
where
\begin{equation}
  \rho(t,x,y) = \Iwmp \Phi (t,x,y,w,\mu,\ph) \, .
\label{dens}
\end{equation}
Hence, the dimensionless Poisson equation writes
\begin{equation}
\label{pois}
 \frac{\partial}{\partial x} \left( \epsilon_{r} \frac{\partial \nV}{\partial x}
 \right)
 +
 \frac{\partial}{\partial y} \left( \epsilon_{r} \frac{\partial \nV}{\partial y}
 \right)
 = c_{p} \left[ \rho(t,x,y) - \mathcal{N}_{D}(x,y) \right]
\end{equation}
with
$$
  \mathcal{N}_{D}(x,y) =
    \left( \frac{\sqrt{2 \,\mass \kbt }}{\hbar} \right)^{\! \! -3}
    N_{D}(\ell_* x, \ell_* y) \,  \mbox{ and }
  c_p = \left( \frac{\sqrt{2 \,\mass \kbt }}{\hbar} \right)^{\! \! 3}
        \frac{\ell_*^{2} q}{\epsilon_{0}} \, .
$$
Choosing the same values of the physical parameters as in \cite{MP},
we  obtain

\begin{center}
\begin{tabular}[t]{|l|l|l|}
\hline $c_0 \approx 0.26531$ & $c_x \approx 0.16857$ & $c_p \approx
1830349.$
\\[5pt]
$c_+ \approx 0.50705$ & $c_k \approx 0.32606$ & $c_v = 10.$
\\[5pt]
$c_- \approx 0.04432$ & $ \qe \approx  2.43723$ & $ \ak \approx
0.01292$
\\[5pt]
$ \epsilon_{r} = 11.7$ & &
\\[5pt]
\hline
\end{tabular}
\end{center}
Moreover, the dimensional $x$-component of the velocity is given by
$$
\frac{\dm \Iwmp g_{1}(w,\mu) \, \Phi (t,x,y,w,\mu,\ph)}{\rho(t,x,y)}
\, ,
$$
the dimensional density by
$$
1.0115 \times 10^{26} \times \rho(t,x,y) \, ,
$$
and the energy by
$$
0.025849 \times \frac{\dm \Iwmp w \, \Phi
(t,x,y,w,\mu,\ph)}{\rho(t,x,y)} \, .
$$

In some simplified models, we consider our device in the
$x-$direction by assuming that the doping profile, the potential and
thus the force field are only $x-$dependent. By cylindrical
symmetry, the resulting distribution function does not depend on
$\ph$. In this case, the Boltzmann transport equation is reduced to

\begin{equation}
\frac{\partial\Phi}{\partial t} + \frac{\partial}{\partial x} (g_1
\Phi) + \frac{\partial}{\partial w} (g_3 \Phi) +
\frac{\partial}{\partial \mu} (g_4 \Phi)  = C(\Phi) \, .
\label{bte1d}
\end{equation}
with
\begin{eqnarray*}
g_1 \argf & = & c_x \frac{\mu \sqrt{w(1+\ak w)}}{1+2 \ak w} \, ,
\\
g_3 \argf & = & \mbox{} - 2 c_k \frac{\sqrt{w(1+\ak w)}}{1+2 \ak w}
\, \mu \, E(t,x)  \, ,
\\
g_4 \argf & = & \mbox{} - c_k \frac{1-\mu^2}{\sqrt{w(1+\ak w)}}
 \, E(t,x)   \,
\end{eqnarray*}
and
\begin{eqnarray*}
&& C(\Phi)(t,x,w,\mu) =  s(w) \left\{ c_{0} \pi \int_{-1}^1 \! d\mu'
\:\Phi(t,x,w,\mu ')
   \right. \\
&& \left. +  \pi \int_{-1}^1 \! d\mu' \: [ c_{+} \Phi(t,x,w +
\qe,\mu ')
   + c_{-} \Phi(t,x,w - \qe,\mu ') ] \right\}
\\[5pt]
&& \mbox{} - 2 \pi [c_0 s(w) +  c_+ s(w - \qe) + c_- s(w + \qe)]
   \Phi(t,x,w,\mu)  \, ,
\end{eqnarray*}

In terms of the new variables the  electron density becomes
$$
  \dm n(t_* t, \ell_* x)
   = \Itre f(t_* t, \ell_* x,  \bk) \: d \bk
   = \left( \frac{\sqrt{2 \,\mass \kbt }}{\hbar} \right)^{\! \! 3}
   \rho(t,x) \, ,
$$
where
\begin{equation}
  \rho(t,x) = \pi \int_{0}^{+ \infty} \! \! dw \int_{-1}^1 \! d\mu
                   \: \Phi (t,x,w,\mu) \, .
\end{equation}
Hence, the dimensionless Poisson equation writes
\begin{equation}
\label{pois1d}
 \frac{\partial}{\partial x} \left( \epsilon_{r} \frac{\partial \nV}{\partial x}
 \right)
  = c_{p} \left[ \rho(t,x) - \mathcal{N}_{D}(x) \right]
\end{equation}
and
$$ E = - c_{v} \frac{\partial \nV}{\partial x}. $$

\section{DG-BTE solver for 1D diodes simulation}

We begin with formulating the DG-BTE solver for  1D diodes. These
examples have been thoroughly studied and tested by WENO in
\cite{cgms03}.

The Boltzmann-Poisson system (\ref{bte1d}) and (\ref{pois1d}) will
be solved on the domain
$$
x \in [0,L], \quad w \in[0, w_{\mbox{max}}], \quad \mu\in[-1,1],
$$
where $L$ is the dimensionless length of the device and
$w_{\mbox{max}}$ is the maximum value of the energy, which is
adjusted in the numerical experiments such that
$$
\Phi(t,x,w,\mu) \approx 0 \qquad \mbox{for}\, w \geq  w_{\mbox{max}}
\quad \mbox{and every} \quad t, x, \mu .
$$

In (\ref{bte1d}), $g_1$ and $g_3$ are completely smooth in the
variable $w$ and $\mu$, assuming $E$ is given and smooth. However,
$g_4$ is singular for the energy $w=0$, although it is compensated
by the $s(w)$ factor in the definition of $\Phi$.

The initial value of $f$ is a locally Maxwellian distribution at the
temperature $T_L$,
$$
\Phi(0,x,w, \mu)=s(w) N_D(x) e^{-w} \mathcal{M}
$$
with the numerical parameter $\mathcal{M}$ chosen so that the
initial value for the density is equal to the doping $N_D(x)$.

We choose to perform our calculations on the following rectangular
grid,
\begin{equation}
\label{grid1}
 \Omega_{ikm} = \left[ x_{i - \ot} , \, x_{i + \ot} \right] \times
                  \left[ w_{k - \ot} , \, w_{k + \ot} \right] \times
                  \left[ \mu_{m - \ot} , \, \mu_{m + \ot} \right]
\end{equation}
where $i=1, \ldots N_x$,  $k=1, \ldots N_w$, $m=1, \ldots N_\mu$,
and
$$
 x_{i \pm \ot} = x_{i} \pm \frac{\Delta x_{i}}{2} \, , \quad
 w_{k \pm \ot} = w_{k} \pm \frac{\Delta w_{k}}{2}\,  , \quad
 \mu_{m \pm \ot} = \mu_{m} \pm \frac{\Delta \mu_{m}}{2}\, .
$$
It is useful that we pick $N_\mu$  to be even, so the function $g_1$
will assume a constant sign in each cell $\Omega_{ikm}$.

The approximation space is thus defined as
\begin{equation}
V_h^{\kpol} = \{ v : v|_{\Omega_{ikm}} \in
P^{\kpol}(\Omega_{ikm})\},
\end{equation}
where $P^{\kpol}(\Omega_{ikm})$ is the set of all polynomials of
degree at most $\kpol$ on $\Omega_{ikm}$.
The DG formulation for the Boltzmann equation (\ref{bte1d}) would
be: to find $\Phi_h \in V_h^\kpol$, such that
\begin{eqnarray}
\label{dgb1D} &&\int_{\Omega_{ikm}} (\Phi_h)_t  \,v_h  \,d \Omega -
\int_{\Omega_{ikm}} g_1 \Phi_h  \,(v_h)_x  \,d \Omega
- \int_{\Omega_{ikm}} g_3 \Phi_h  \,(v_h)_w  \,d \Omega \nonumber \\
&&- \int_{\Omega_{ikm}} g_4 \Phi_h  \,(v_h)_\mu  \,d \Omega +F_x^+ -
F_x^- +F_w^+-F_w^- +F_\mu^+-F_\mu^-
\\&&
=\int_{\Omega_{ikm}} C(\Phi_h)  \,v_h \,d \Omega . \nonumber
\end{eqnarray}
for any test function $v_h \in V_h^\kpol$. In (\ref{dgb1D}),
$$
F_x^+=\iw \imu  g_1 \, \check{\Phi} \, v_h^-(x_{i+ \ot},  w, \mu)
dw \, d\mu  ,
$$
$$
F_x^-= \iw \imu  g_1 \, \check{\Phi} \, v_h^+ (x_{i- \ot}, w, \mu)
dw \, d\mu  ,
$$

$$
F_w^+=\ix \imu   g_3 \, \hat{ \Phi} \, v_h^- (x, w_{k+ \ot}, \mu )dx
\,  d\mu ,
$$
$$
F_w^-=\ix  \imu  g_3 \, \hat{ \Phi} \, v_h^+  (x,w_{k- \ot},\mu )dx
\,  d\mu ,
$$

$$
F_\mu^+=\ix  \iw  g_4 \, \tilde{\Phi} \, v_h^- (x, w, \mu_{m+ \ot}
)dx \,  dw  ,
$$
$$
F_\mu^-=\ix  \iw  g_4 \, \tilde{\Phi} \, v_h^+ ( x,  w,\mu_{m-
\ot})dx \,  dw ,
$$

where the upwind numerical fluxes $\check{\Phi},  \hat{ \Phi} ,
\tilde{\Phi}$ are chosen according to the following rules,
\begin{itemize}
\item The sign of $g_1$ only depends on $\mu$,
if $\mu_m >0$,  $ \check{\Phi} =\Phi^-$; otherwise, $ \check{\Phi}
=\Phi^+.$

\item The sign of $g_3$ only depends on $\mu E(t,x) $,
if $\mu_m E(t,x_i) <0$,   $ \hat{\Phi} =\Phi^-$; otherwise, $
\hat{\Phi} =\Phi^+.$

\item The sign of $g_4$ only depends on $E(t,x) $,
if $ E(t,x_i) <0$,  $ \tilde{\Phi} =\Phi^-$; otherwise, $
\tilde{\Phi} =\Phi^+.$

\end{itemize}

At the source and drain contacts, we implement the same boundary
condition as proposed in \cite{cgms06} to realize neutral charges.
In the $(w,\mu)$-space, non boundary condition is necessary, since
\begin{itemize}
\item at $w=0$, $g_3=0$. At $w=w_{\mbox{max}}$, $\Phi$ is machine zero.
\item At $\mu=\pm 1$, $g_4=0$,
\end{itemize}
 $F_w^+, F_w^-, F_\mu^+, F_\mu^- $ are always zero. This saves us the effort
of constructing ghost elements in comparison with WENO.

The Poisson equation (\ref{pois1d}) is solved by the LDG method on a
consistent grid of (\ref{grid1}) in the $x-$direction. It involves
rewriting the equation into the following form,
\begin{equation}
\label{pois1} \left\{\begin{array} {l}
        \displaystyle    q= \frac{\partial \nV}{\partial x}   \\
    \displaystyle     \frac{\partial}{\partial x} \left( \epsilon_{r} q \right)
 = R(t,x)
         \end{array}
   \right.
\end{equation}
where $ R(t,x)=c_{p} \left[ \rho(t,x) - \mathcal{N}_{D}(x) \right]$
is a known function that can be computed at each time step once
$\Phi$ is solved from (\ref{dgb1D}), and the coefficient
$\epsilon_r$ here is a constant. The grid we use is $I_{i}=\left[
x_{i - \ot} , \, x_{i + \ot} \right]$, with $i=1,\ldots, N_x$. The
approximation space  is
$$
W_h^\kpol=\{ v : v|_{I_{i}} \in P^\kpol(I_{i})\},
$$
with $P^\kpol(I_{i})$ denoting the set of all polynomials of degree
at most $\kpol$ on $I_{i}$.  The LDG scheme for (\ref{pois1}) is
given by: to find $q_h, \nV_h \in V_h^\kpol$, such that
\begin{eqnarray}
\label{ldgpois1}
 \int_{I_{i}} q_h v_h dx + \int_{I_{i}} \nV_h (v_h)_x dx
-\hat{\nV}_h v_h^-( x_{i + \ot})
+\hat{\nV}_h v_h^+( x_{i - \ot})  =0, &&\nonumber \\
 -\int_{I_{i}} \epsilon_{r} q_h (p_h)_x dx
+\widehat{ \epsilon_{r} q}_h p_h^-( x_{i + \ot}) -
\widehat{\epsilon_{r}  q}_h p_h^+( x_{i - \ot}) =\int_{I_{i}} R(t,x)
p_h dx &&
\end{eqnarray}
hold true for any $v_h, p_h \in W_h^\kpol$. In the above
formulation,   the flux is chosen as follows,
 $\hat{\nV}_h=\nV^-_h$, $\widehat{ \epsilon_{r} q}_h=\epsilon_{r} q_h^+ -[\nV_h]$, where $[\nV_h]=\nV_h^+ - \nV_h^-$.
 At $x=L$ we need to flip the flux to  $\hat{\nV}_h=\nV^+_h$, $\widehat{ \epsilon_{r} q}_h=\epsilon_{r} q_h^- -[\nV_h]$ to adapt to the Dirichlet boundary conditions. Solving  (\ref{ldgpois1}), we can obtain the numerical approximation of the  electric potential $\nV_h$ and electric field $E_h = - c_{v} q_h$ on each cell $I_i$.

To summarize, start with an initial condition for $\Phi_h$, the
DG-LDG algorithm advances from $t^n$ to $t^{n+1}$ in the following
steps:
\begin{description}
\item[Step 1] Compute $\rho_h(t,x)= \pi \, \Iwm \Phi_h (t,x,w,\mu)$.
\item[Step 2] Use $\rho_h(t,x)$ to solve from (\ref{ldgpois1})  the electric field, and compute $g_i$, $i=1, 3, 4$.
\item[Step 3] Solve (\ref{dgb1D}) and get a method of line ODE for $\Phi_h$.
\item[Step 4]Evolve this ODE by proper time stepping from $t^n$ to $t^{n+1}$, if partial time step is necessary, then repeat Step 1 to 3 as needed.
\end{description}

\begin{figure}[htb]
\centering
\includegraphics[width=4in,angle=0]{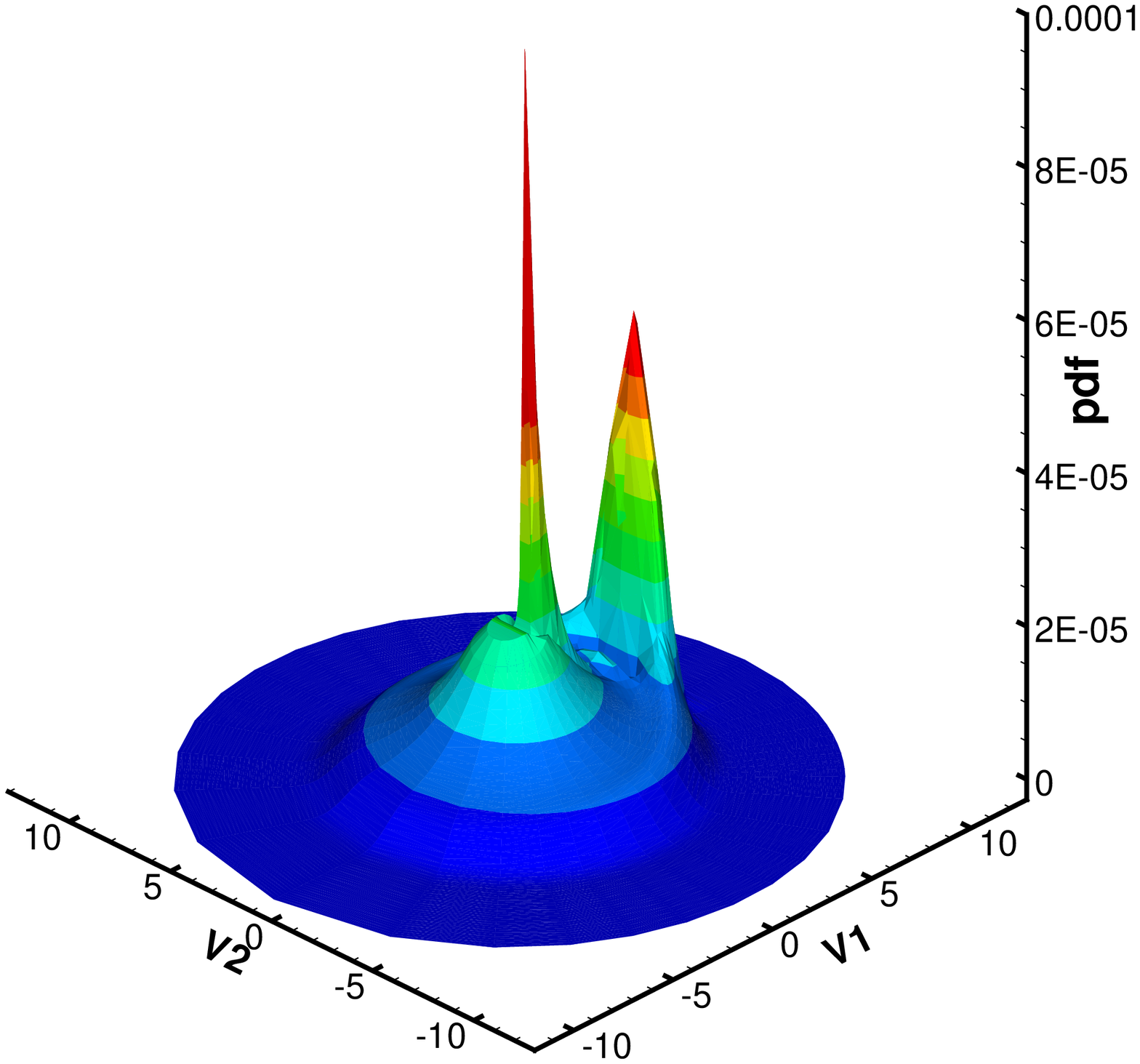}
\caption{PDF for 50nm channel at $t=3.0$, $x=0.149$.}
\label{50cart1}
\end{figure}

\begin{figure}[htb]
\centering
\includegraphics[width=4in,angle=0]{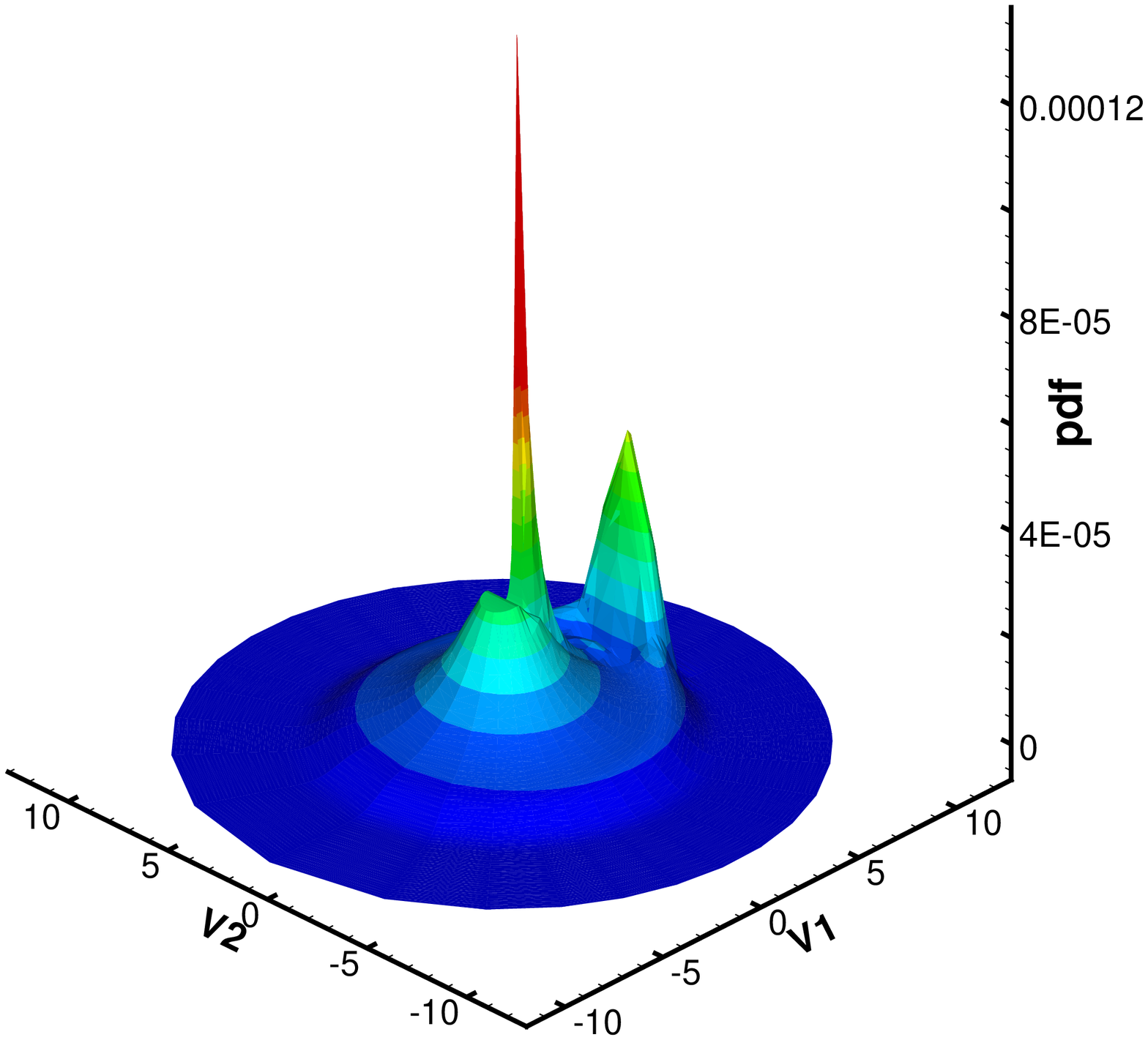}
\caption{PDF for 50nm channel at $t=3.0$, $x=0.15$.} \label{50cart2}
\end{figure}

\begin{figure}[htb]
\centering
\includegraphics[width=4in,angle=0]{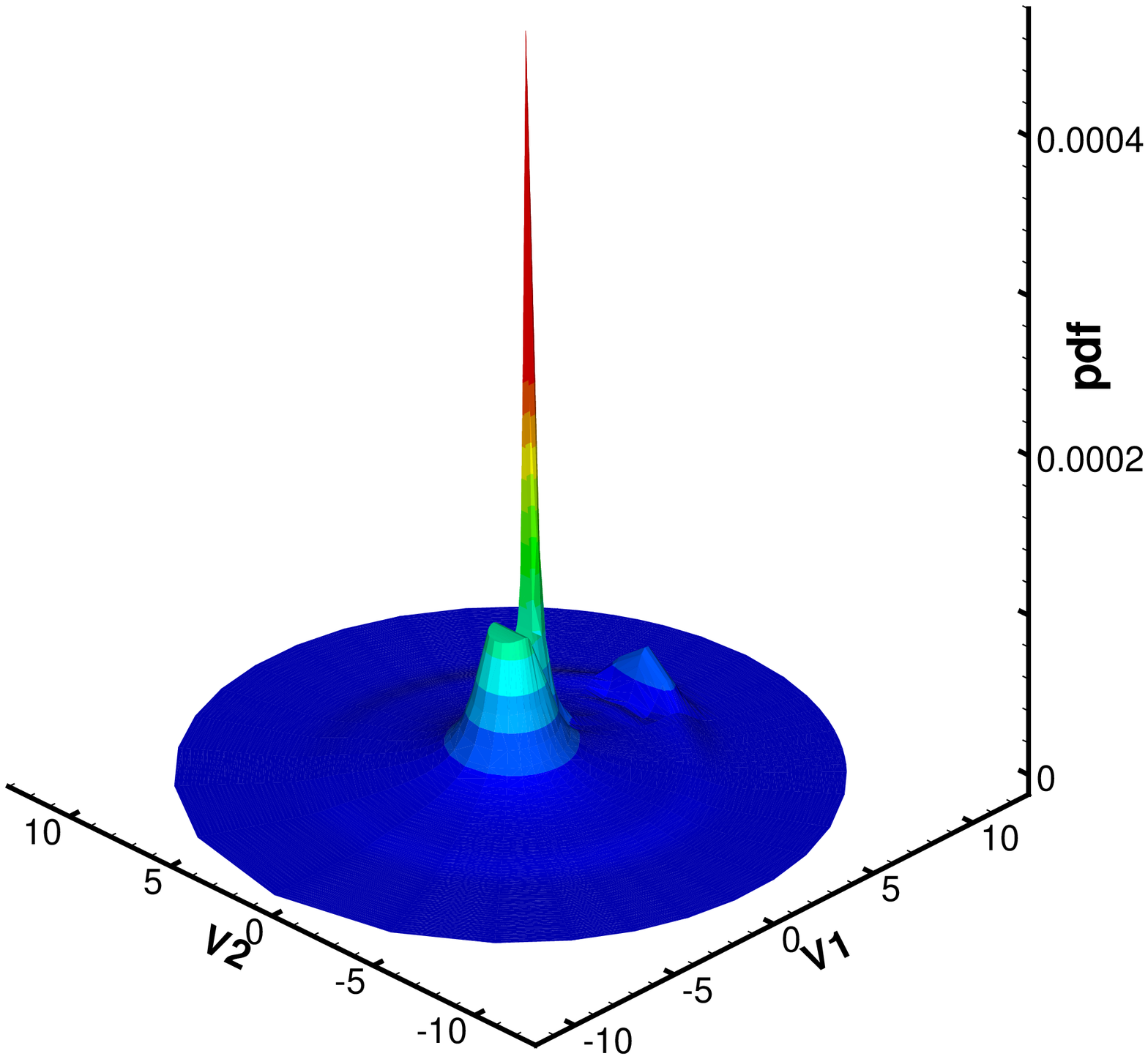}
\caption{PDF for 50nm channel at $t=3.0$, $x=0.152$.}
\label{50cart3}
\end{figure}

\begin{figure}[htb]
\centering
\includegraphics[width=4in,angle=0]{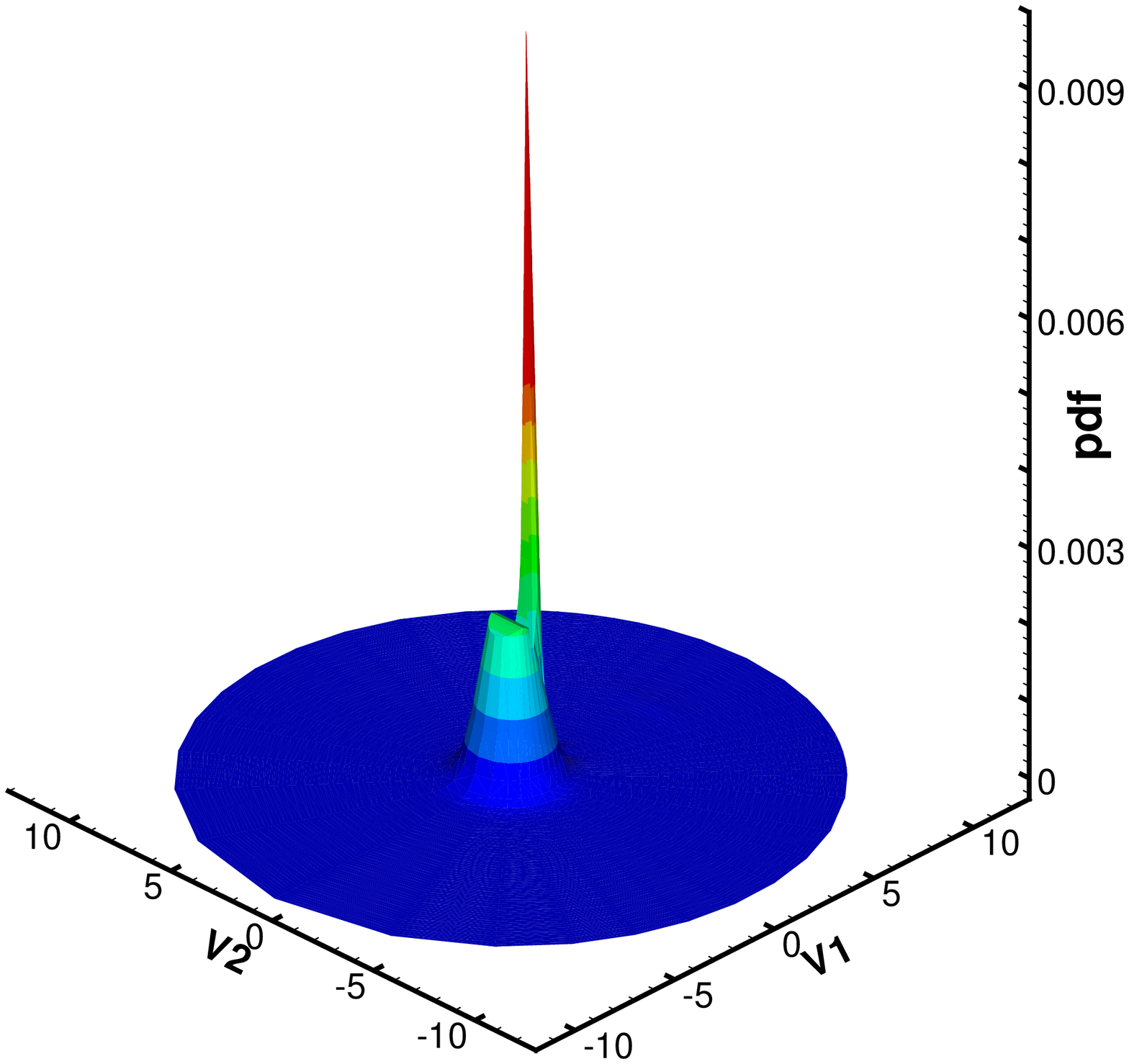}
\caption{PDF for 50nm channel at $t=3.0$, $x=0.16$.} \label{50cart4}
\end{figure}

We want to remark that, unlike WENO, the DG formulation above has no
restriction on the mesh size.
 In fact,  nonuniform meshes would be more desirable in practice. For small semiconductor devices, it is common to have nonsmooth doping profiles and
strong applied electric fields. The nonsmooth doping profile will
create a  distribution function $f$ with high densities in some
regions but low densities in other regions. Only a nonuniform grid
may guarantee accurate results without using a large number of grid
points. The strong electric fields give high energy to the charge
particles. Hence the distribution function has, for some fixed
points in the physical domain, a shape that is very different from
the Maxwellian distribution (the equilibrium distribution in the
absence of electric field), see for example Figures \ref{50cart1} to
\ref{50cart4}. Moreover, taking into account the exponential decay
of $f$ for large value of $|{\bf k}|$, a nonuniform grid can save us
tremendous amount of computational time without sacrificing accuracy
of the calculation. In the following simulations, we   use a
nonuniform mesh and refine locally near the junction of the channel
and near $\mu=1$ where most of the interesting phenomena happen.

We  consider two test examples: Si $n^+ -n-n^+$ diodes of a total
length of $1$ and $0.25 \mu m$, with $400$ and $50 nm$ channels
located in the middle of the device, respectively. For the $400 nm$
channel device, the dimensional doping is given by $N_D=5 \times
10^{17} cm^{-3}$ in the $n^+$ region and $N_D=2 \times 10^{15}
cm^{-3}$ in the $n^-$ region. For the $50 nm$ channel device, the
dimensional doping is given by $N_D=5 \times 10^{18} cm^{-3}$ in the
$n^+$ region and $N_D=1 \times 10^{15} cm^{-3}$ in the $n^-$ region.
Both examples were computed by WENO in \cite{carr03}.

In our simulation, we use piecewise linear polynomials, i.e.
$\ell=1$, and second-order Runge-Kutta time discretization. The
doping $N_D$ is smoothened in the following way near the channel
junctions to obtain non-oscillatory solutions. Suppose $N_D=N_D^+$
in the $n^+$ region,  $N_D=N_D^-$ in the $n^-$ region and the length
of the transition region is 2 cells, then the smoothened function is
$(N_D^+-N_D^-)(1-y^3)^3+N_D^-$, where $y=(x-x_0+\triangle x)/(2
\triangle x+10^{-20})$ is the coordinate transformation that makes
the transition region $(x_0-\triangle x, x_0+\triangle x)$ varies
from $0$ to $1$ in $y$.

The nonuniform mesh we use for $400 nm$ channels is defined as
follows. In the $x$-direction, if $x<0.2$ or $x>0.4$, $\triangle x=
0.01$. In the region $0.2<x<0.4$, $\triangle x= 0.005$. Thus, the
total number of cells in $x$ direction is $120$. In the
$w$-direction, we use $60$ uniform cells. In the $\mu$-direction, we
use $24$ cells, $12$ in the region $\mu<0.7$, $12$ in $\mu>0.7$.
Thus, the grid  consists of  $120 \times 60 \times 24$ cells,
compared to the WENO grid of $180 \times 60 \times 24$ uniform
cells.

\begin{figure}[htb]
\centering
\includegraphics[width=3in,angle=0]{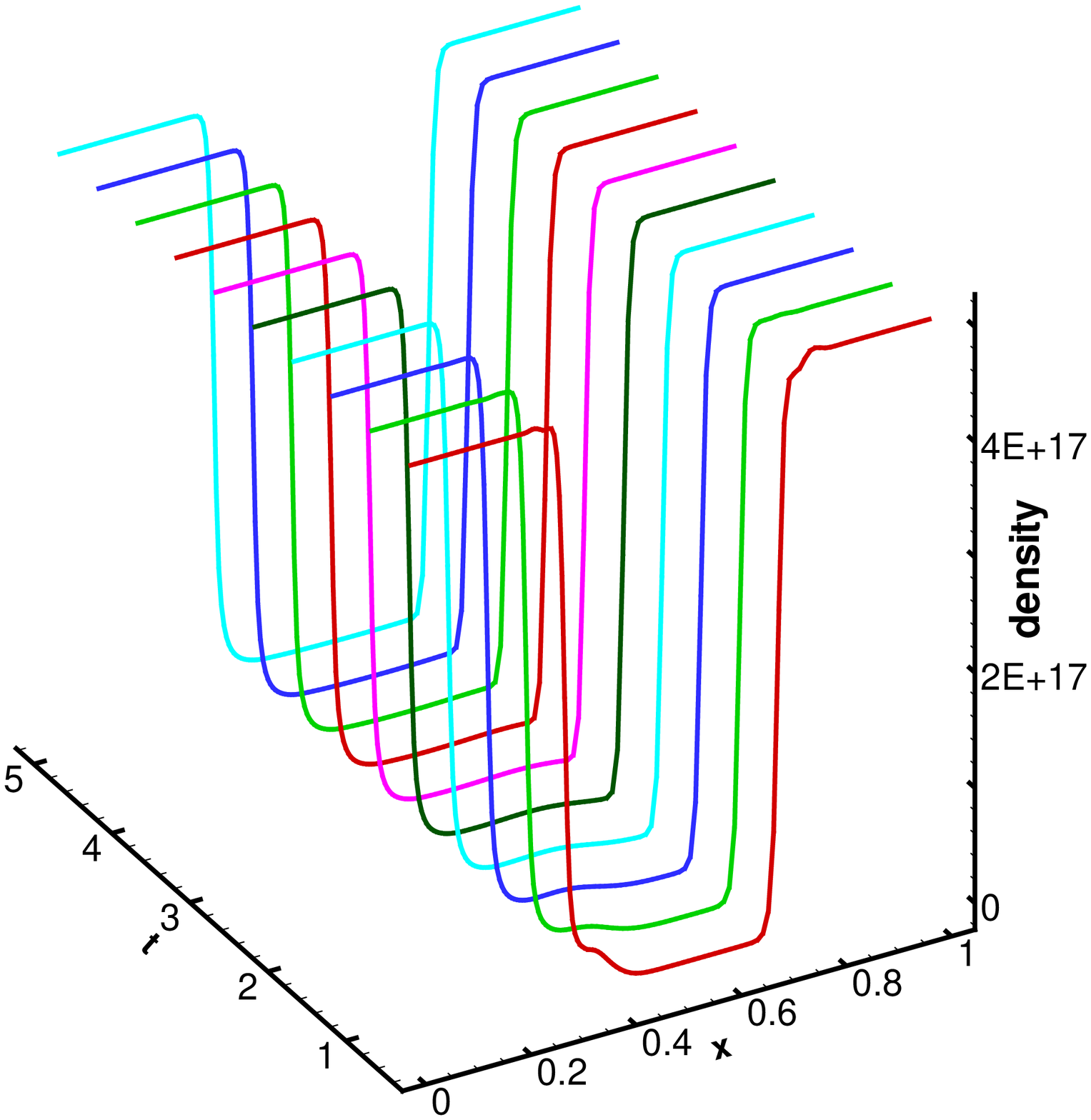}
\includegraphics[width=3in,angle=0]{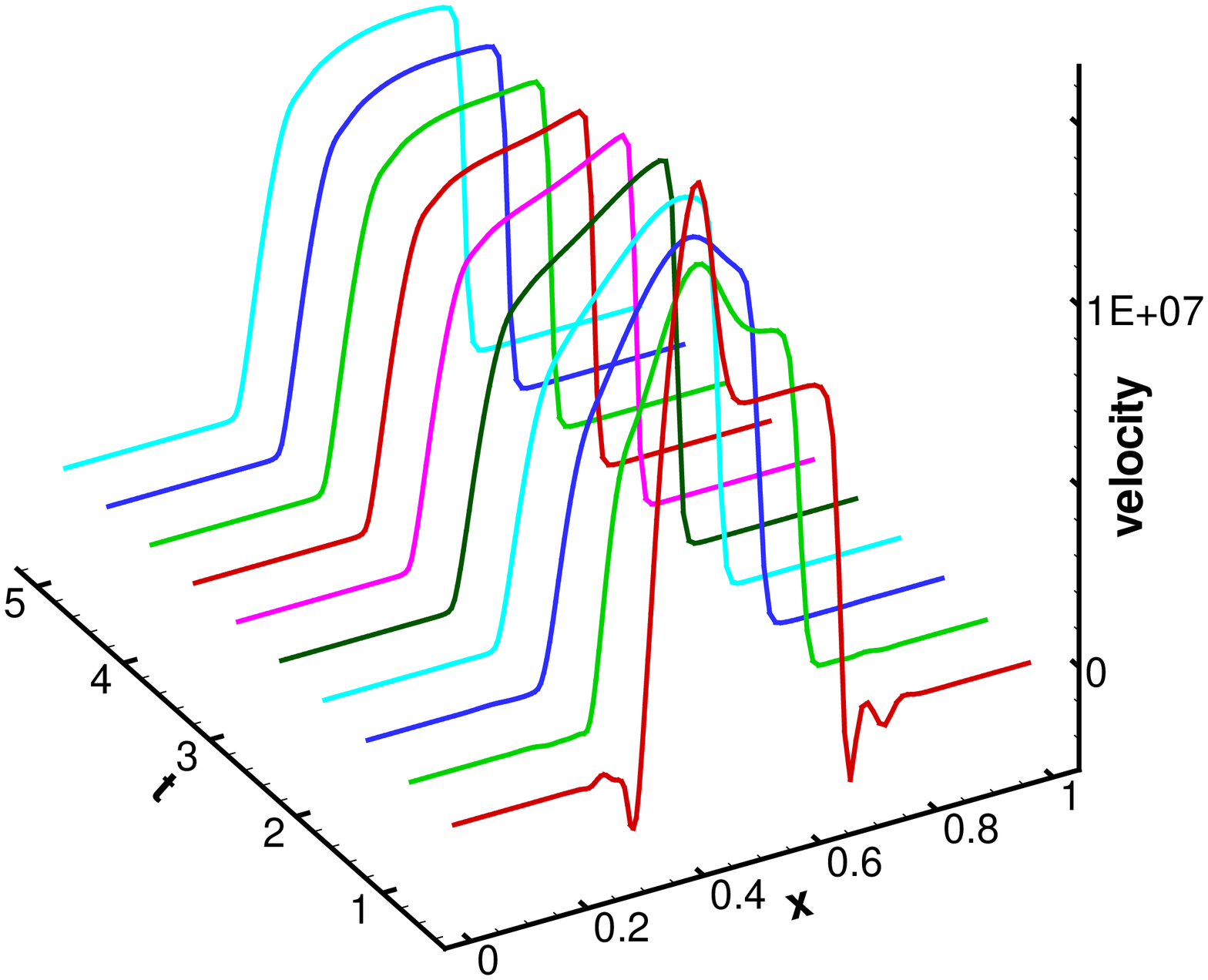}\\
\includegraphics[width=3in,angle=0]{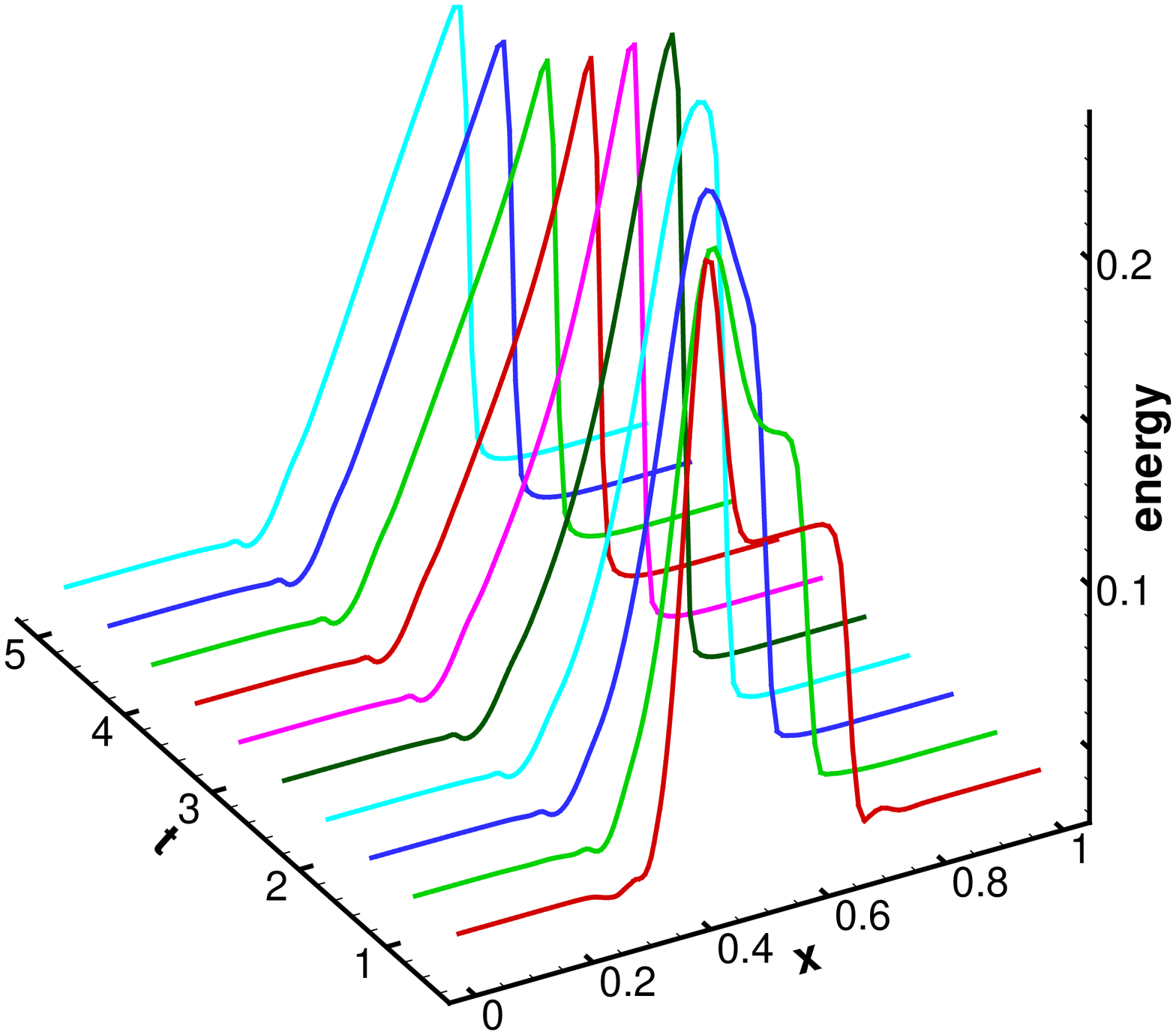}
\includegraphics[width=3in,angle=0]{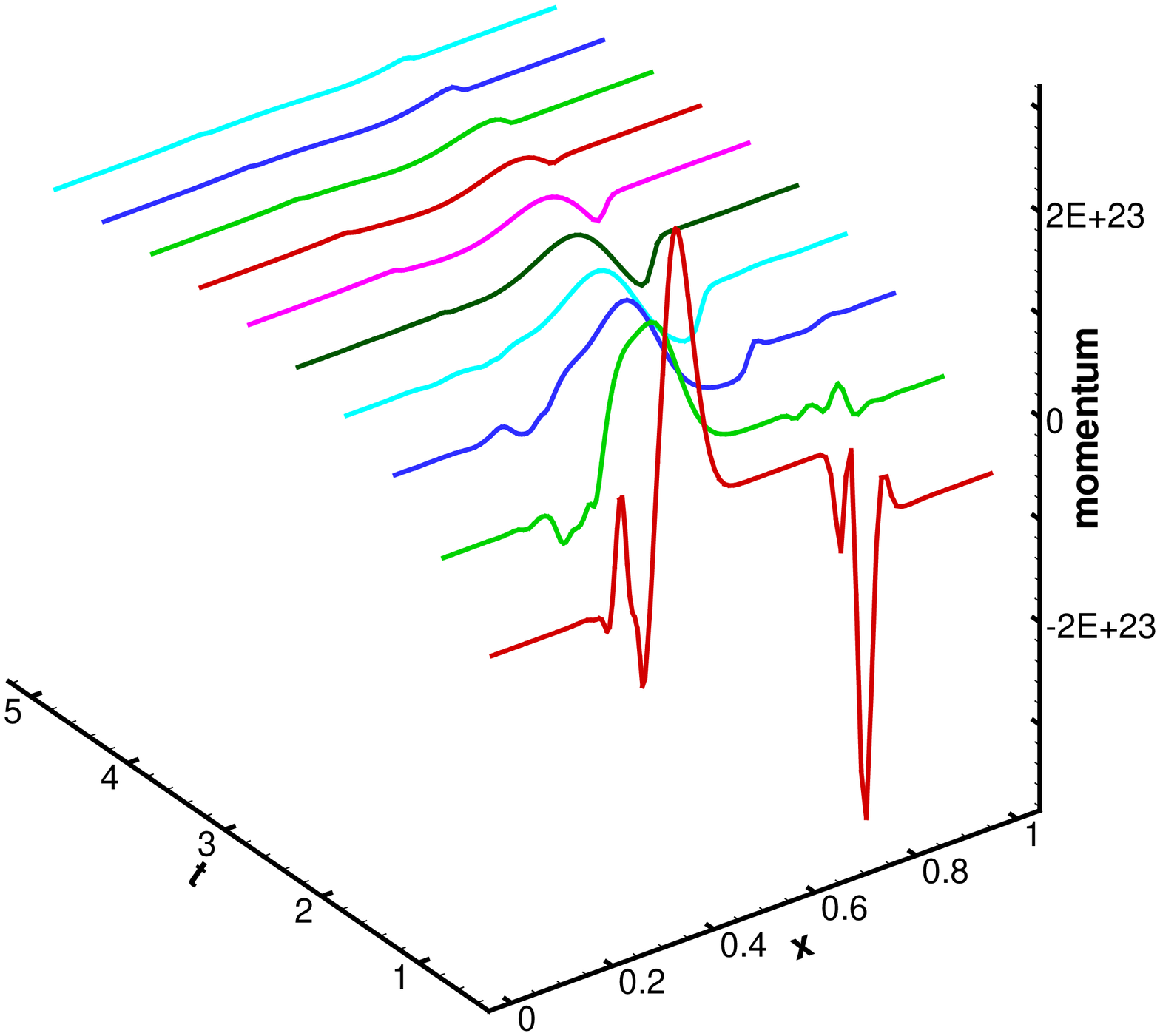}
\caption{Time evolution of macroscopic quantities using DG method
for $400$nm channel at  $V_{\mbox{bias}}=1.0$. Top left: density in
${cm}^{-3}$ ; top right: mean velocity in $cm/s$; bottom left:
energy in $eV$;  bottom right: momentum in ${cm}^{-2} \, s^{-1}$. }
\label{400evo}
\end{figure}

\begin{figure}[htb]
\centering
\includegraphics[width=2.93in,angle=0]{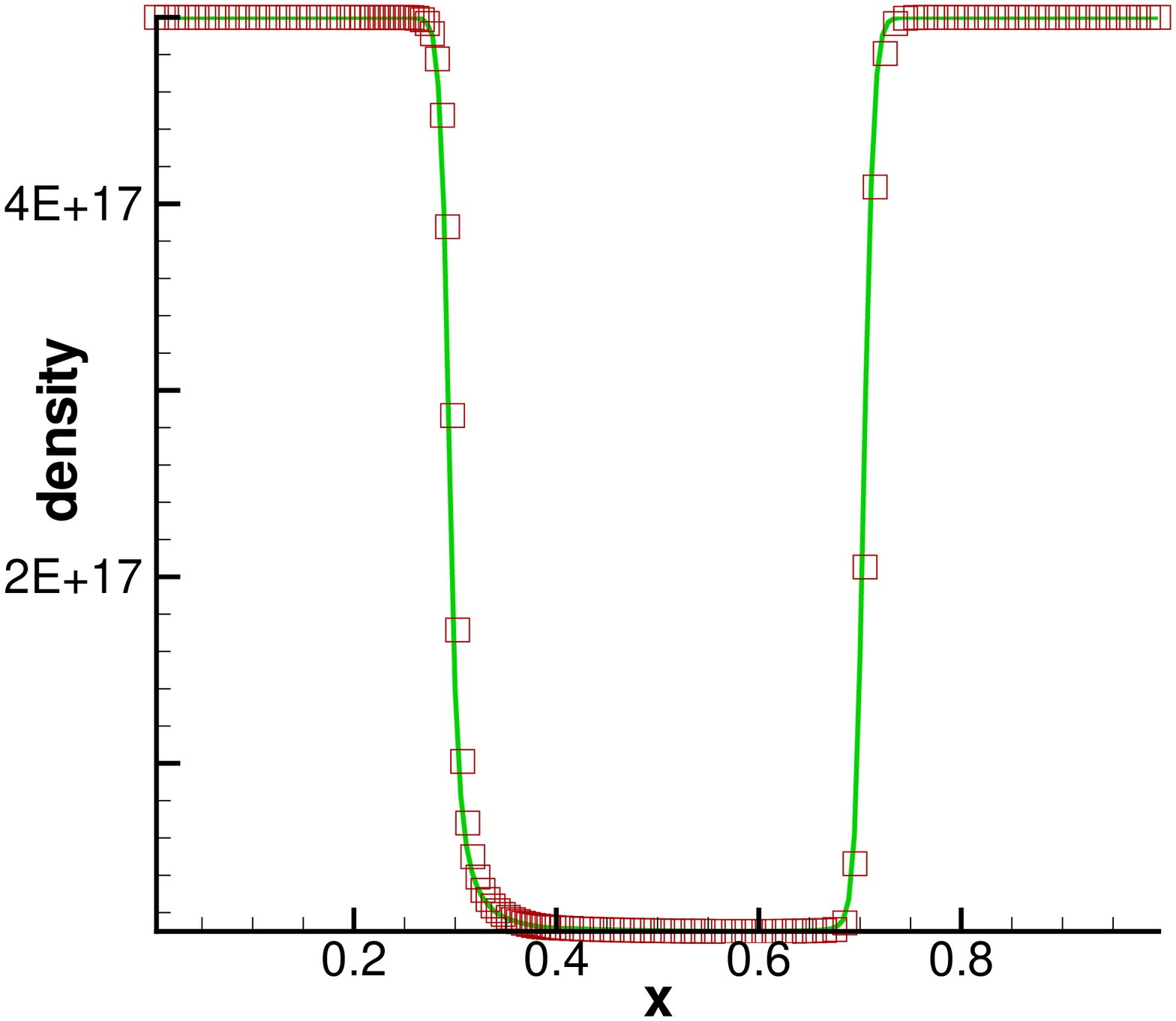}
\includegraphics[width=2.93in,angle=0]{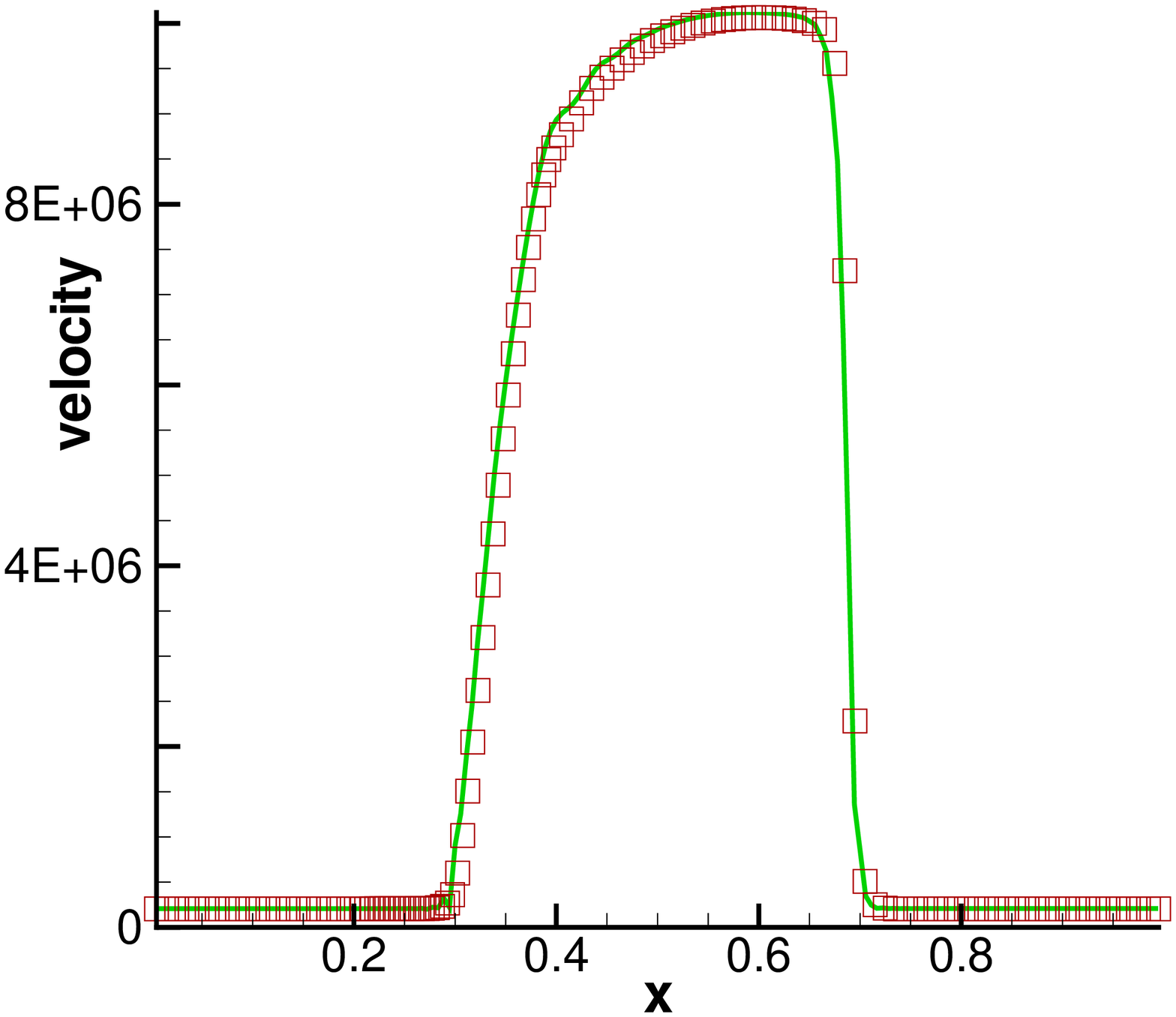}\\
\includegraphics[width=2.93in,angle=0]{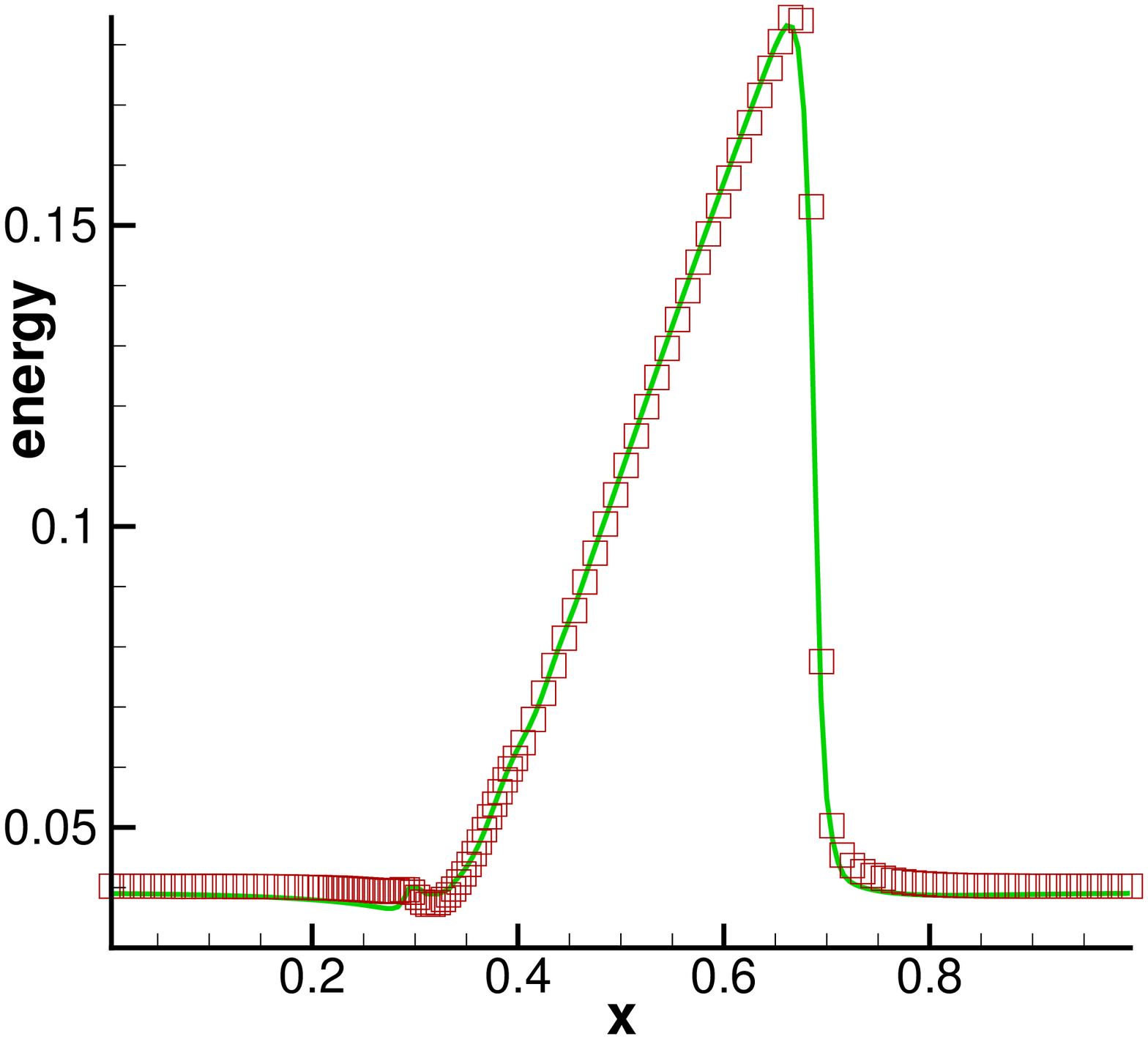}
\includegraphics[width=2.93in,angle=0]{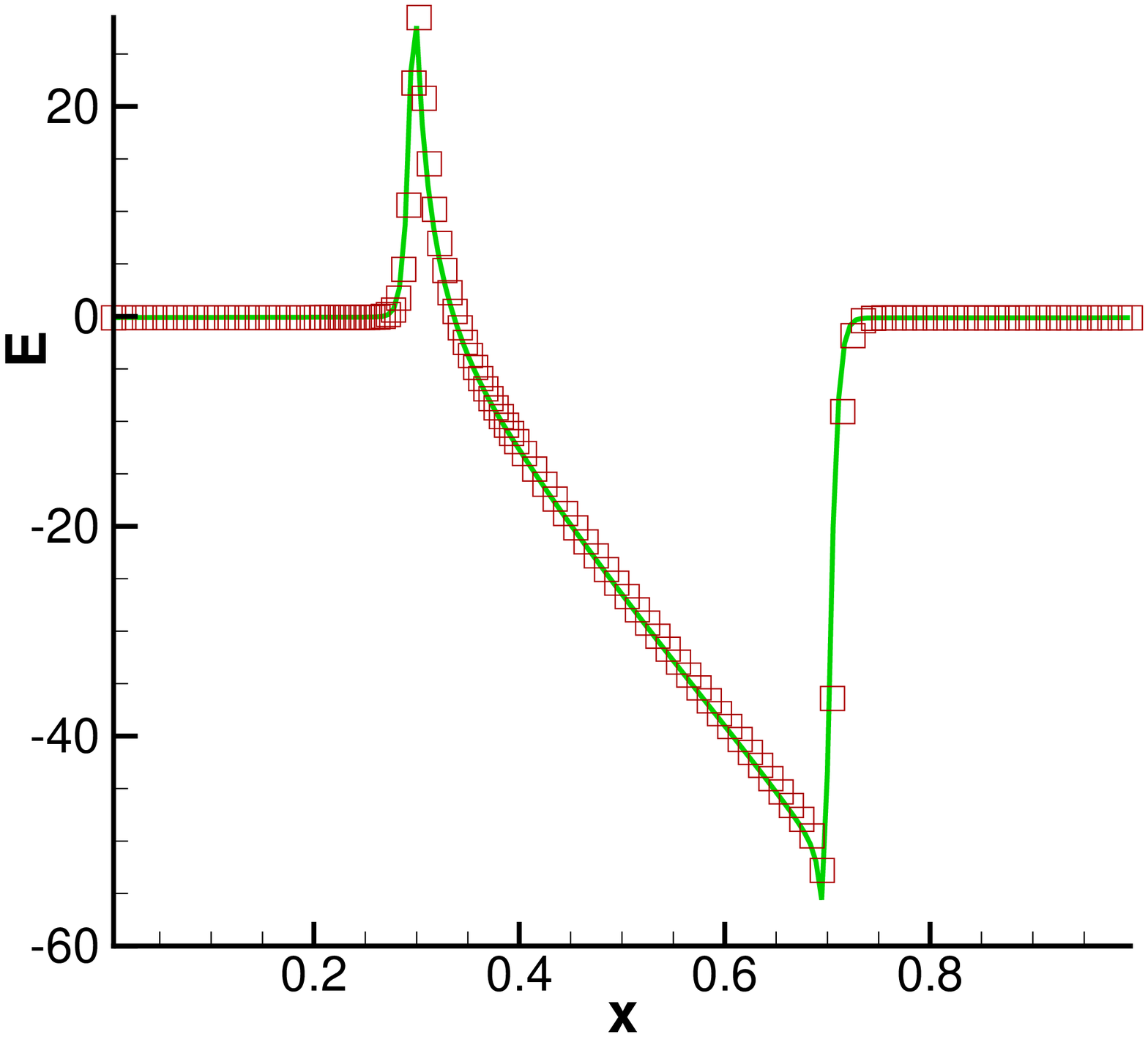}\\
\includegraphics[width=2.93in,angle=0]{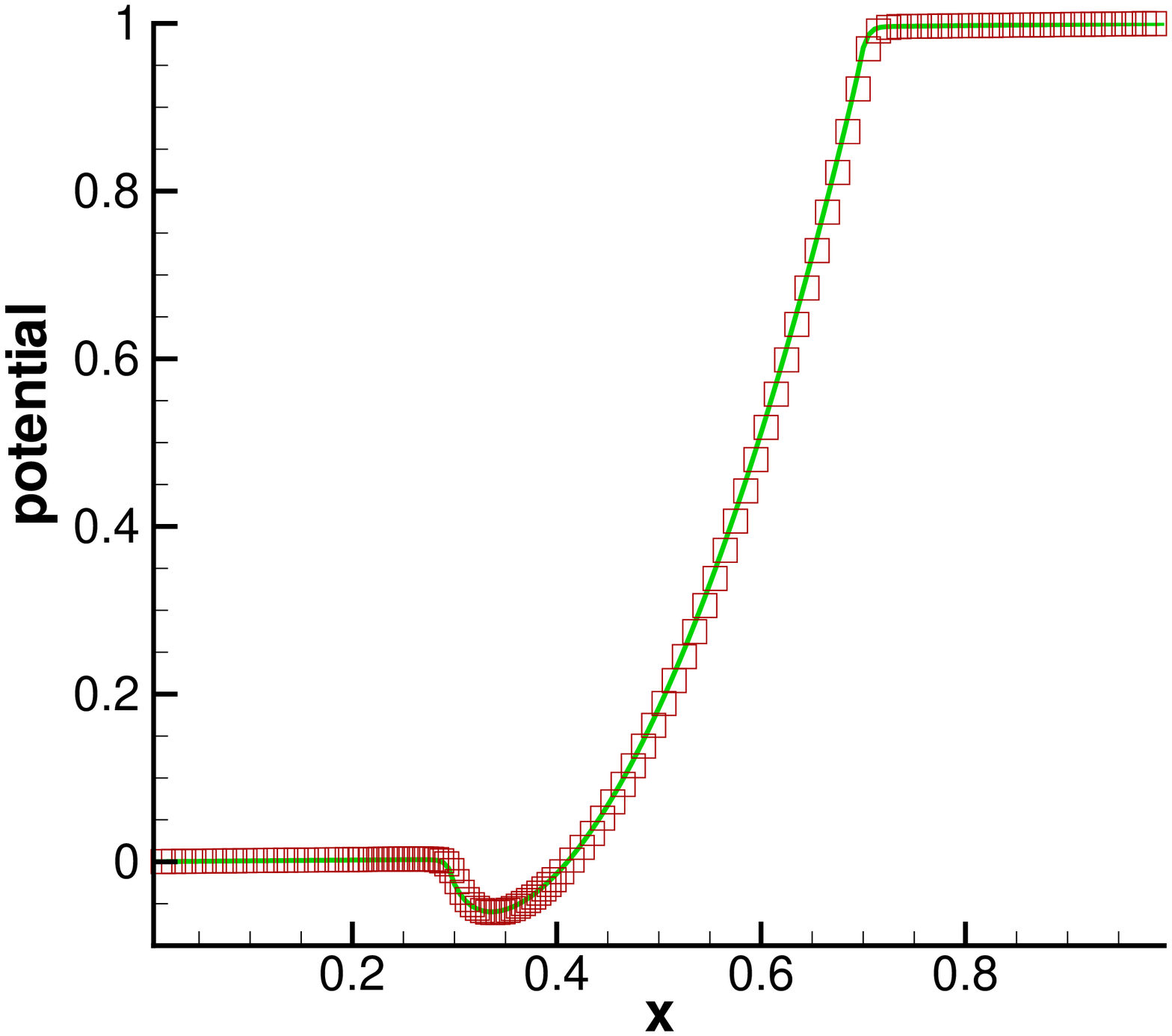}
\includegraphics[width=2.93in,angle=0]{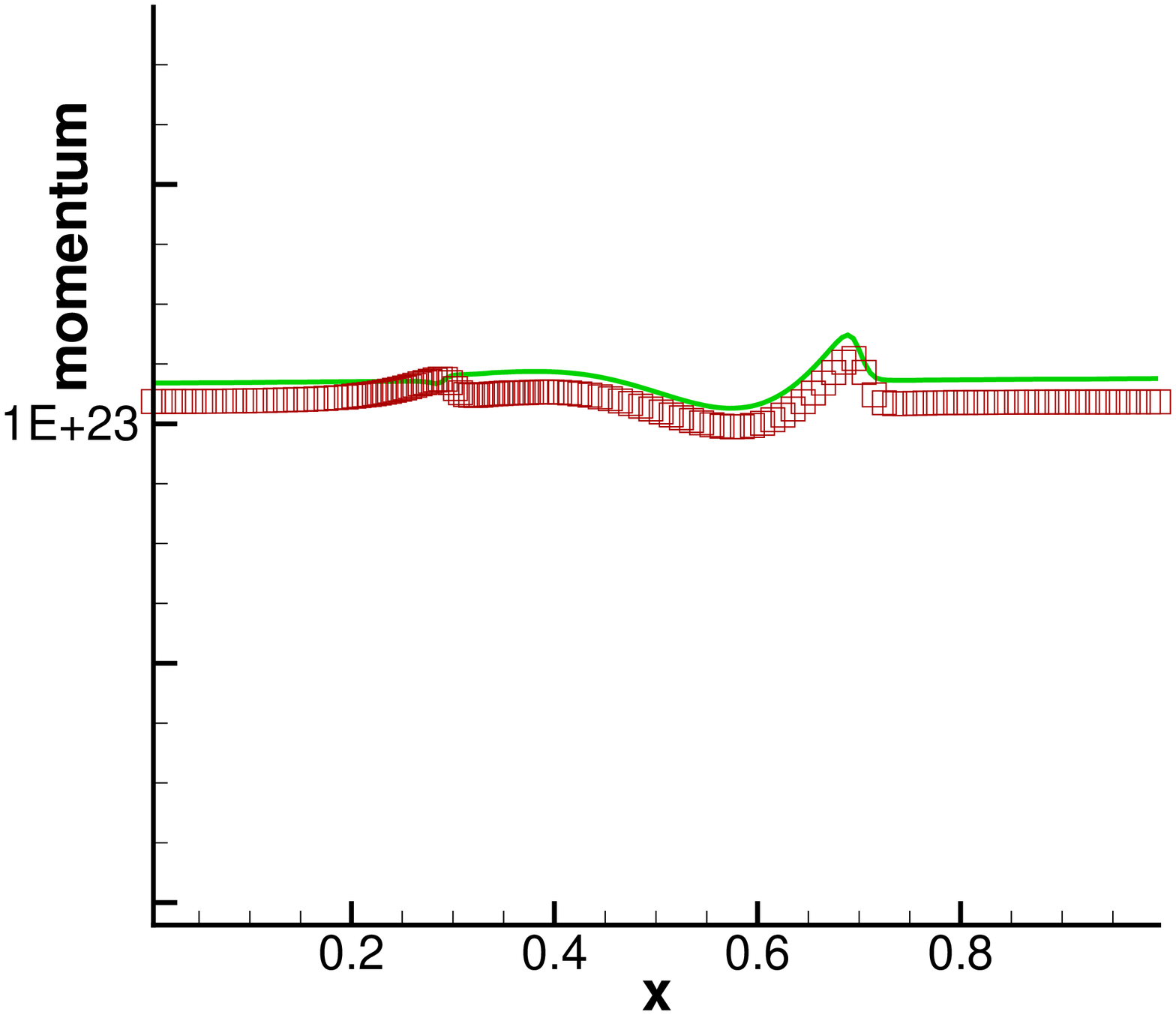}
\caption{Comparison of macroscopic quantities using DG (symbols) and
WENO (solid line) for $400$nm channel at $t=5.0$,
$V_{\mbox{bias}}=1.0$. Top left: density in ${cm}^{-3}$; top right:
mean velocity in $cm/s$; middle left: energy in $eV$; middle right:
electric field in $kV/cm$; bottom left: potential in $V$; bottom
right: momentum in ${cm}^{-2} \, s^{-1}$. Solution has reached
steady state.} \label{400m2}
\end{figure}

\begin{figure}[htb]
\centering
\includegraphics[width=3in,angle=0]{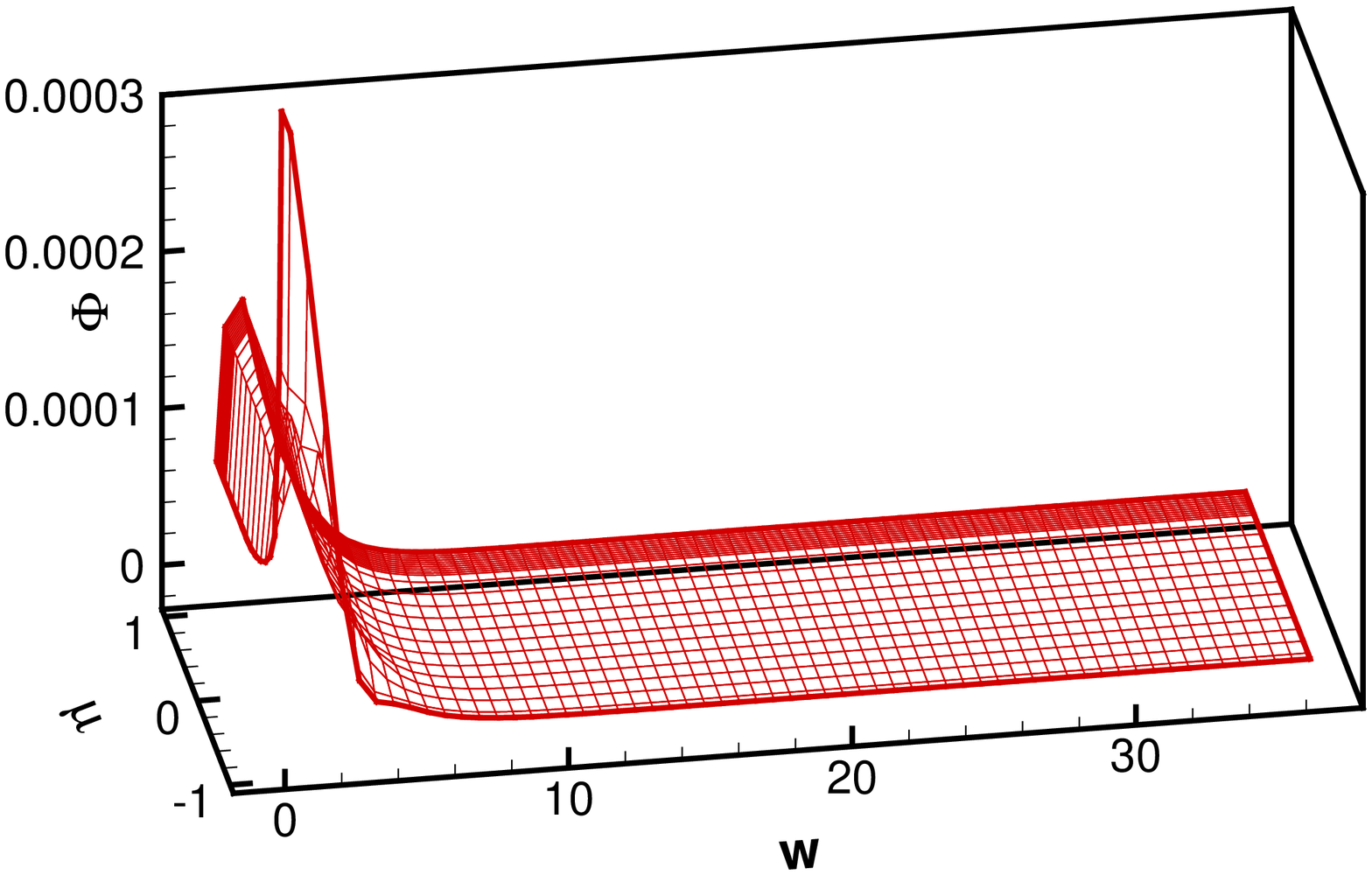}
\includegraphics[width=3in,angle=0]{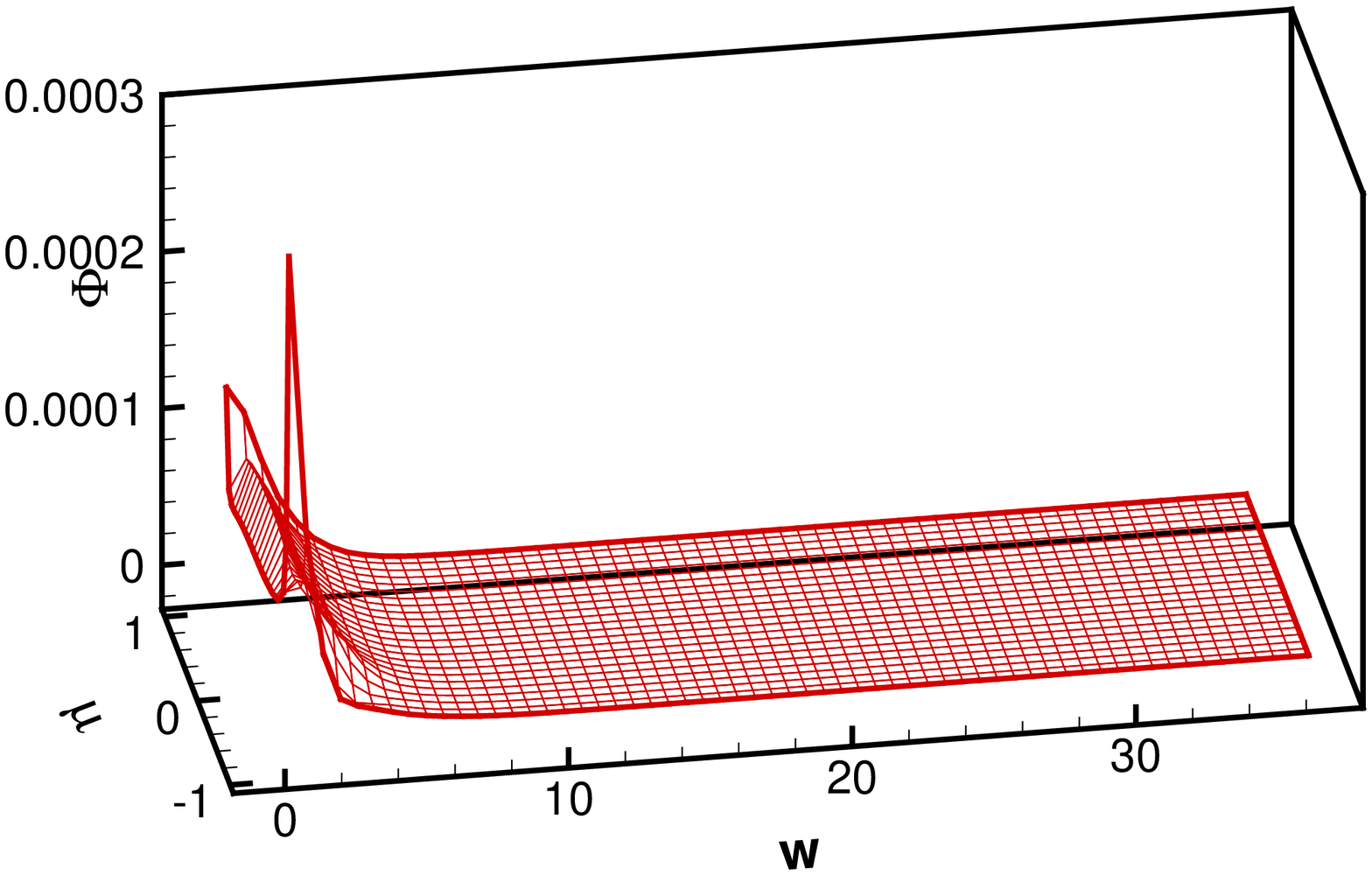}\\
\includegraphics[width=3in,angle=0]{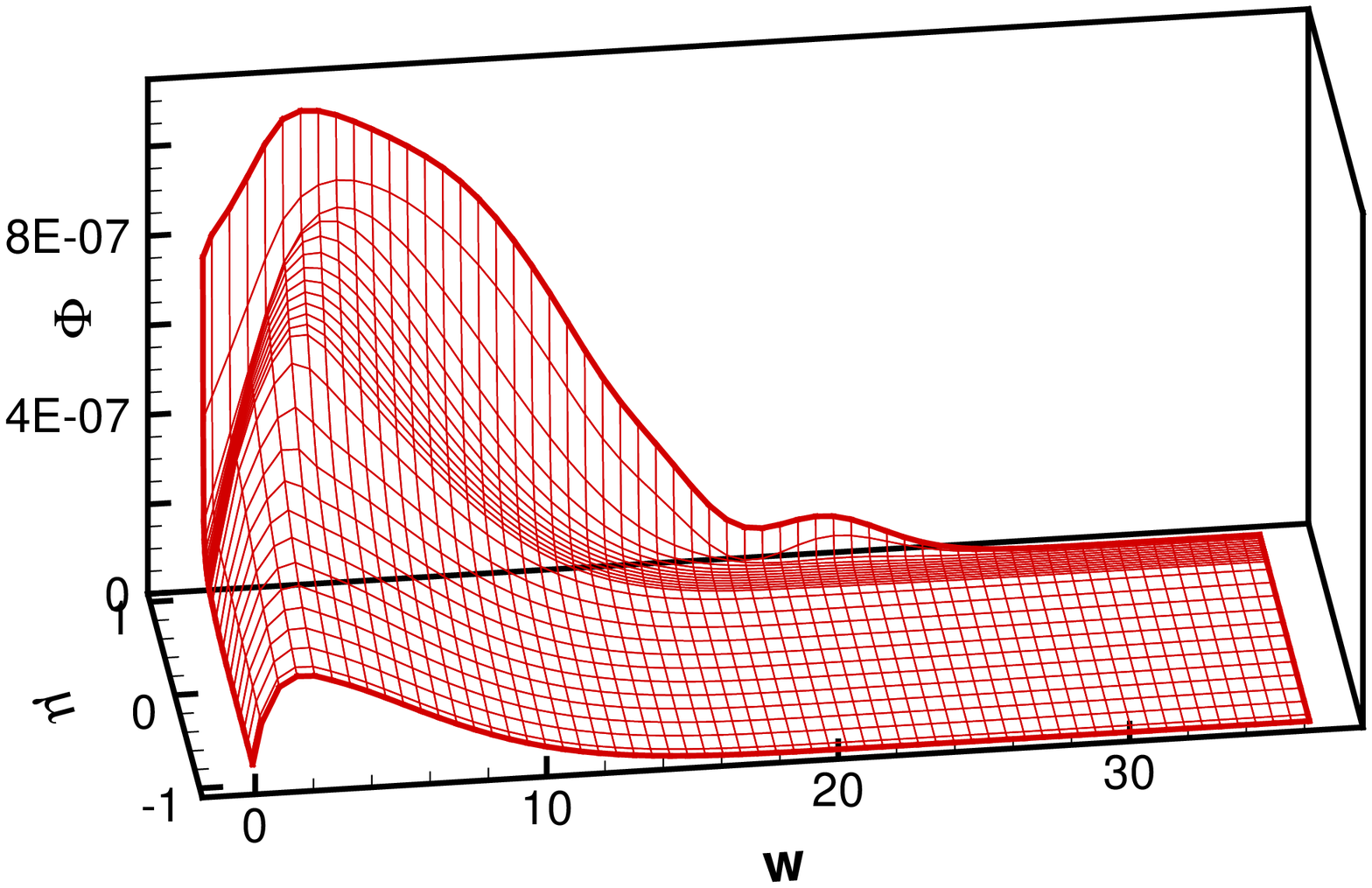}
\includegraphics[width=3in,angle=0]{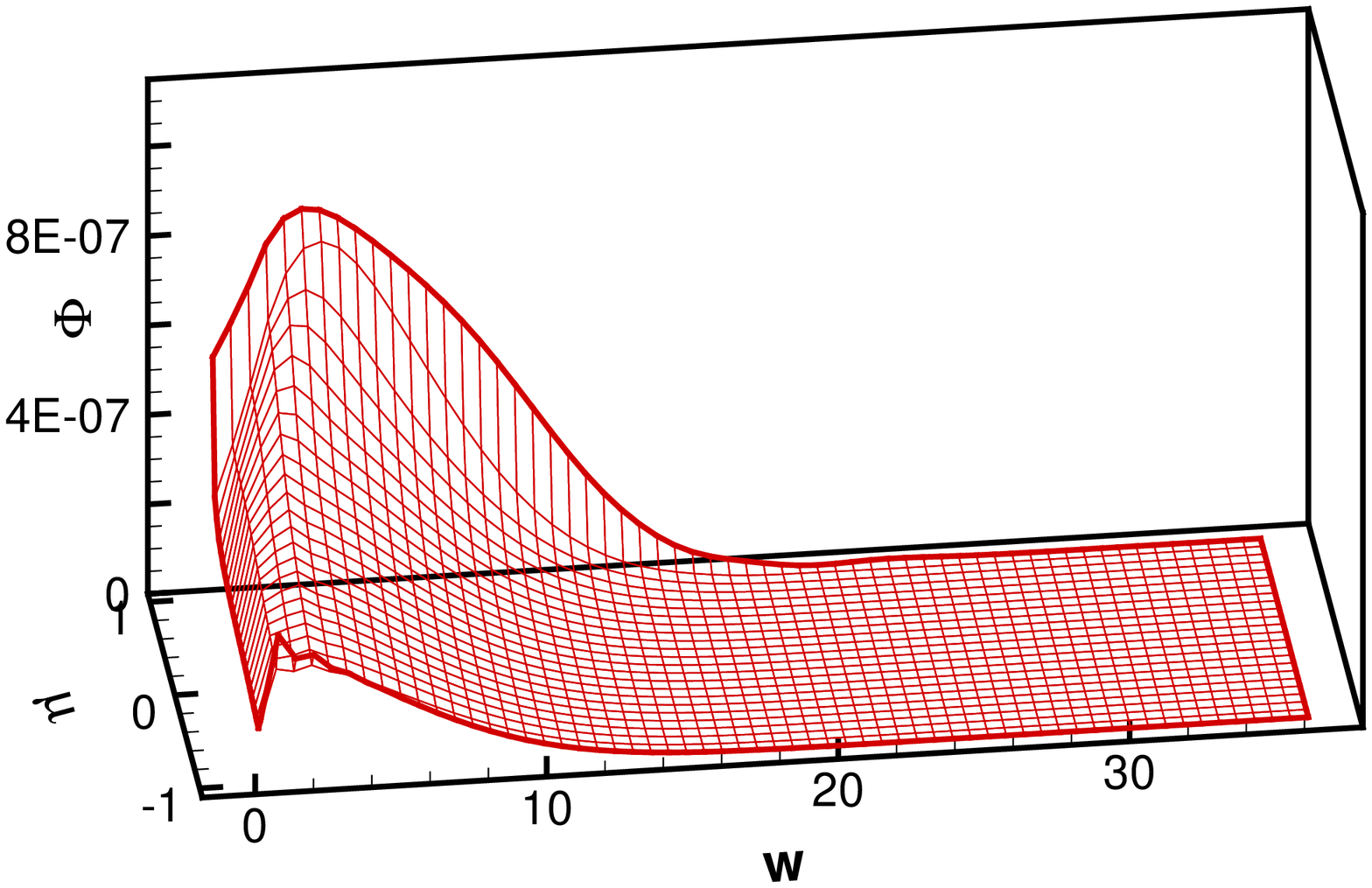}\\
\includegraphics[width=3in,angle=0]{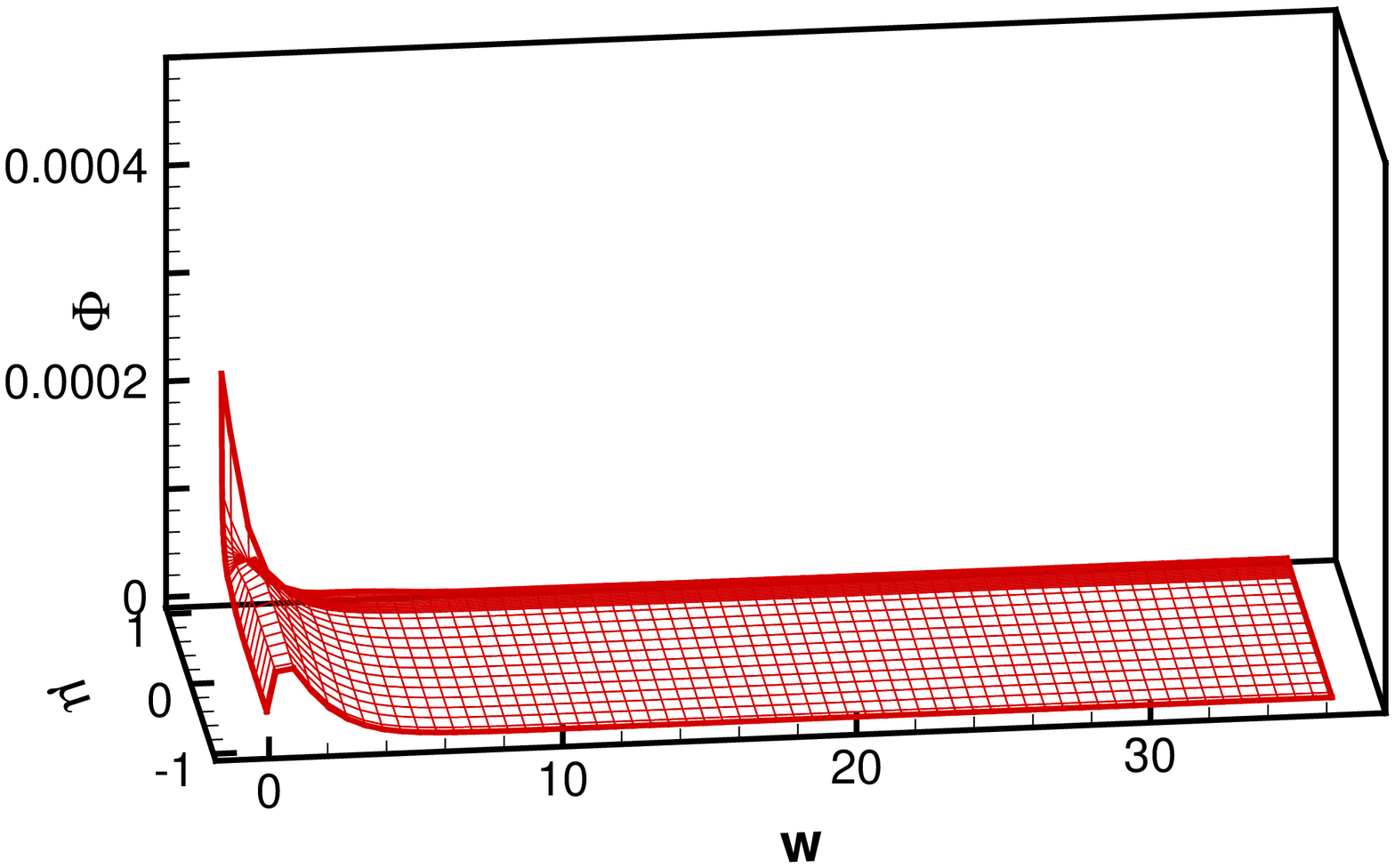}
\includegraphics[width=3in,angle=0]{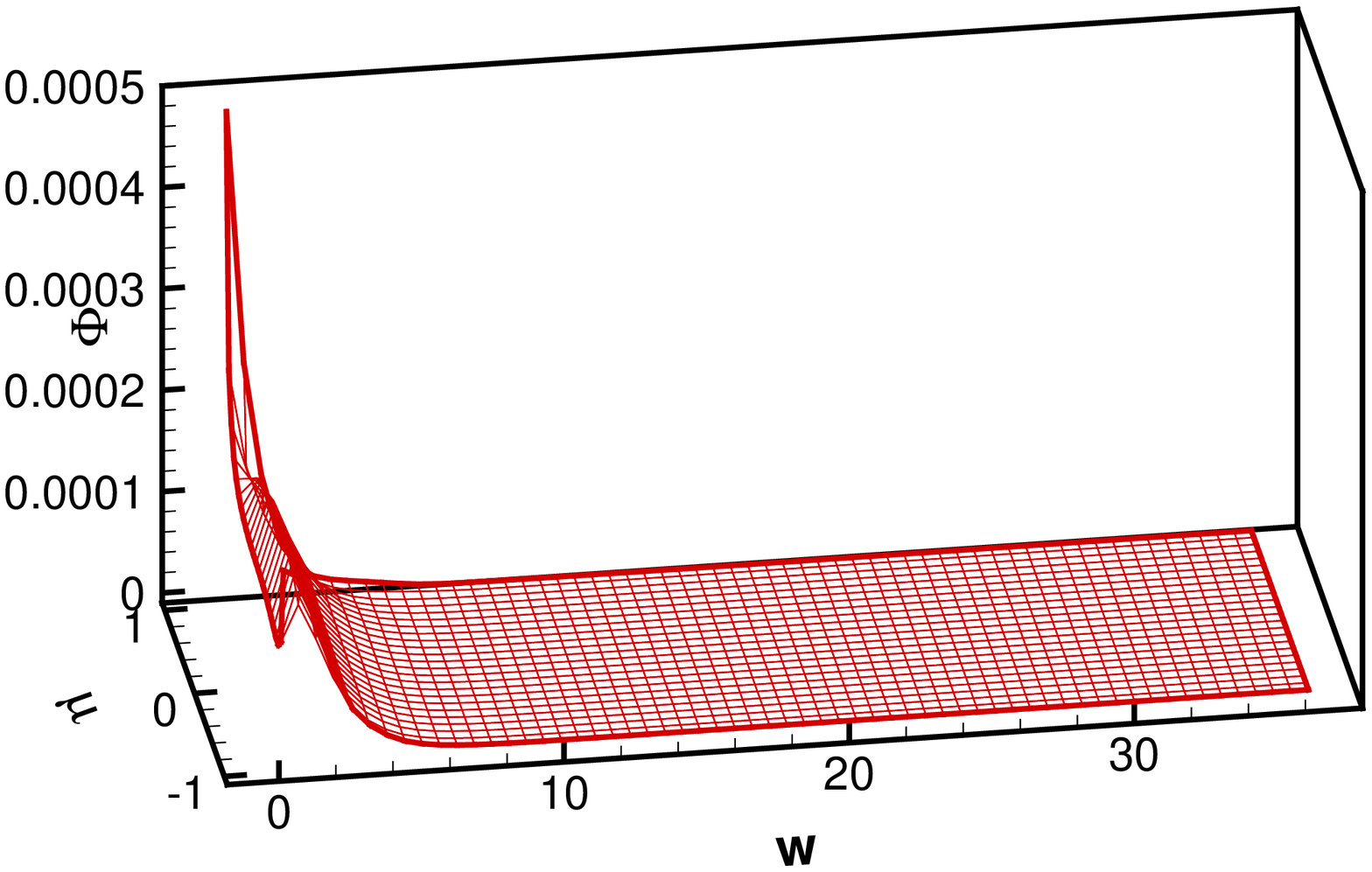}\\
\caption{Comparison of the snapshot for $\Phi(x_0,w,\mu)$ using DG
(left) and WENO (right) for $400$nm channel at $t=0.5$,
$V_{\mbox{bias}}=1.0$. Top:  $x_0=0.3$; middle:  $x_0=0.5$; bottom:
$x_0=0.7$. Solution has not yet reached steady state.} \label{400p1}
\end{figure}

\begin{figure}[htb]
\centering
\includegraphics[width=3in,angle=0]{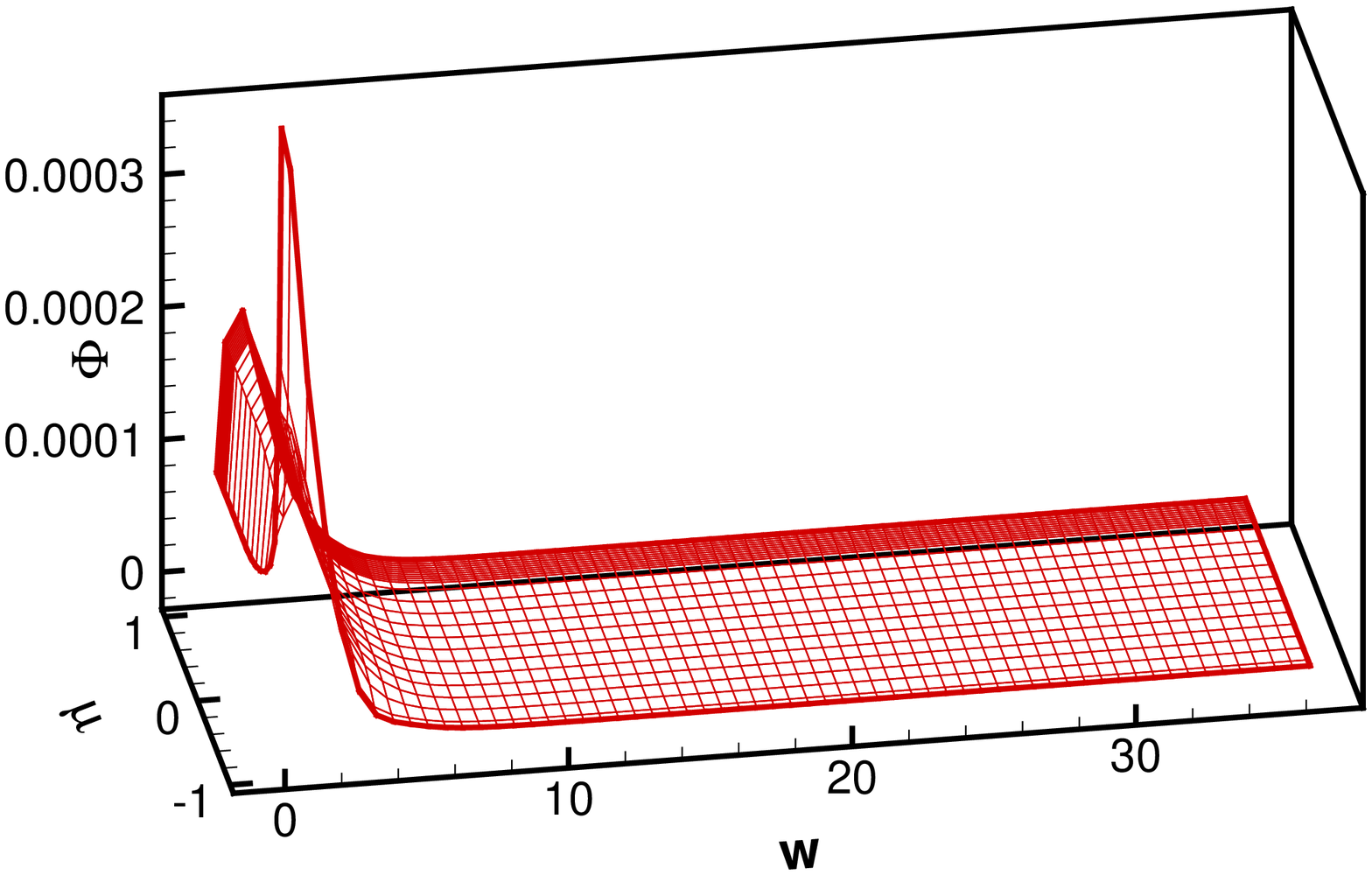}
\includegraphics[width=3in,angle=0]{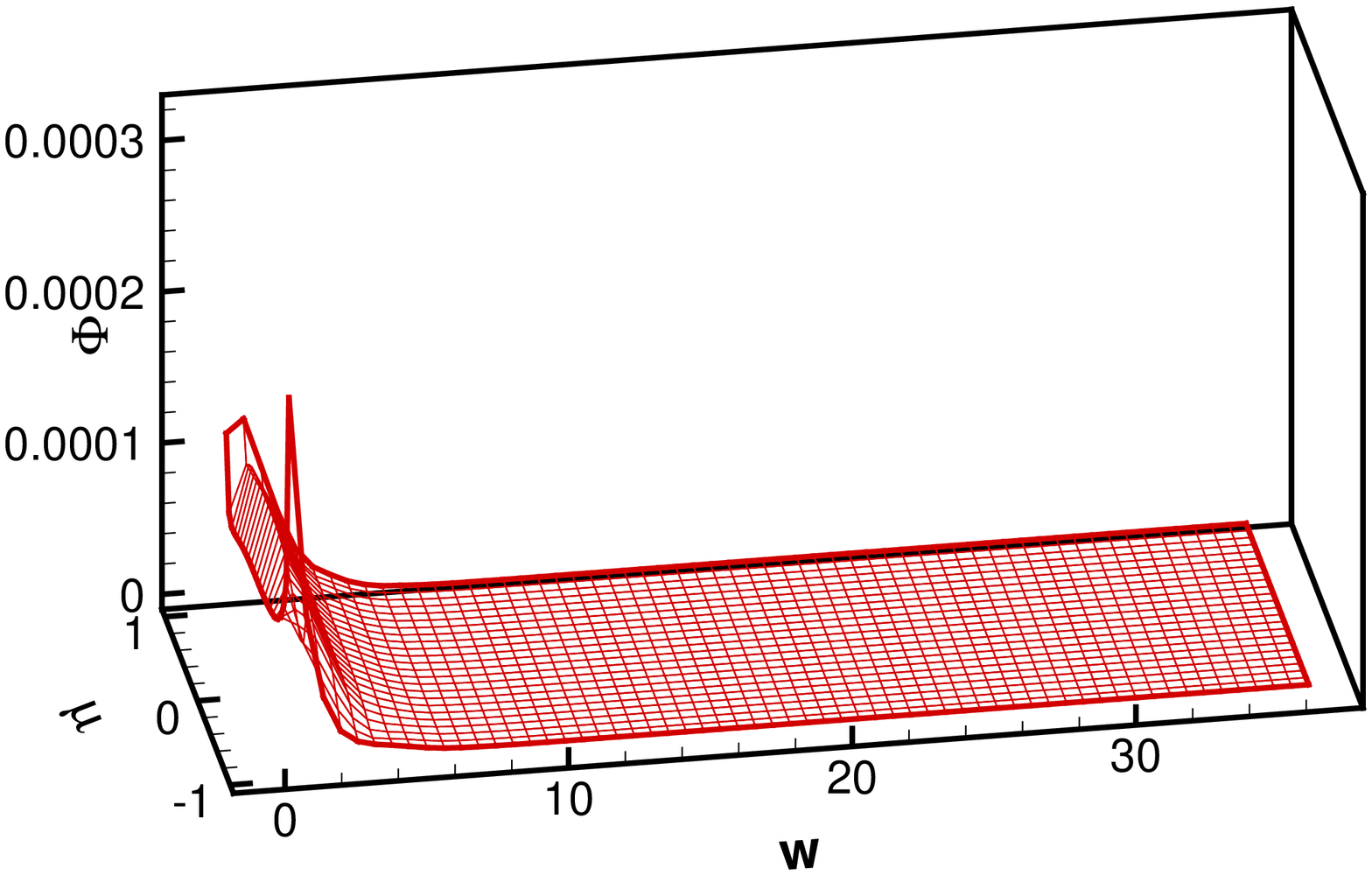}\\
\includegraphics[width=3in,angle=0]{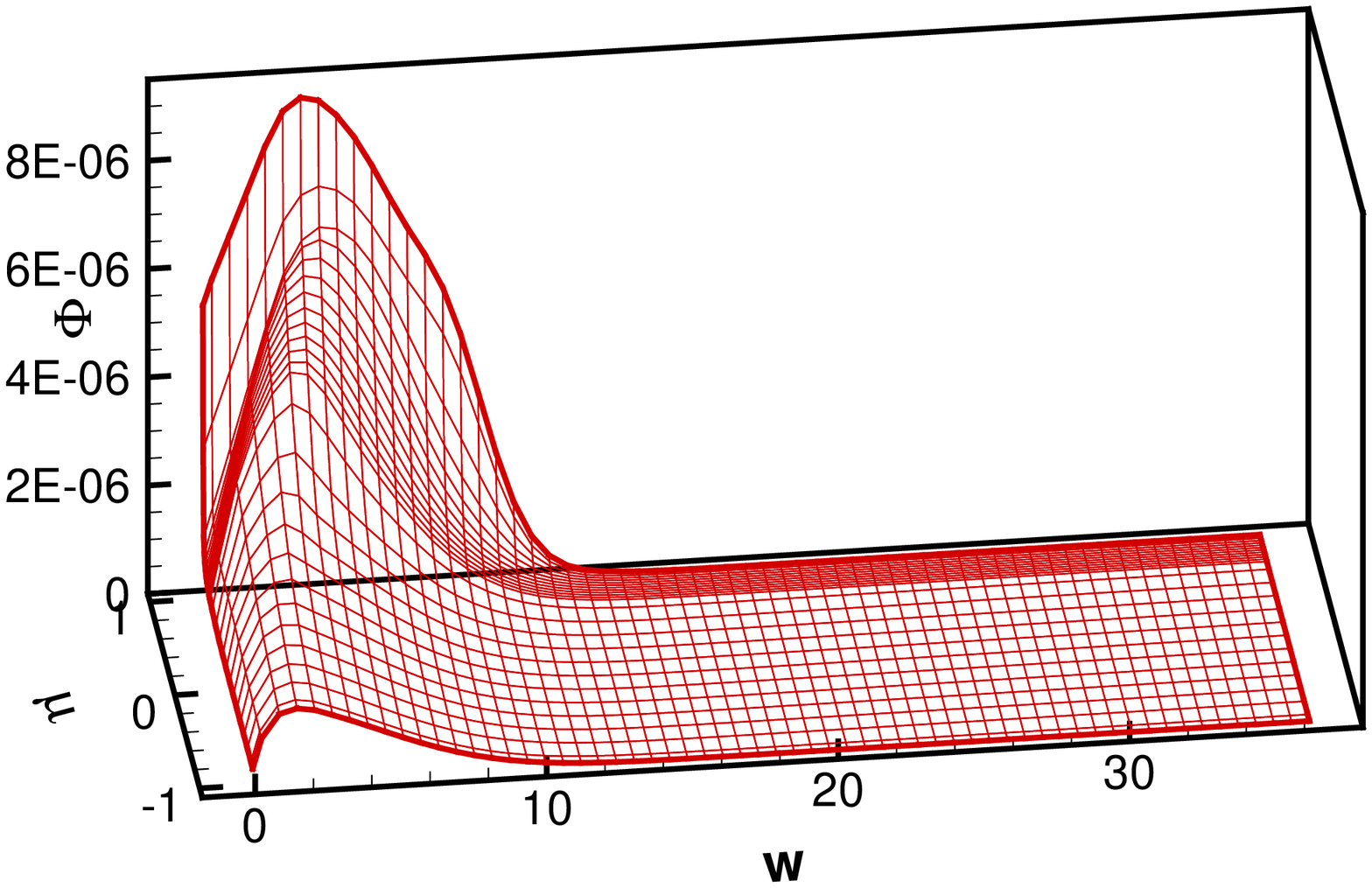}
\includegraphics[width=3in,angle=0]{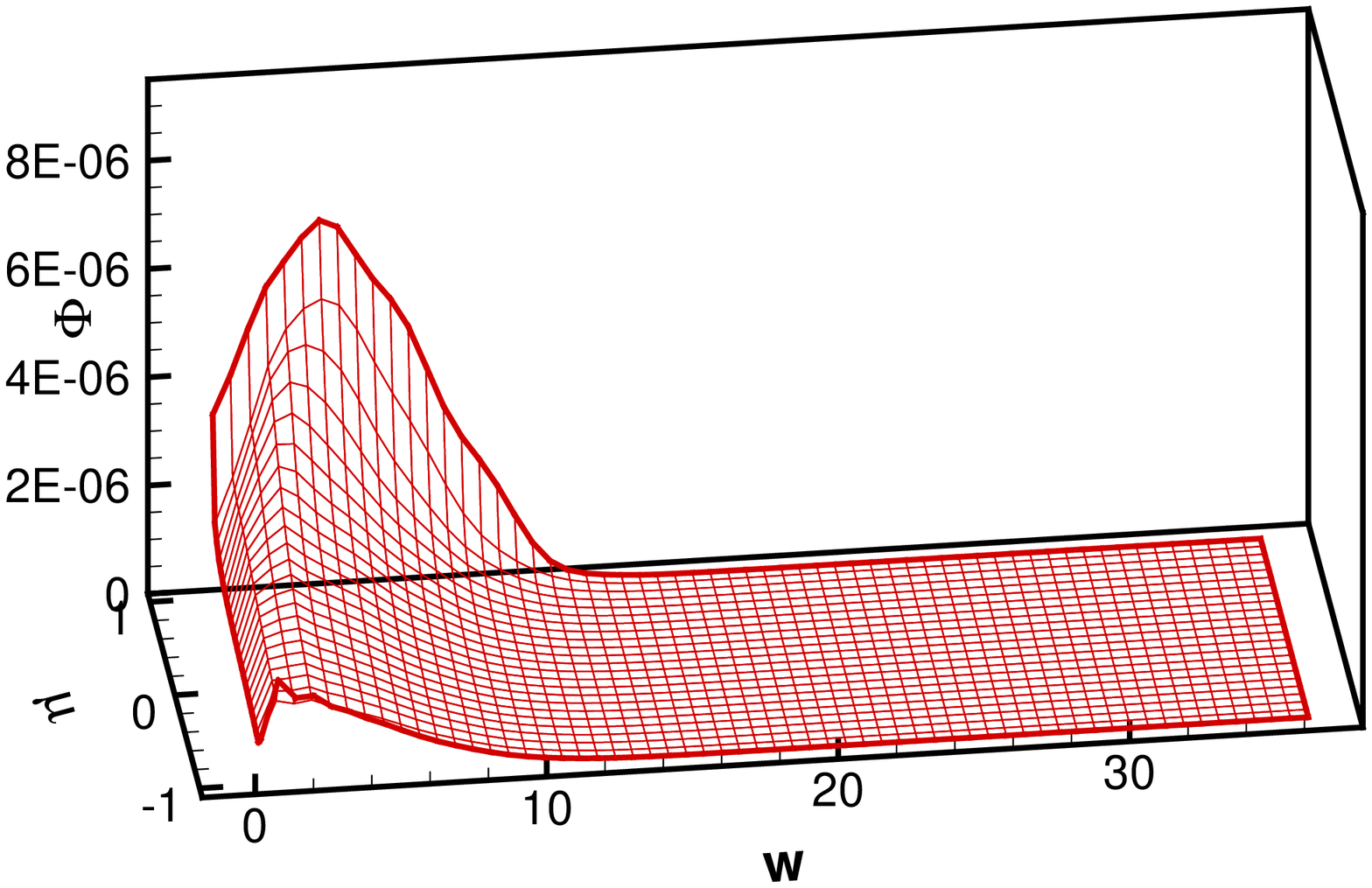}\\
\includegraphics[width=3in,angle=0]{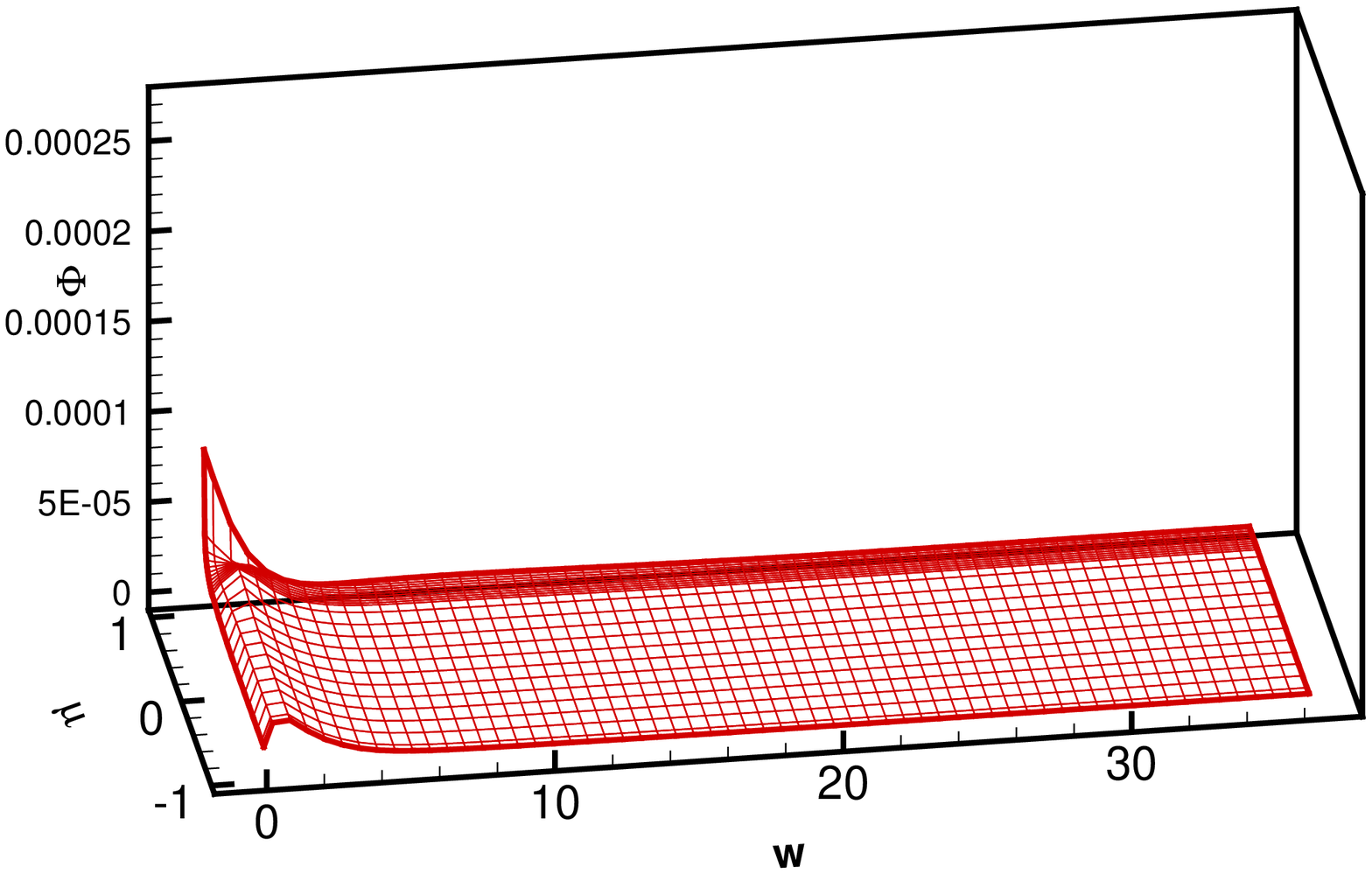}
\includegraphics[width=3in,angle=0]{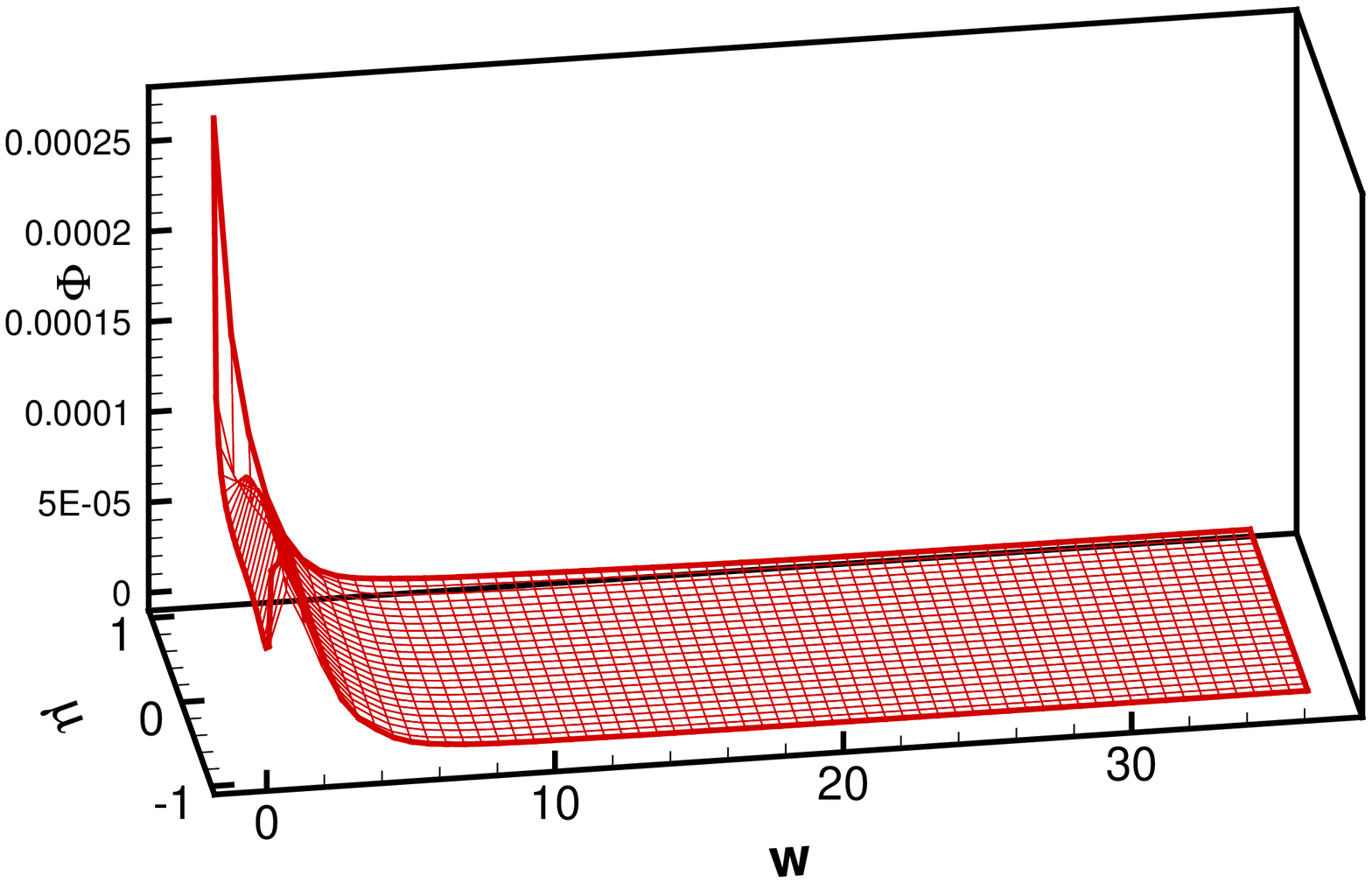}\\
\caption{Comparison of the snapshot for $\Phi(x_0,w,\mu)$ using DG
(left) and WENO (right) for $400$nm channel at $t=5.0$,
$V_{\mbox{bias}}=1.0$. Top:  $x_0=0.3$; middle:  $x_0=0.5$; bottom:
$x_0=0.7$.  Solution has reached steady state.} \label{400p2}
\end{figure}

\begin{figure}[htb]
\centering
\includegraphics[width=4in,angle=0]{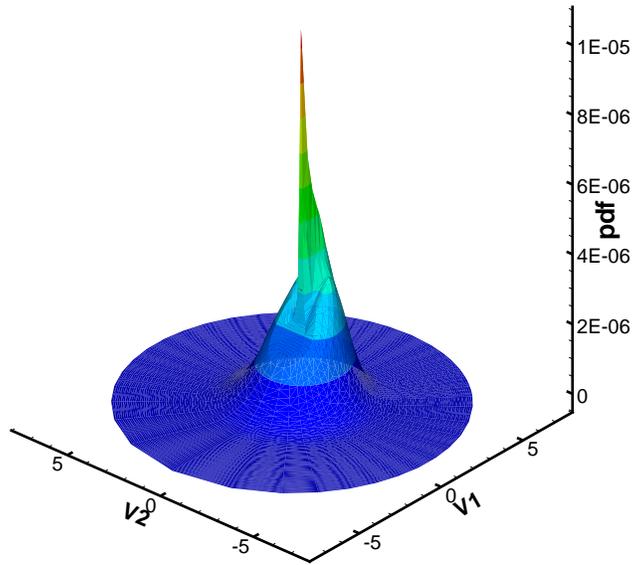}
\caption{PDF for 400nm channel at $t=5.0$, $x=0.5$.} \label{400cart}
\end{figure}

We plot the evolution of density, mean velocity, energy and momentum
in Figure \ref{400evo}. The solution has already stabilized at
$t=5.0$ from the momentum plots. The macroscopic quantities at
steady state are plotted in Figure  \ref{400m2}. The results are
compared with the WENO calculation. They agree with each other in
general, with DG offering more resolution and a higher peak in
energy near the junctions. Figures \ref{400p1} and \ref{400p2} show
comparisons for the \emph{pdf} at transient and steady state. We
plot at different position of the device, namely, the left, center
and right of the channel. We notice a  larger value of \emph{pdf}
especially at the  center  of the channel, where the \emph{pdf} is
no longer  Maxwellian. Moreover, at $t=0.5$, $x_0=0.5$, the
\emph{pdf} shows a double hump structure, which is not captured by
the WENO solver. All of these advantages come from the fact that we
are refining more near $\mu=1$. To have a better idea of the shape
of the \emph{pdf}, we  plot $\Phi(t=5.0, x=0.5)$ in the cartesian
coordinates in Figure \ref{400cart}. The coordinate $V1$ in the plot
is the momentum parallel to the force field  $k_1$, $V2$ is the
modulus of the orthogonal component. The peak is captured very
sharply compared to WENO.

\begin{figure}[htb]
\centering
\includegraphics[width=3in,angle=0]{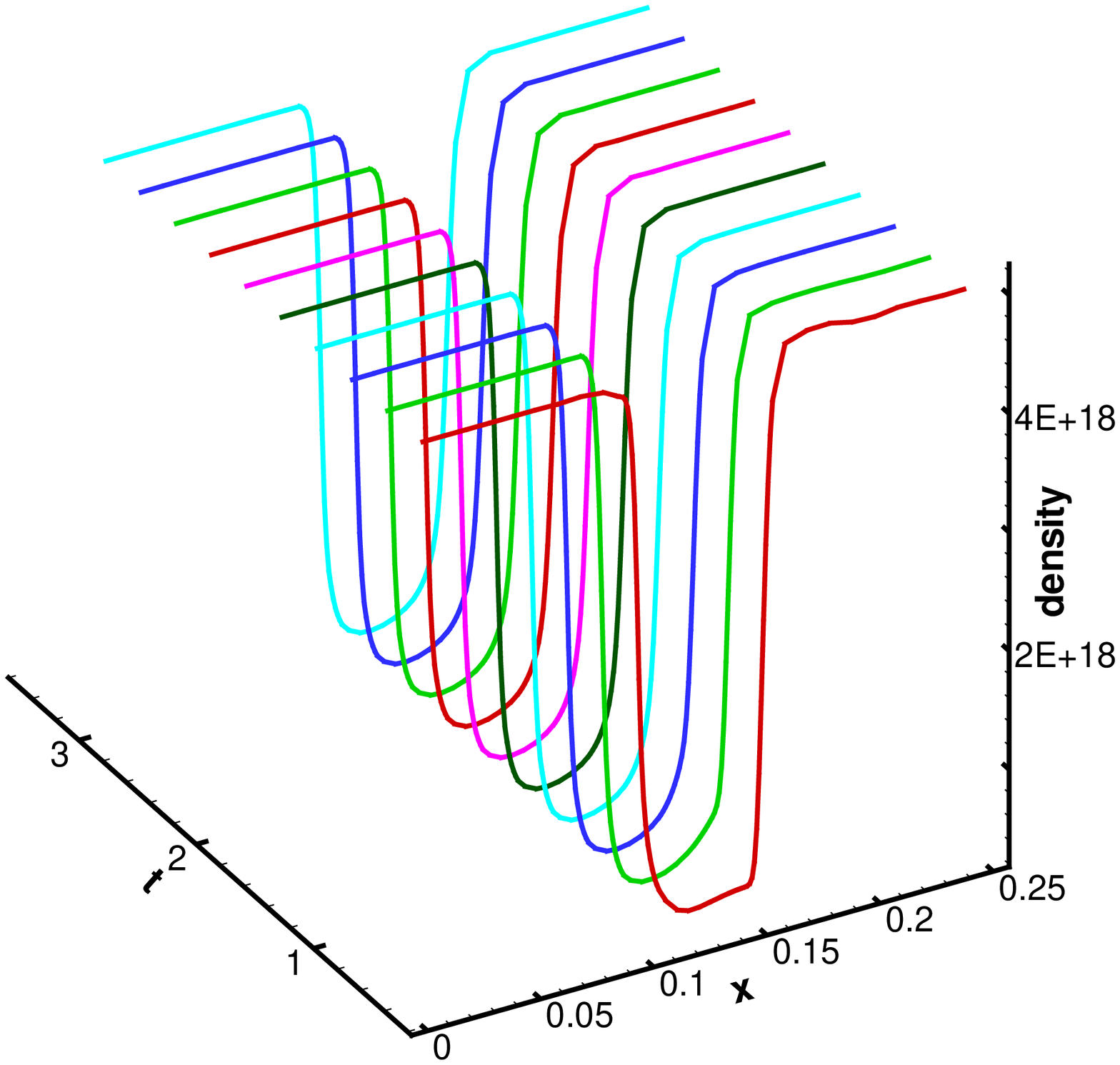}
\includegraphics[width=3in,angle=0]{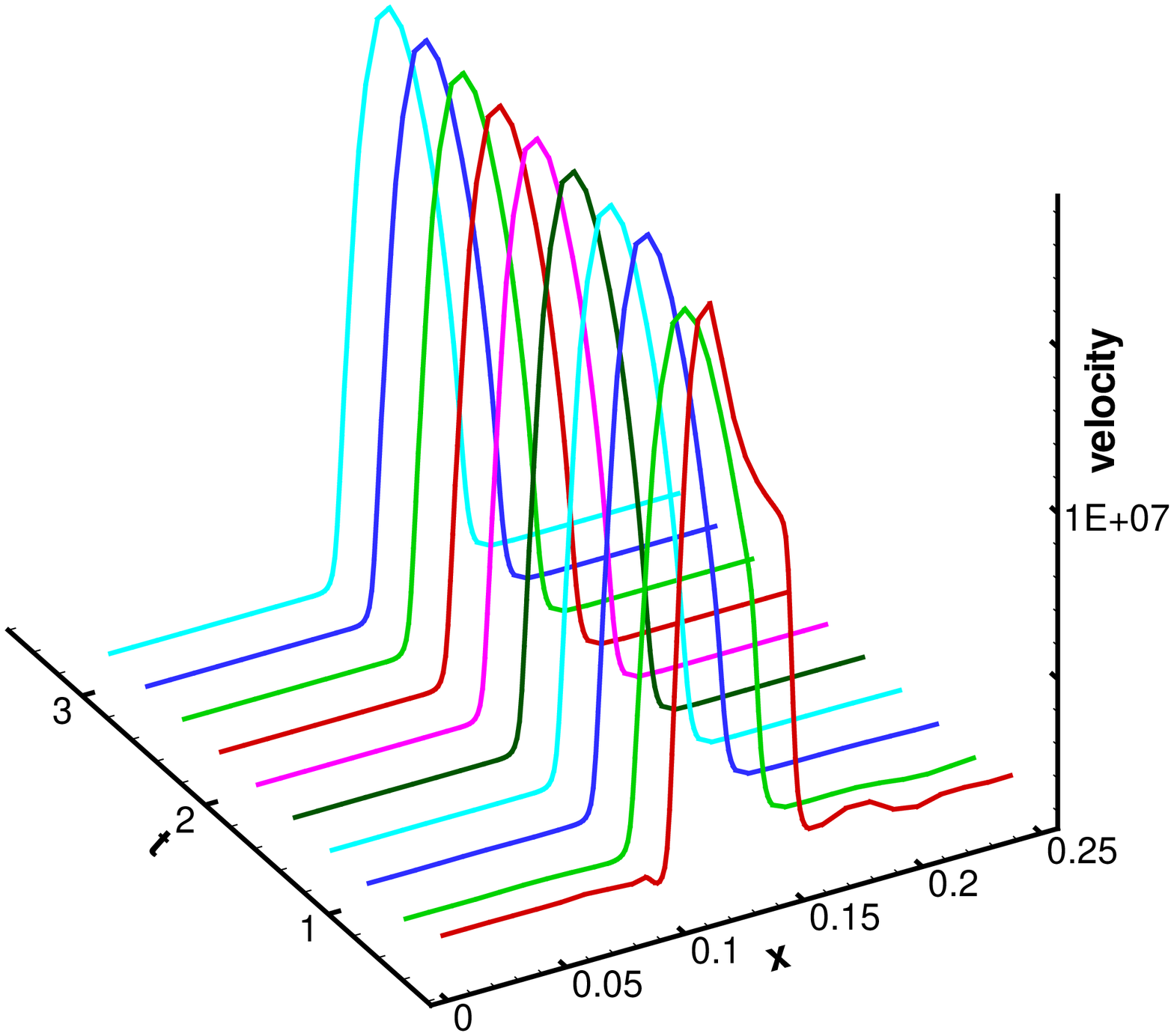}\\
\includegraphics[width=3in,angle=0]{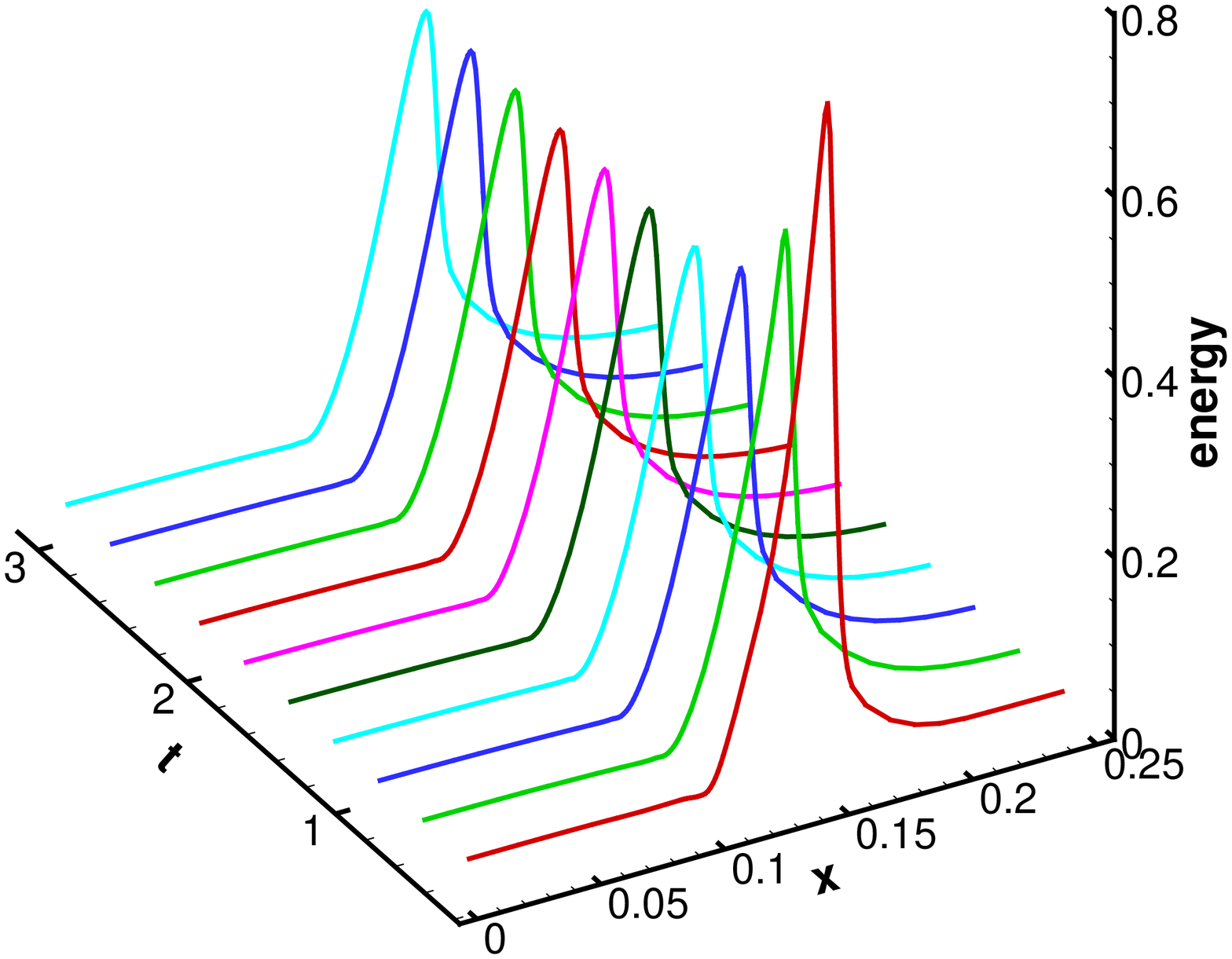}
\includegraphics[width=3in,angle=0]{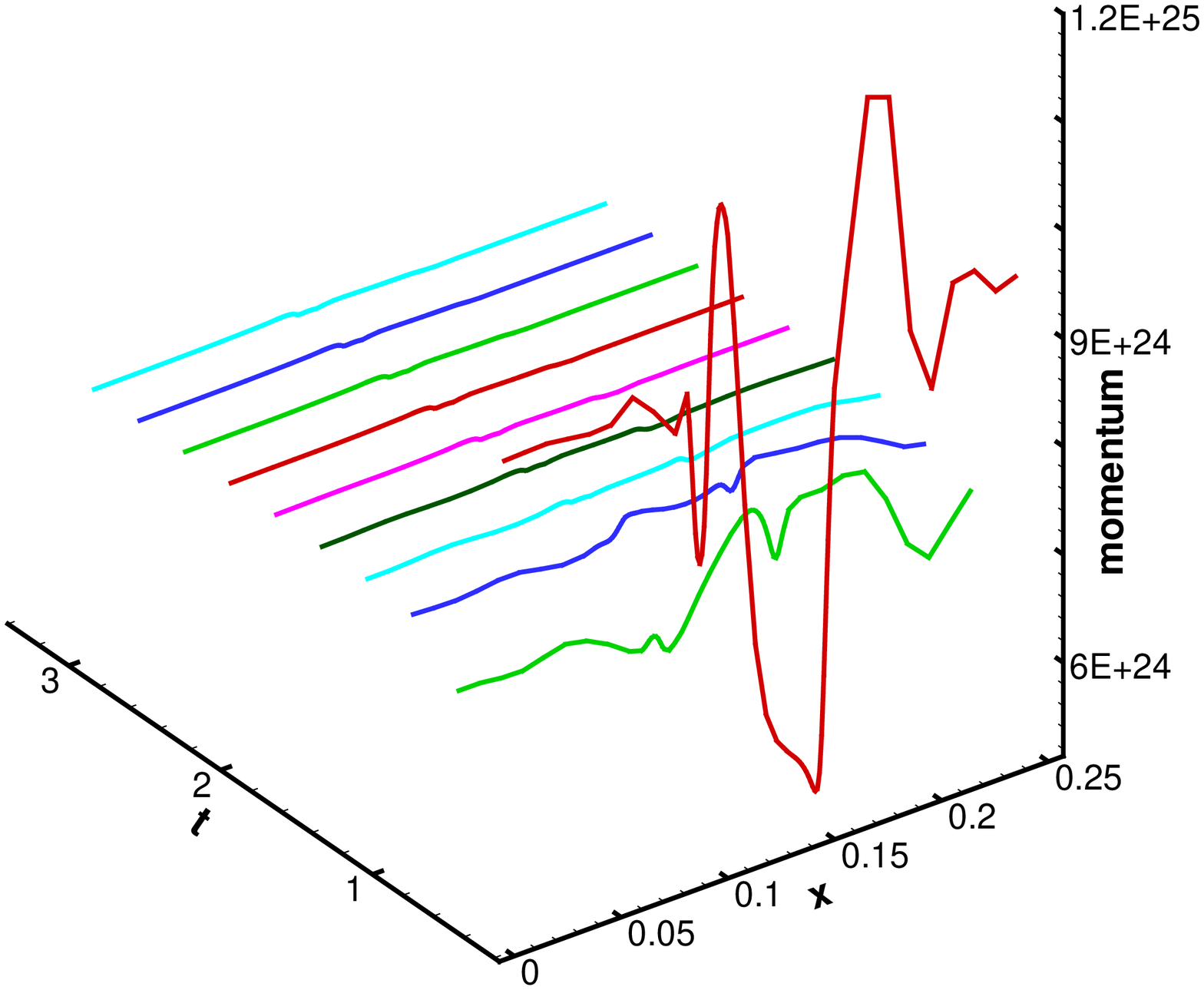}
\caption{Time evolution of macroscopic quantities using DG method
for $50$nm channel at  $V_{\mbox{bias}}=1.0$. Top left: density in
${cm}^{-3}$; top right: mean velocity in $cm/s$; bottom left: energy
in $eV$;  bottom right: momentum in ${cm}^{-2} \, s^{-1}$. }
\label{50evo}
\end{figure}

\begin{figure}[htb]
\centering
\includegraphics[width=2.93in,angle=0]{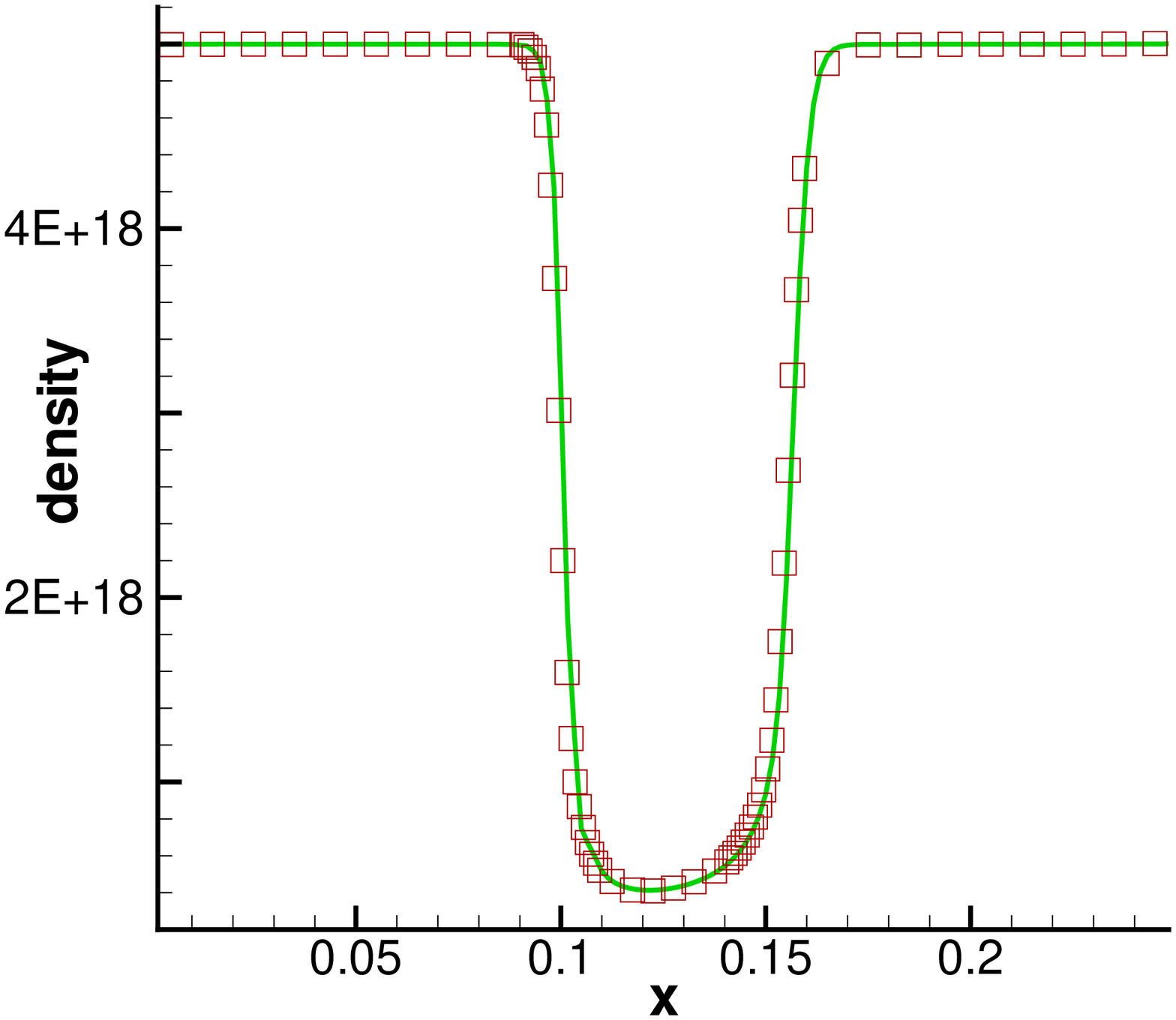}
\includegraphics[width=2.93in,angle=0]{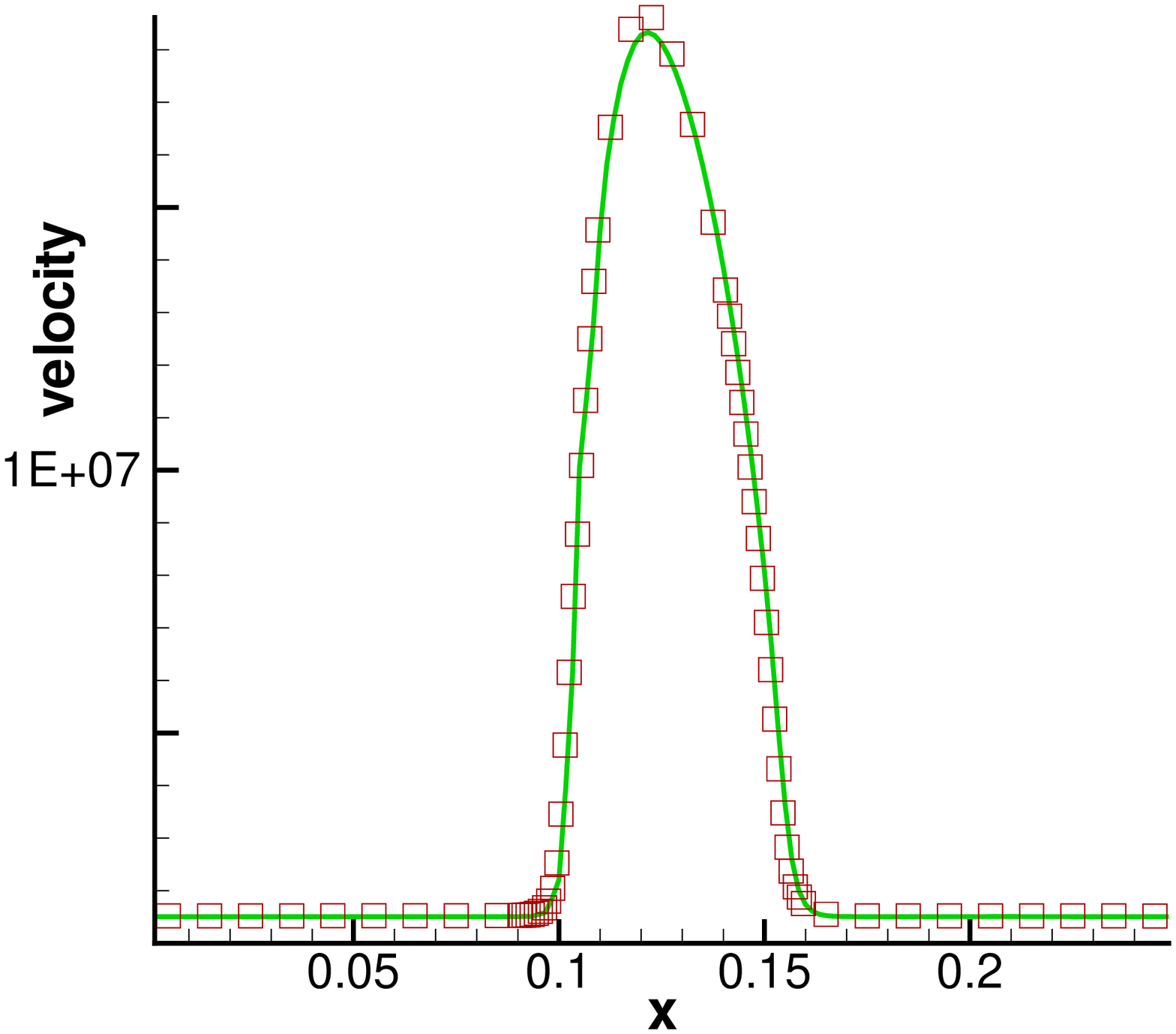}\\
\includegraphics[width=2.93in,angle=0]{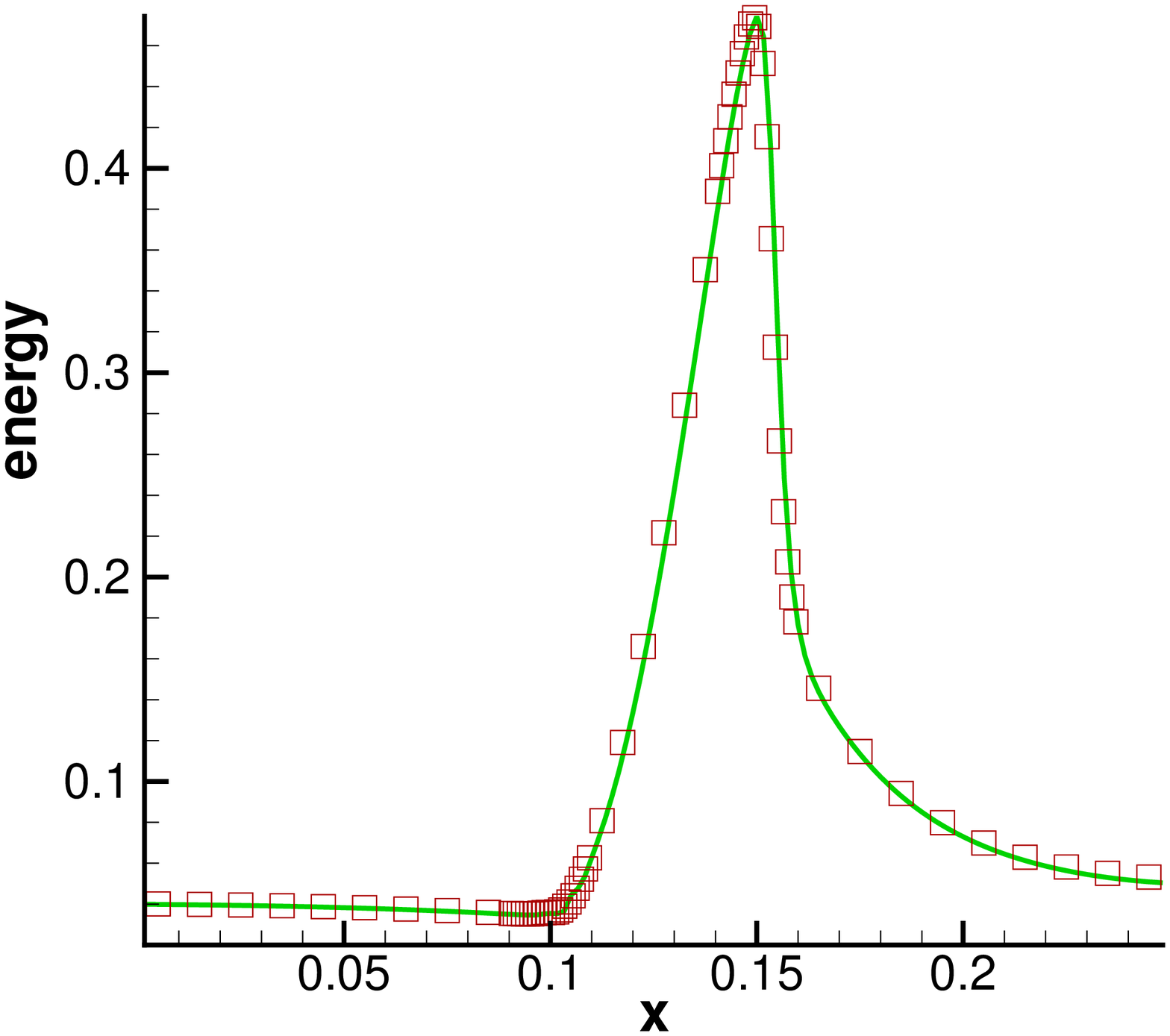}
\includegraphics[width=2.93in,angle=0]{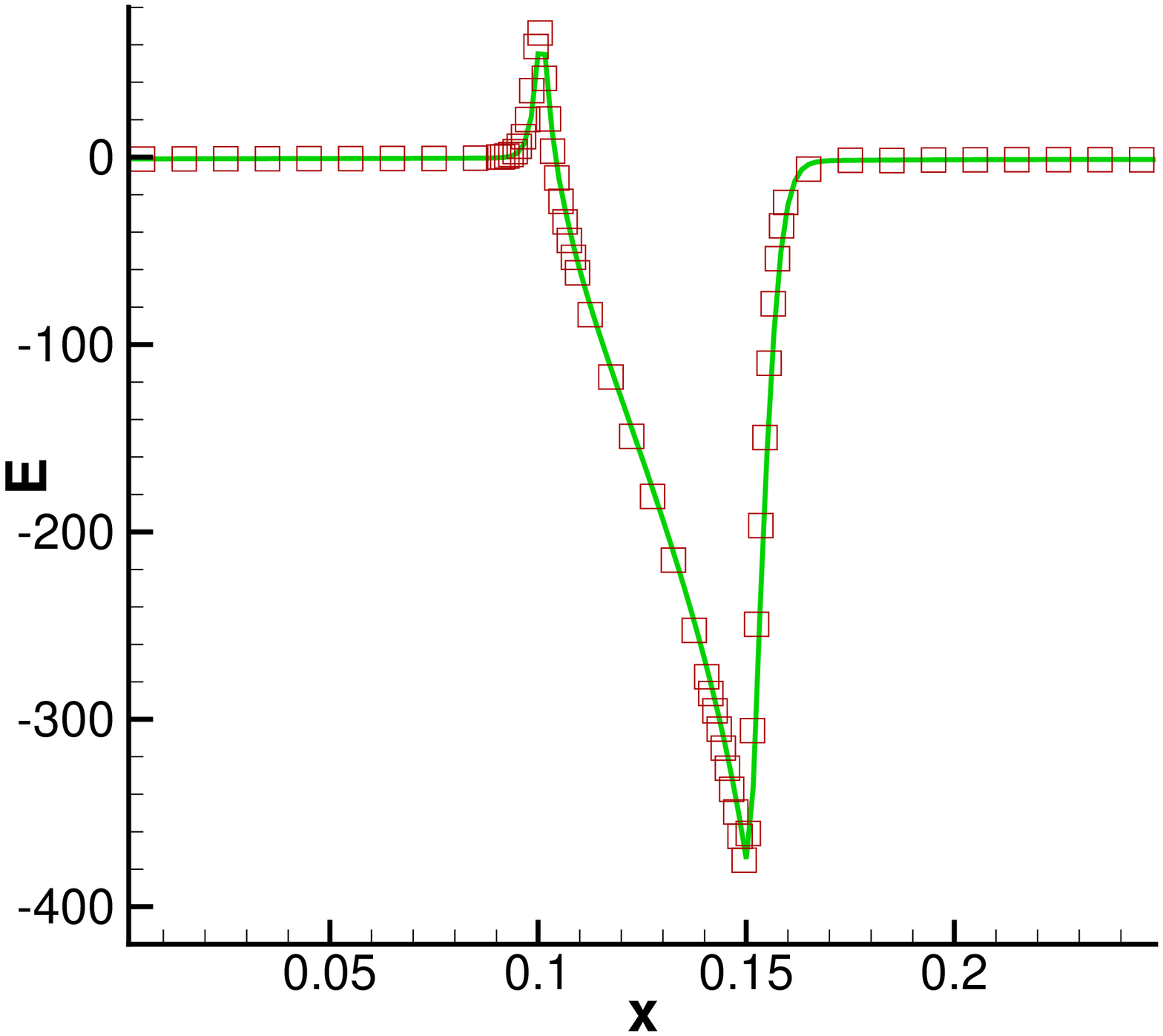}\\
\includegraphics[width=2.93in,angle=0]{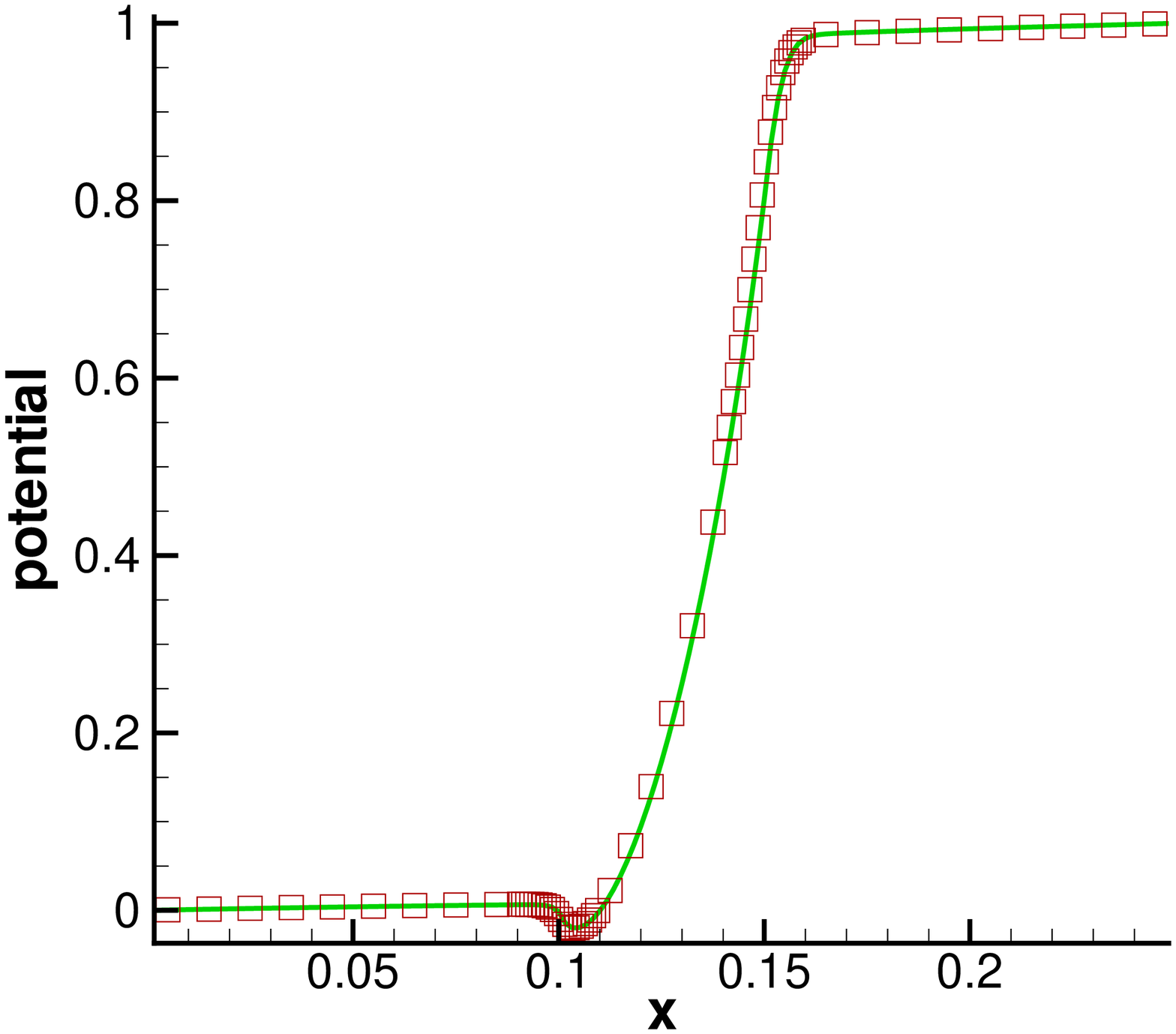}
\includegraphics[width=2.93in,angle=0]{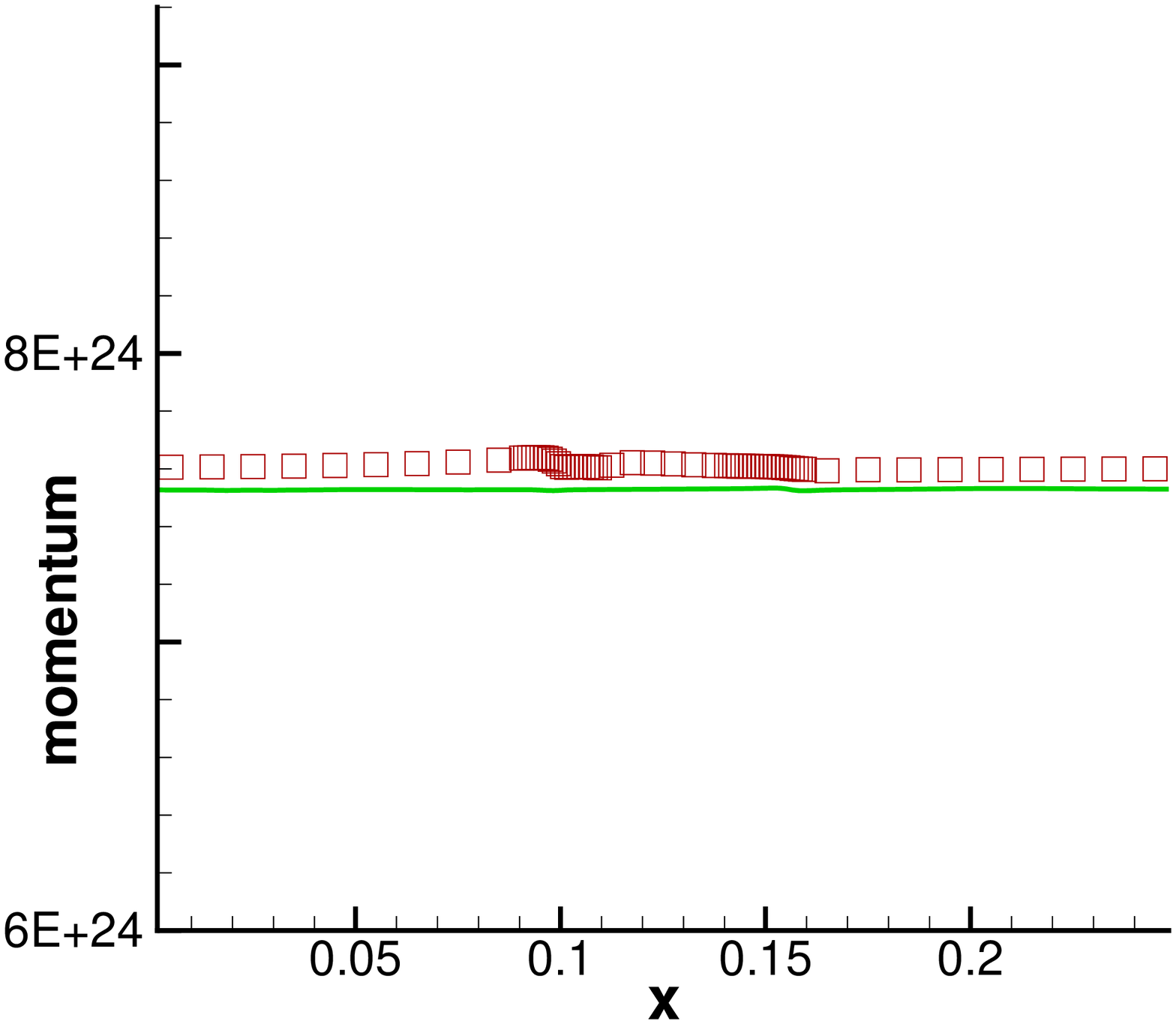}
\caption{Comparison of macroscopic quantities using DG (symbols) and
WENO (solid line) for $50$nm channel at $t=3.0$,
$V_{\mbox{bias}}=1.0$. Top left: density in ${cm}^{-3}$; top right:
mean velocity in $cm/s$; middle left: energy in $eV$; middle right:
electric field in $kV/cm$; bottom left: potential in $V$; bottom
right: momentum in ${cm}^{-2} \, s^{-1}$. Solution has reached
steady state.} \label{50m2}
\end{figure}

\begin{figure}[htb]
\centering
\includegraphics[width=3in,angle=0]{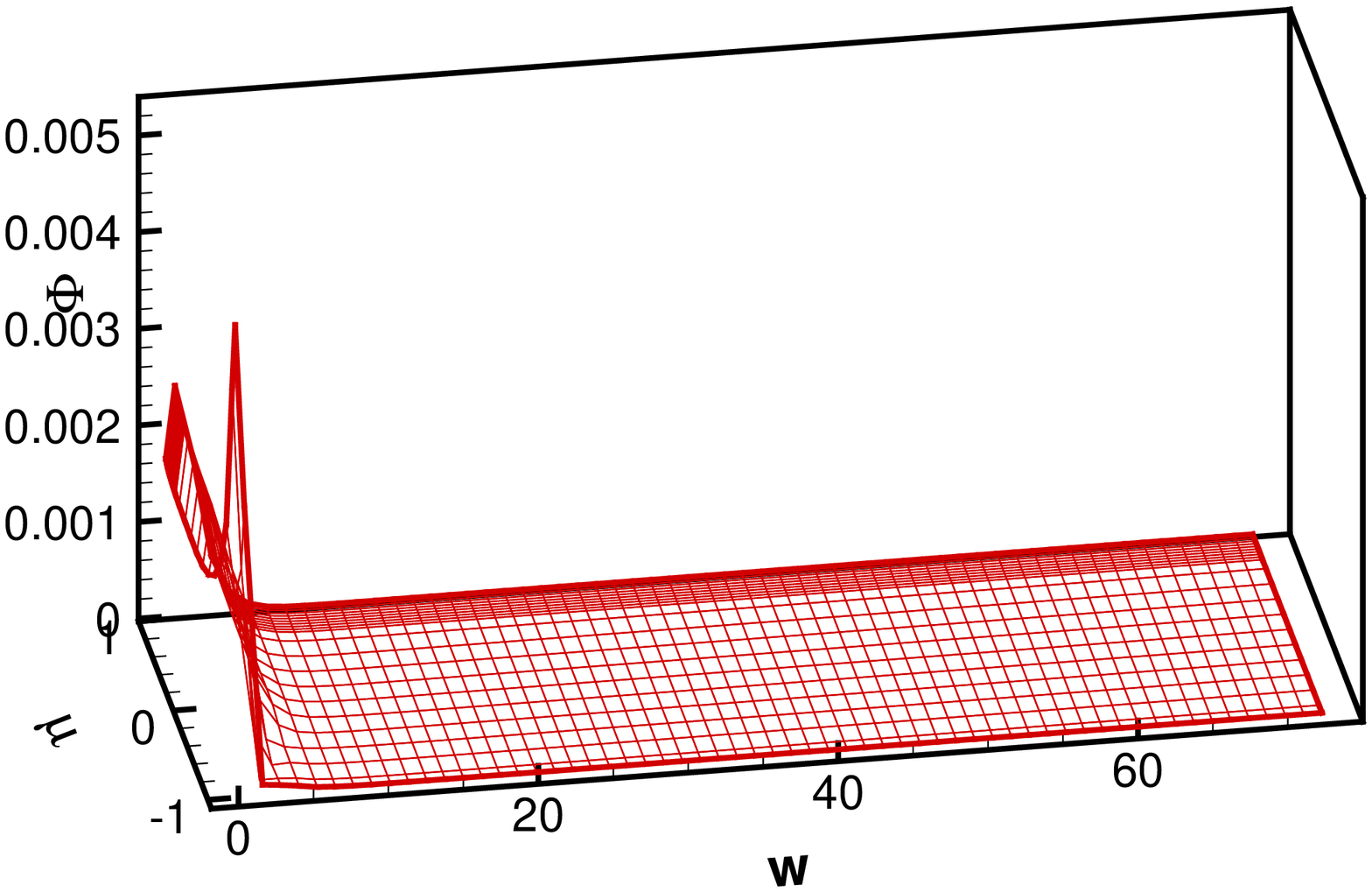}
\includegraphics[width=3in,angle=0]{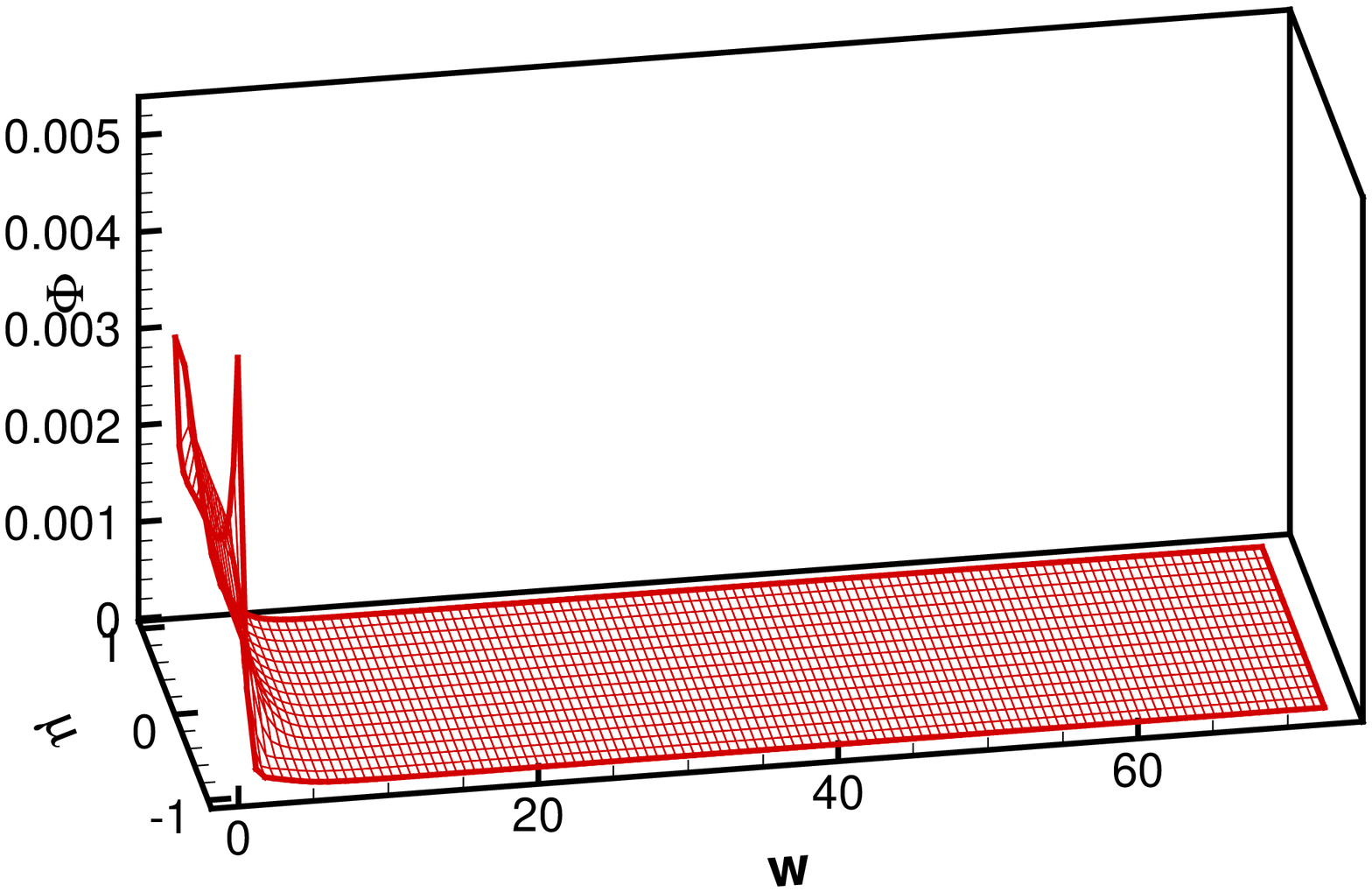}\\
\includegraphics[width=3in,angle=0]{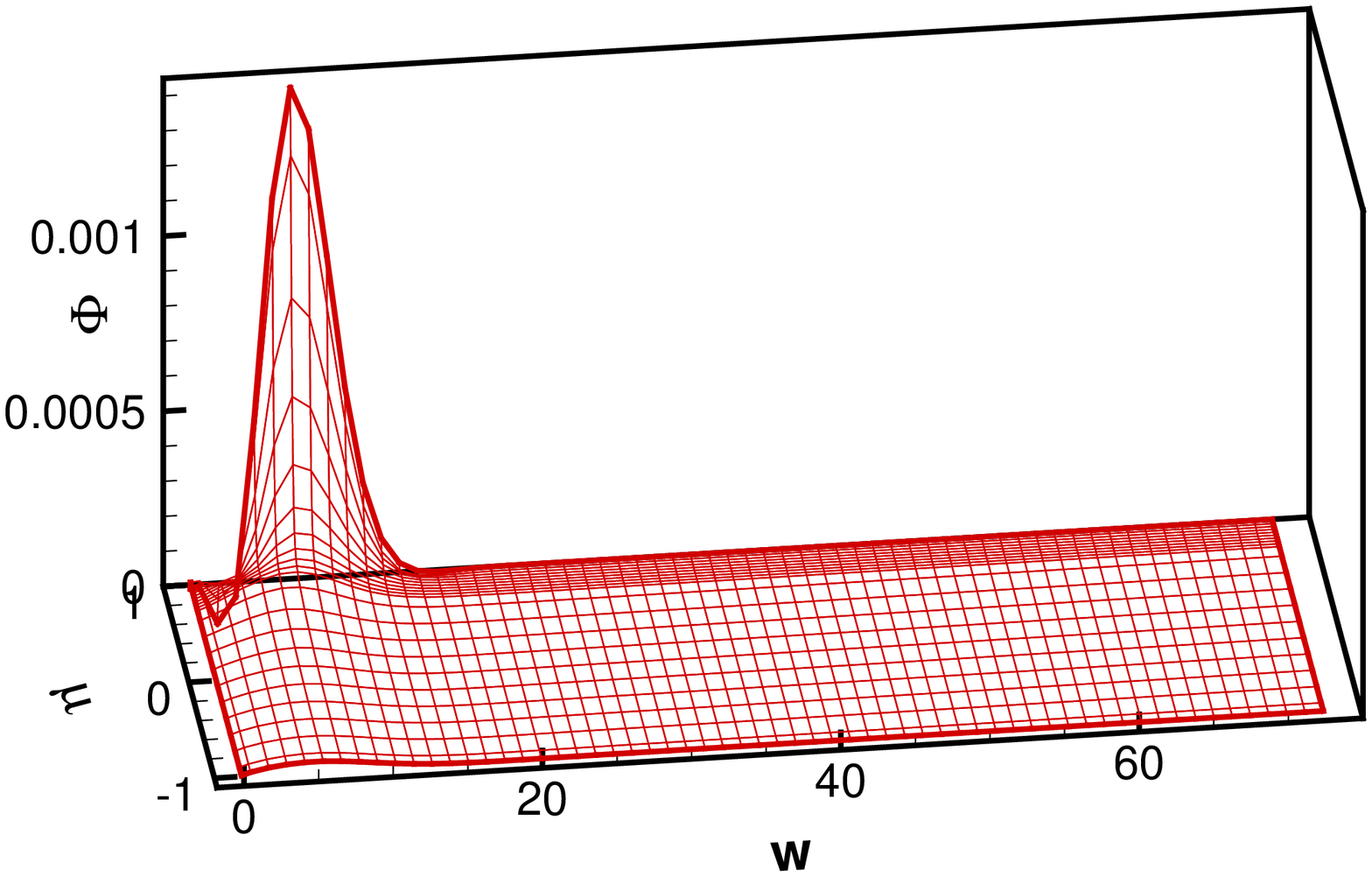}
\includegraphics[width=3in,angle=0]{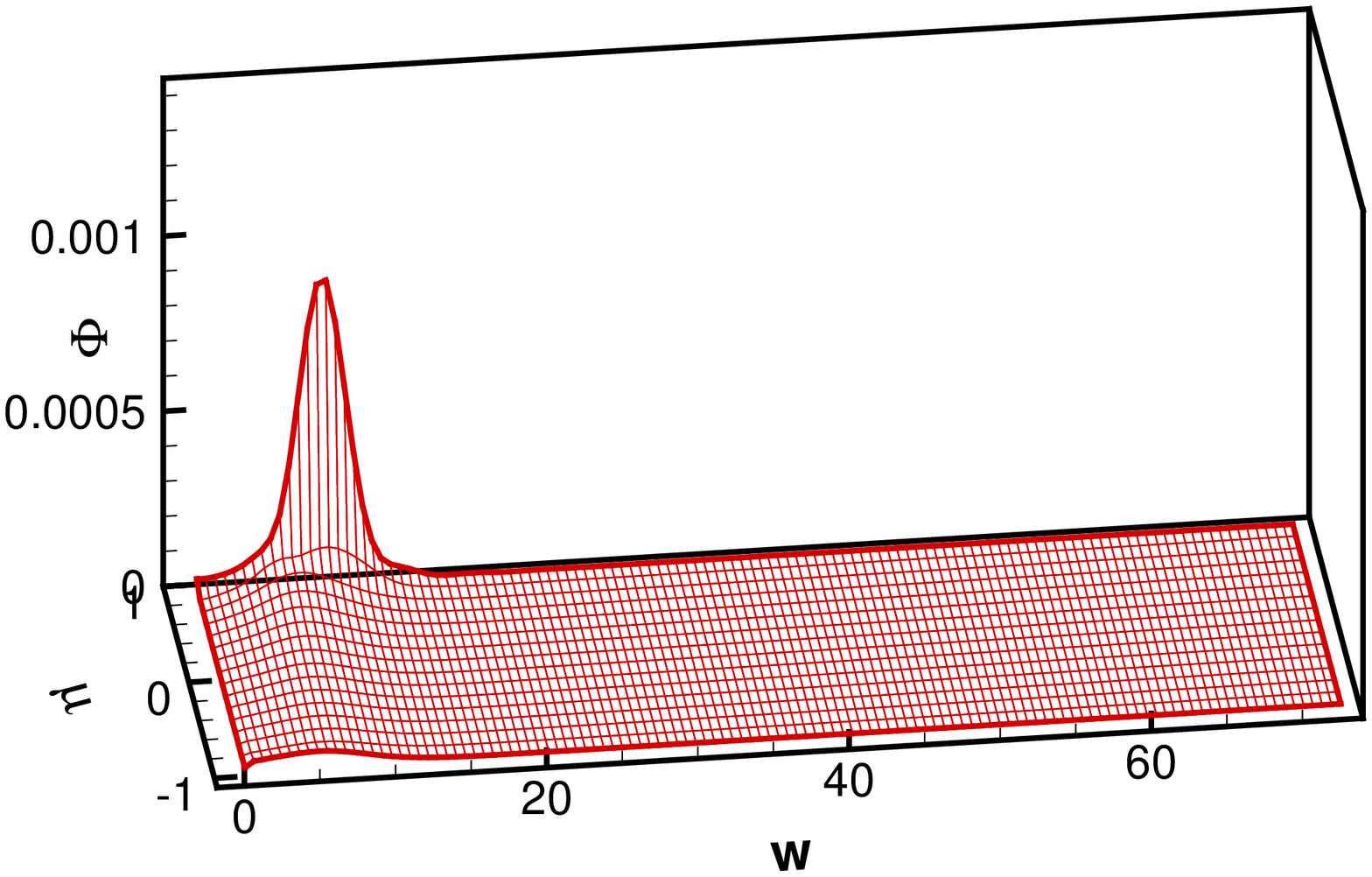}\\
\includegraphics[width=3in,angle=0]{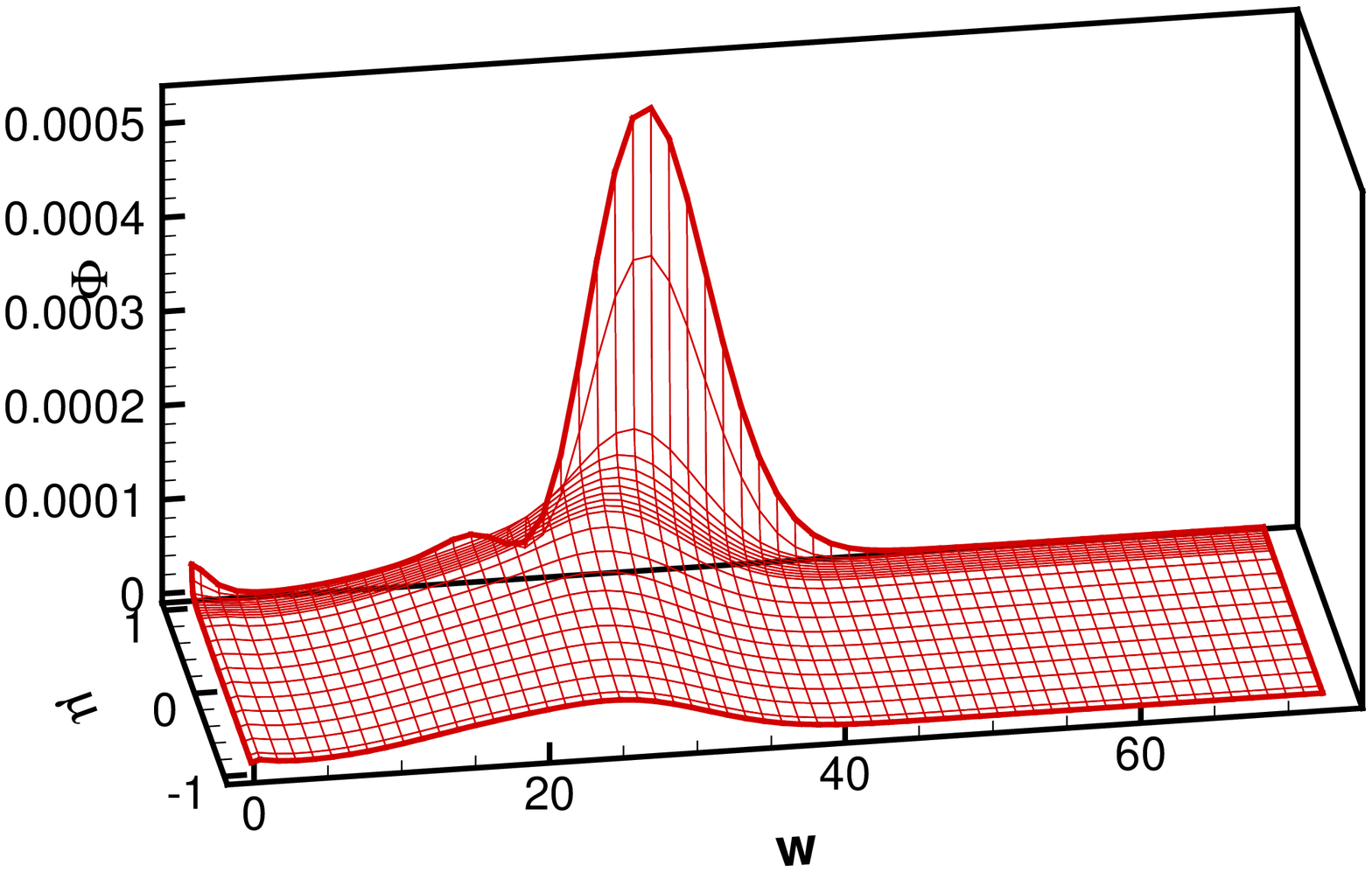}
\includegraphics[width=3in,angle=0]{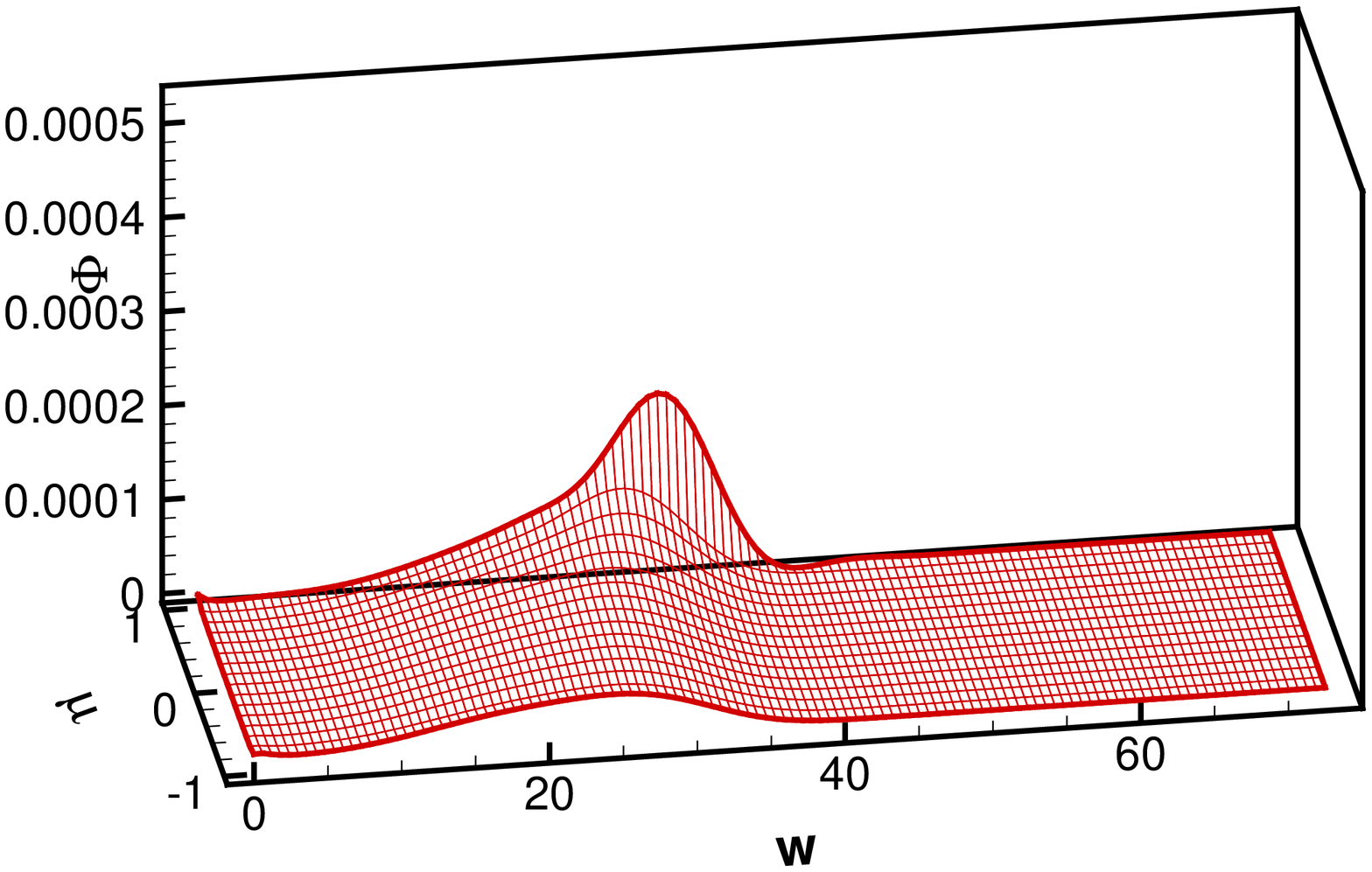}\\
\caption{Comparison of the snapshot for $\Phi(x_0,w,\mu)$ using DG
(left) and WENO (right) for $50$nm channel at $t=0.5$,
$V_{\mbox{bias}}=1.0$. Top:  $x_0=0.1$ ; middle: $x_0=0.125$;
bottom: $x_0=0.15$. Solution has not yet reached steady state.}
\label{50p1}
\end{figure}

\begin{figure}[htb]
\centering
\includegraphics[width=3in,angle=0]{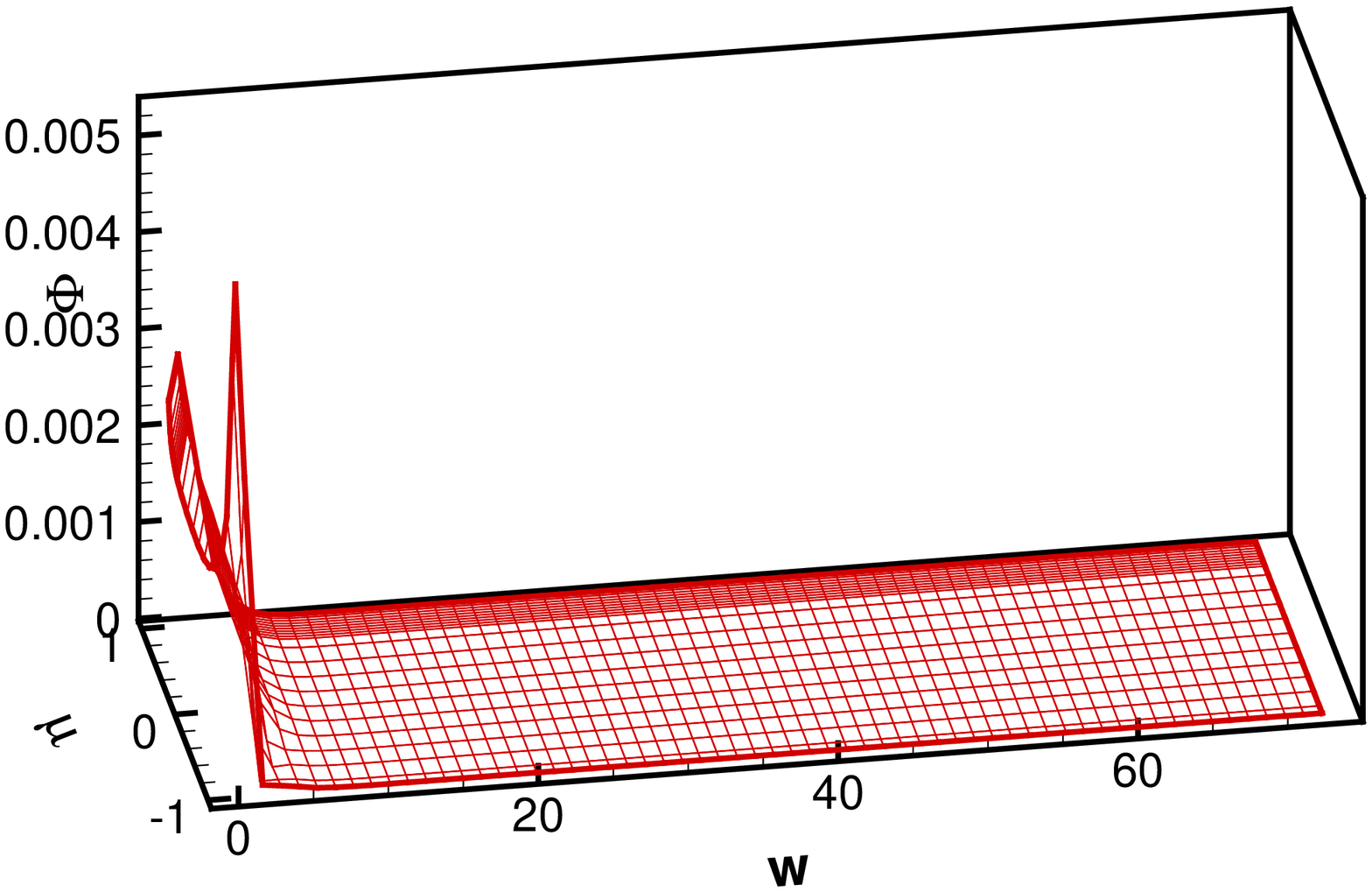}
\includegraphics[width=3in,angle=0]{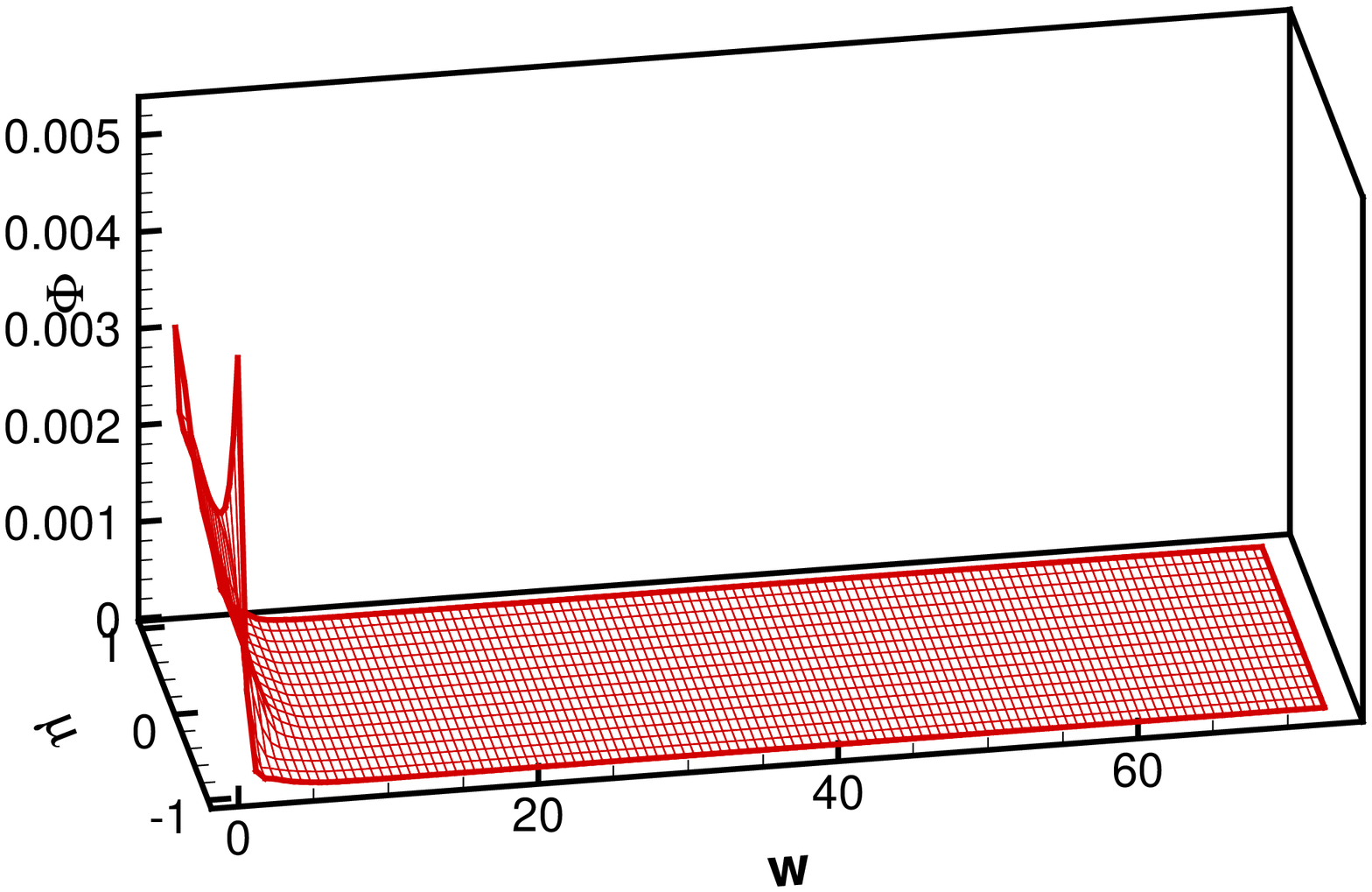}\\
\includegraphics[width=3in,angle=0]{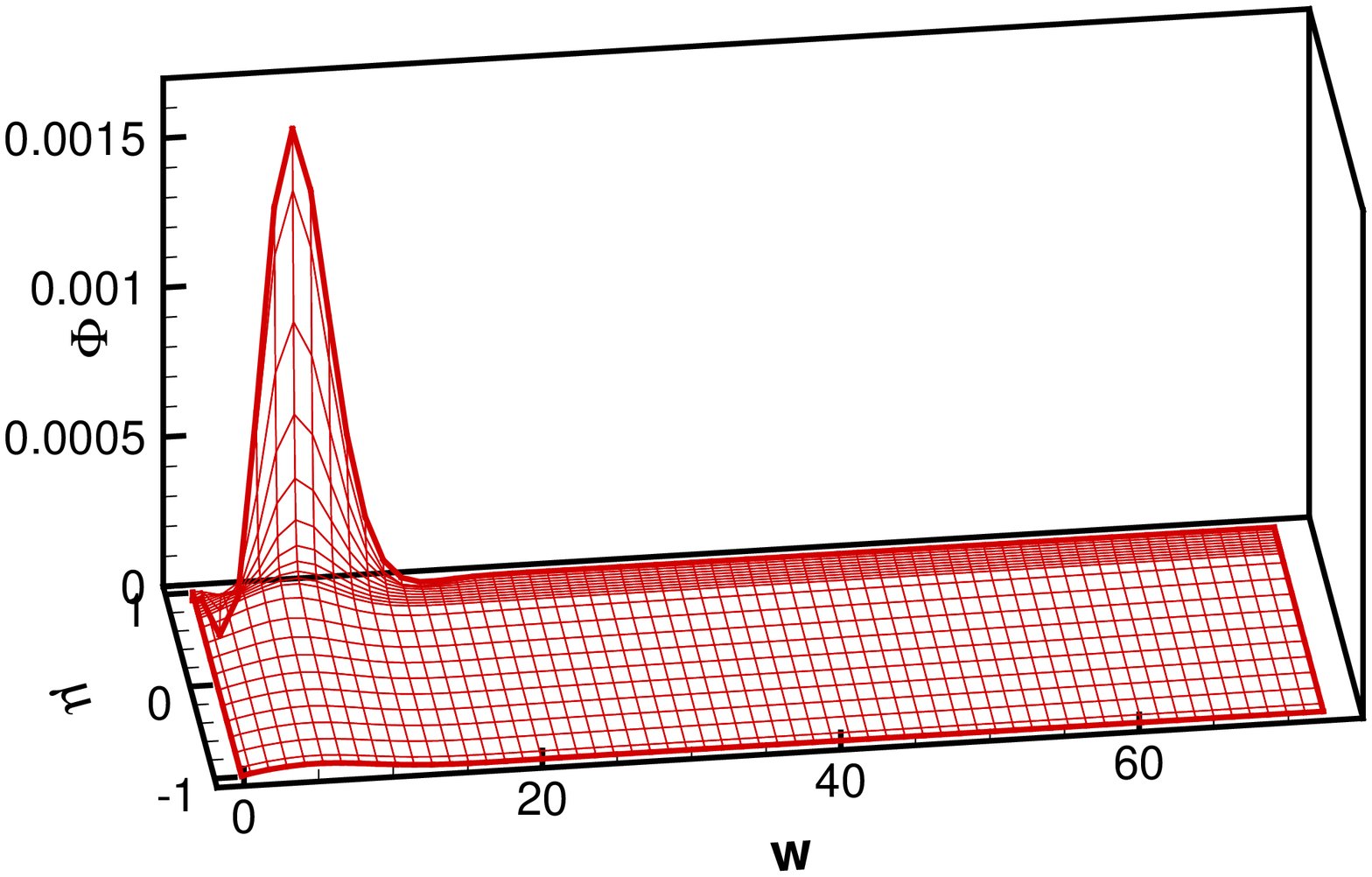}
\includegraphics[width=3in,angle=0]{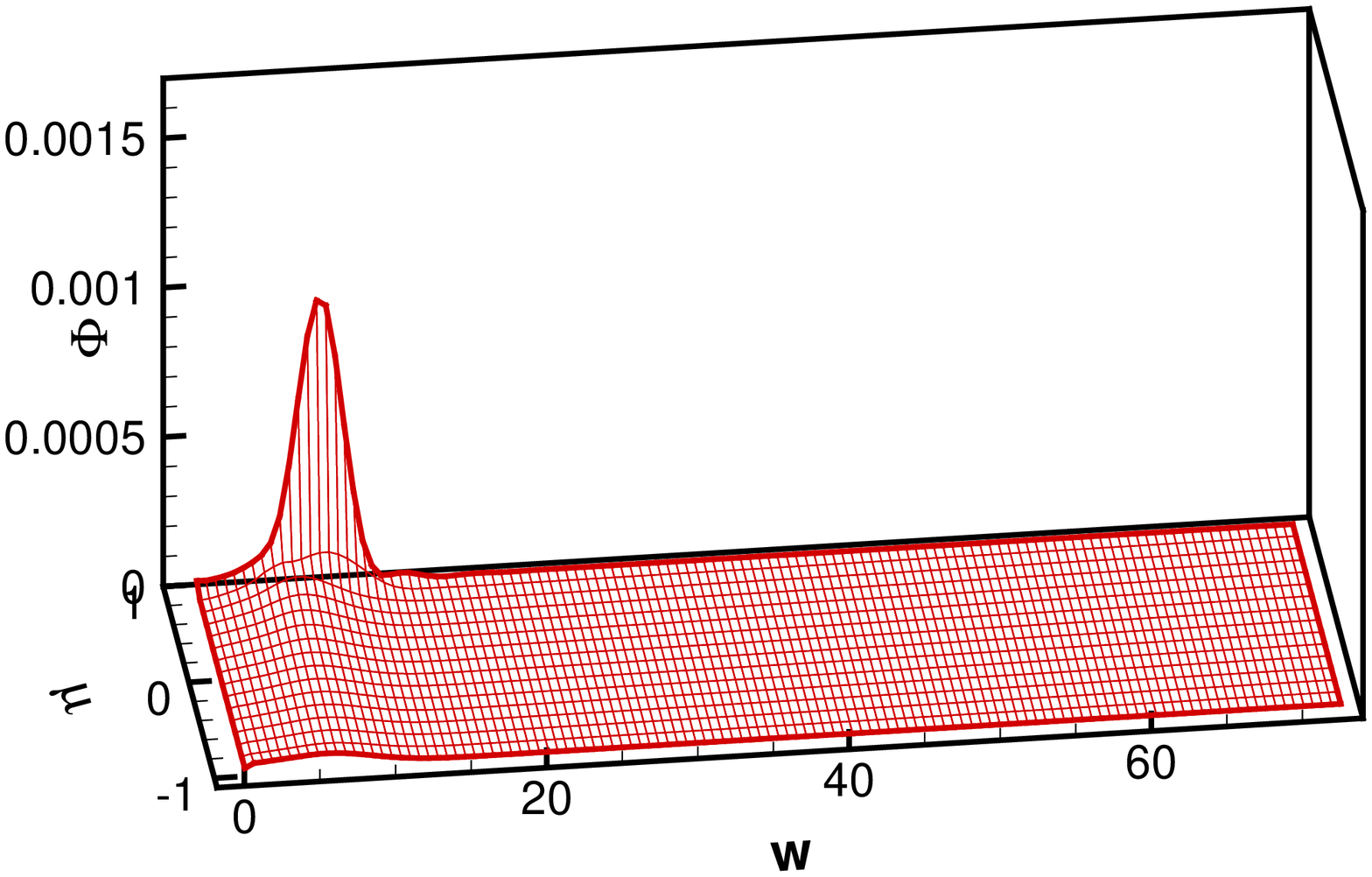}\\
\includegraphics[width=3in,angle=0]{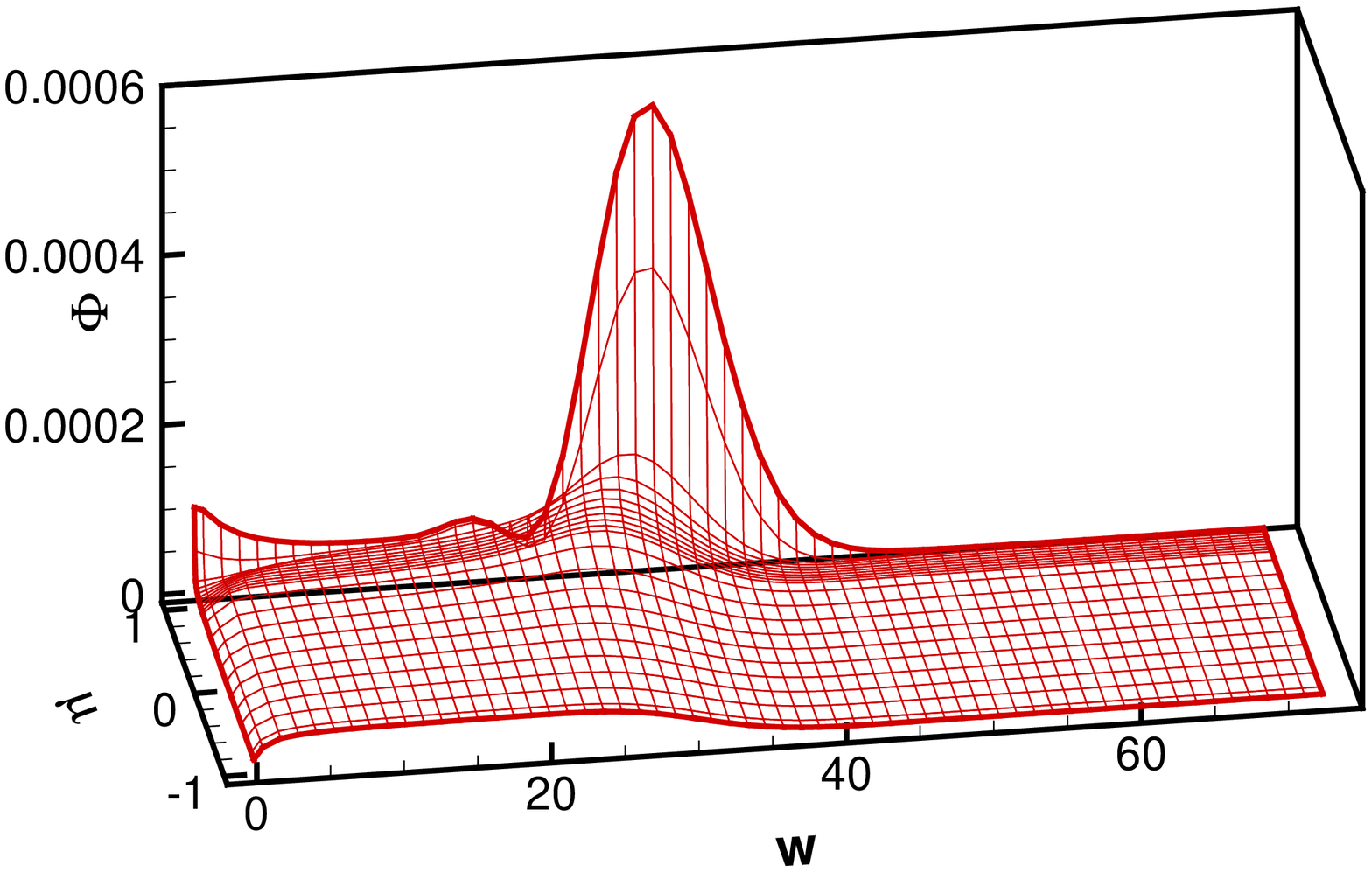}
\includegraphics[width=3in,angle=0]{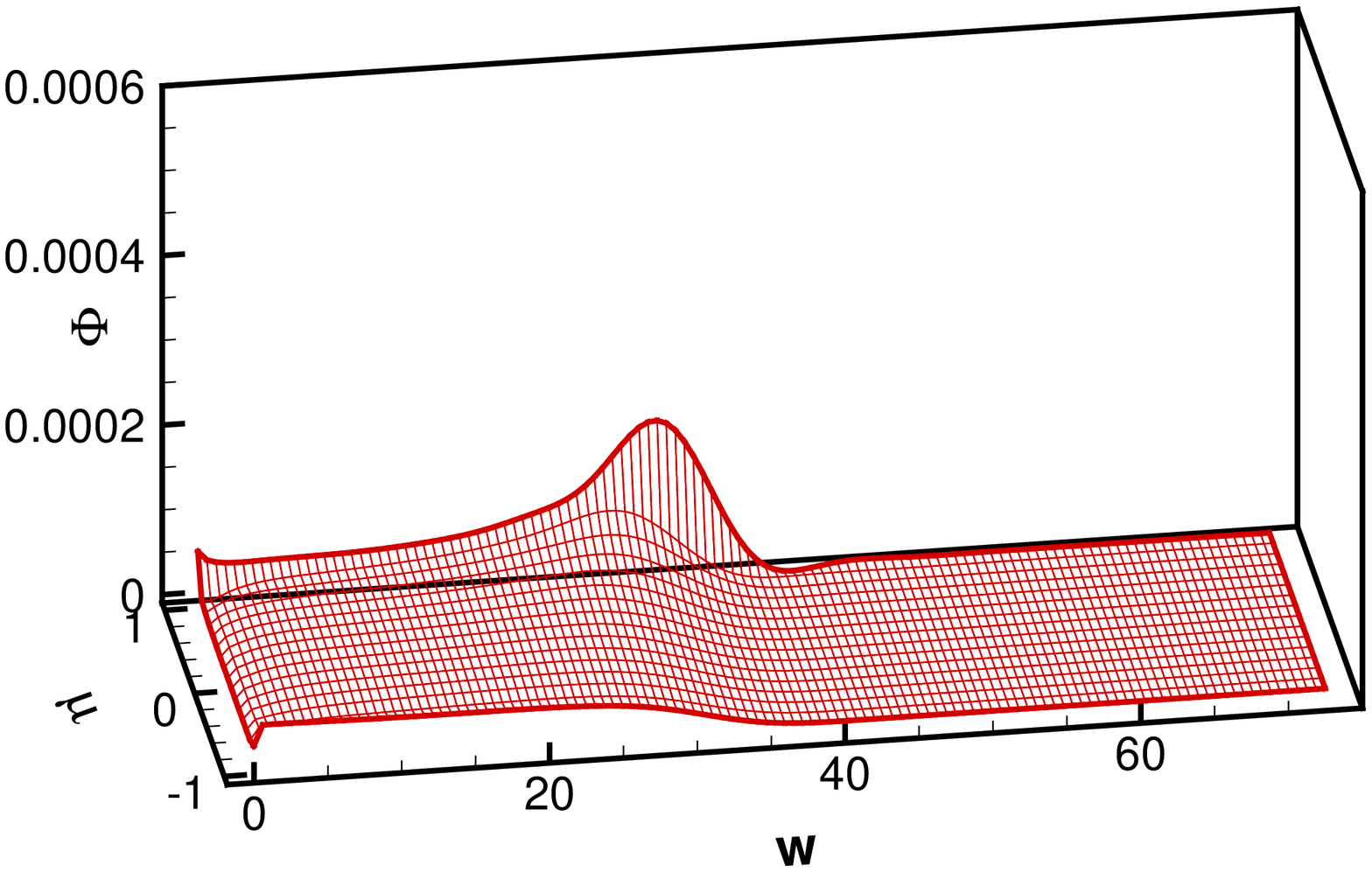}\\
\caption{Comparison of the snapshot for $\Phi(x_0,w,\mu)$ using DG
(left) and WENO (right) solution for $50$nm channel at $t=3.0$,
$V_{\mbox{bias}}=1.0$. Top:  $x_0=0.1$; middle: $x_0=0.125$; bottom:
$x_0=0.15$. Solution has reached steady state.} \label{50p2}
\end{figure}

\begin{figure}[htb]
\centering
\includegraphics[width=4in,angle=0]{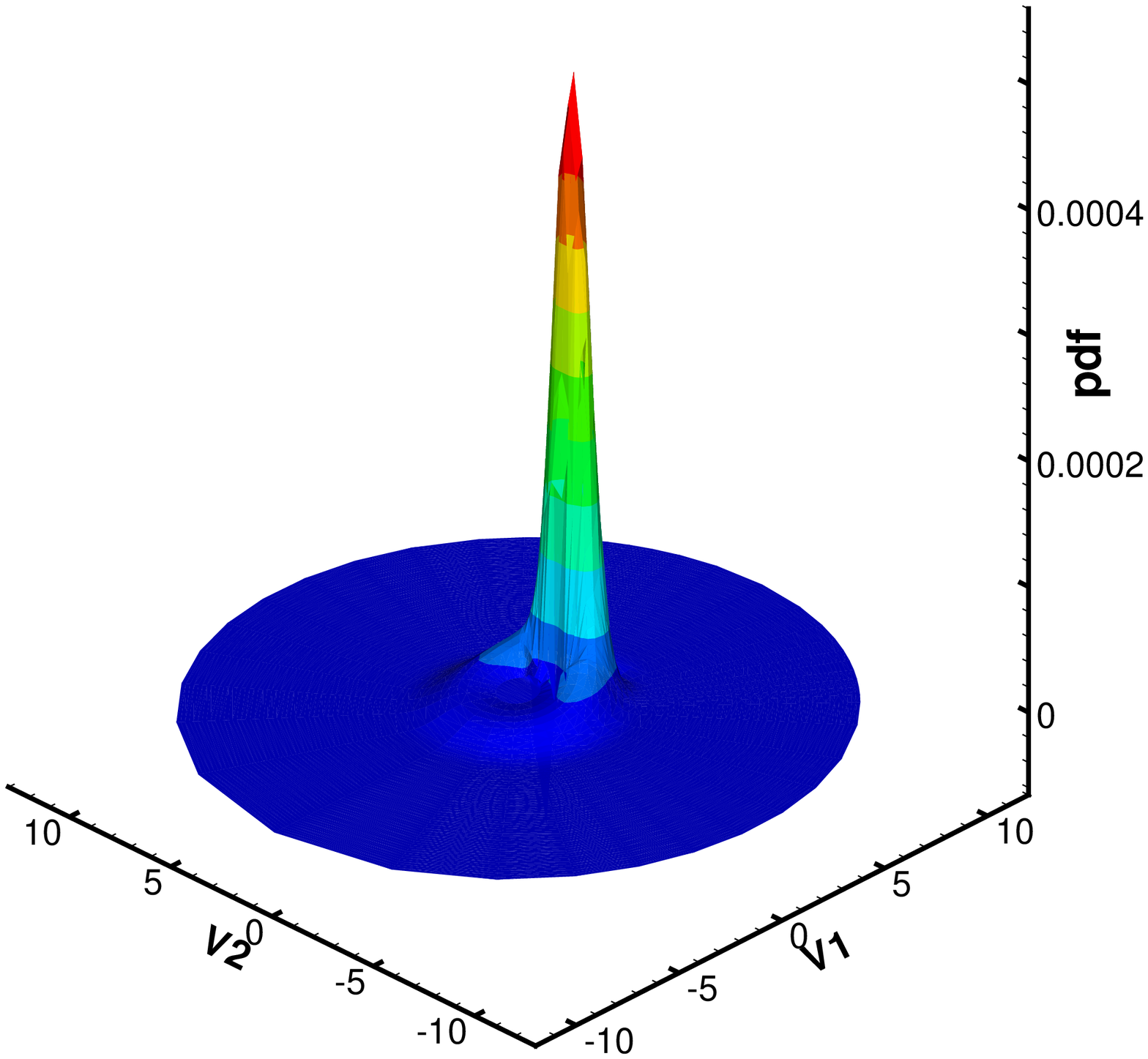}
\caption{PDF for 50nm channel at $t=3.0$, $x=0.125$.} \label{50cart}
\end{figure}

The nonuniform mesh we use for $50 nm$ channels is defined as
follows. In the $x$-direction, near the junctions, in $0.09<x<0.11$
and $0.14<x<0.16$ , $\triangle x= 0.001$; in center of the channel
$0.11<x<0.14$, $\triangle x= 0.005$; at everywhere else, $\triangle
x= 0.01$.  Thus, the total number of cells in $x$ direction is $64$.
In the $w$-direction, we use $60$ uniform cells. In the
$\mu$-direction, we use $20$ cells, $10$ in the region $\mu<0.7$,
$10$ in $\mu>0.7$. Thus, the grid  consists of  $64 \times 60 \times
20$ cells, compared to the WENO calculation of $150 \times 120
\times 16$ uniform cells. The evolution and steady state plots are
listed in Figures \ref{50evo} to \ref{50p2}. The conclusions are
similar with $400 nm$, that we obtain better resolutions near the
channel junctions and the peak for \emph{pdf}  is much higher.
Figure \ref{50cart} plots $\Phi(t=5.0, x=0.125)$ in the cartesian
coordinates. The peak is twice the height of WENO and is very sharp.
Figure \ref{50cart1} to \ref{50cart4} plot the \emph{pdf} near
$x=0.15$, the drain junction. We obtain distributions far away from
statistical equilibrium, that reflects the lack of suitability of
the classical hydrodynamical models for the drain region of a small
gated device under even moderate voltage bias.

\begin{figure}[htb]
\centering
\includegraphics[width=2.93in,angle=0]{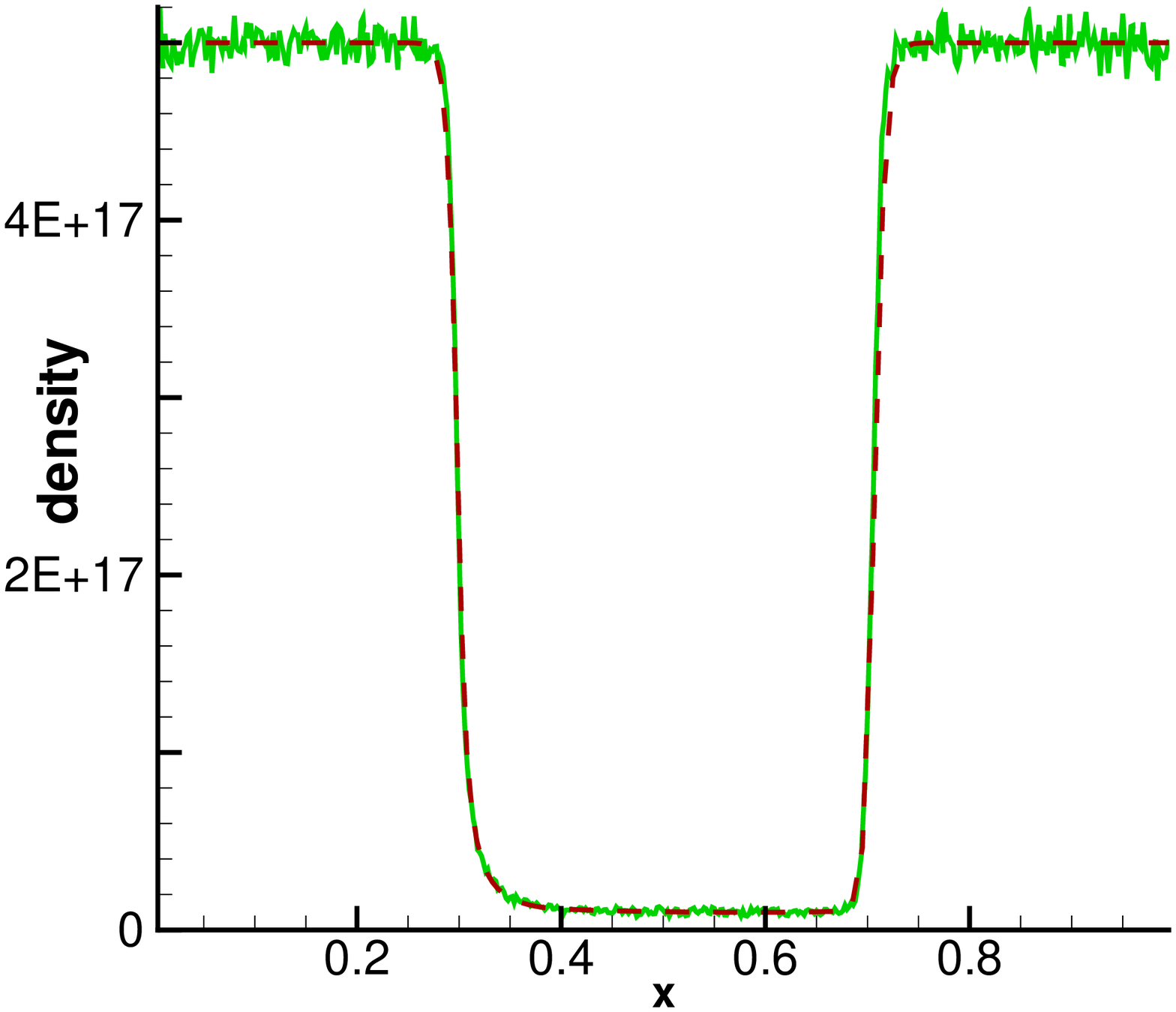}
\includegraphics[width=2.93in,angle=0]{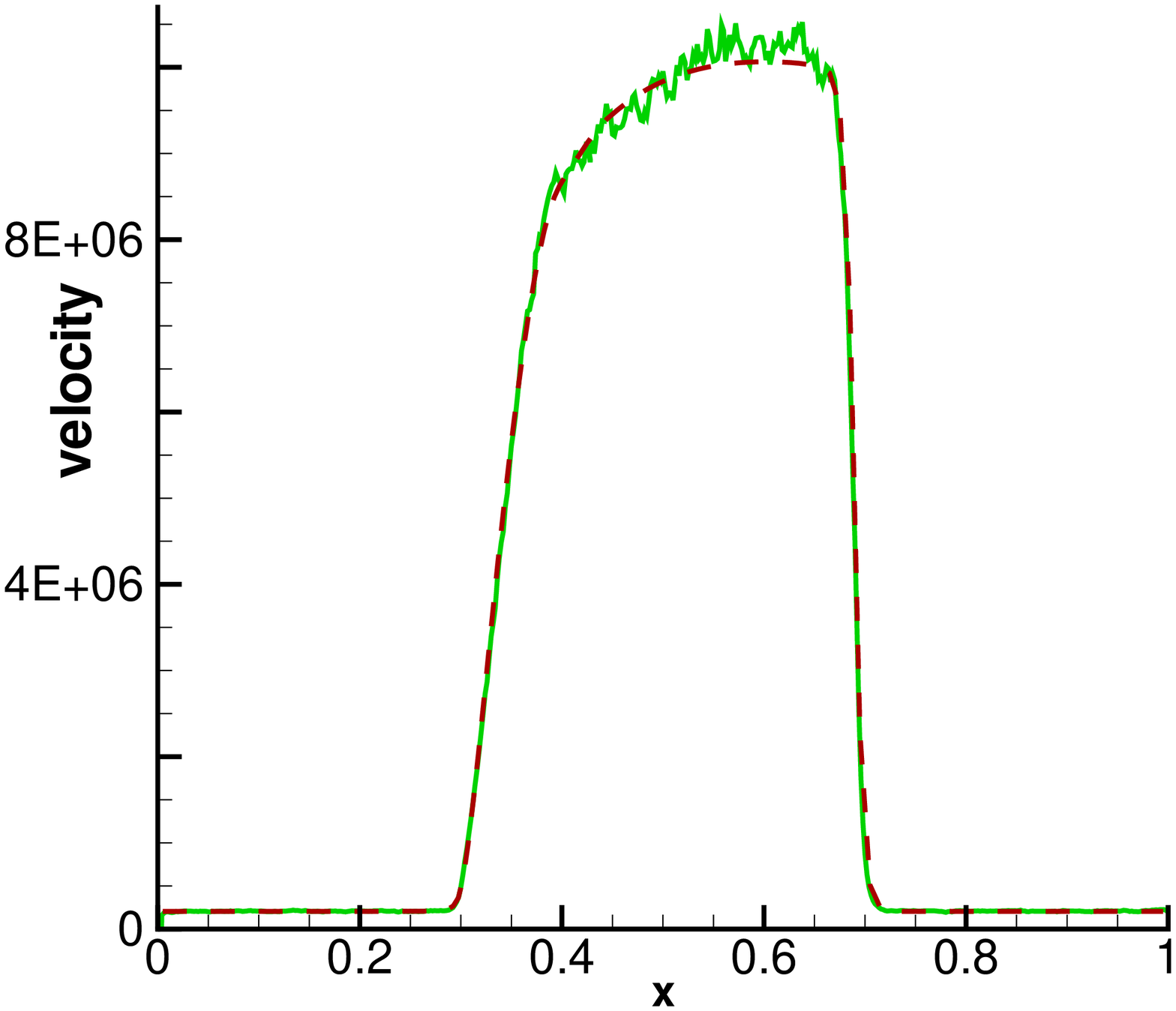}\\
\includegraphics[width=2.93in,angle=0]{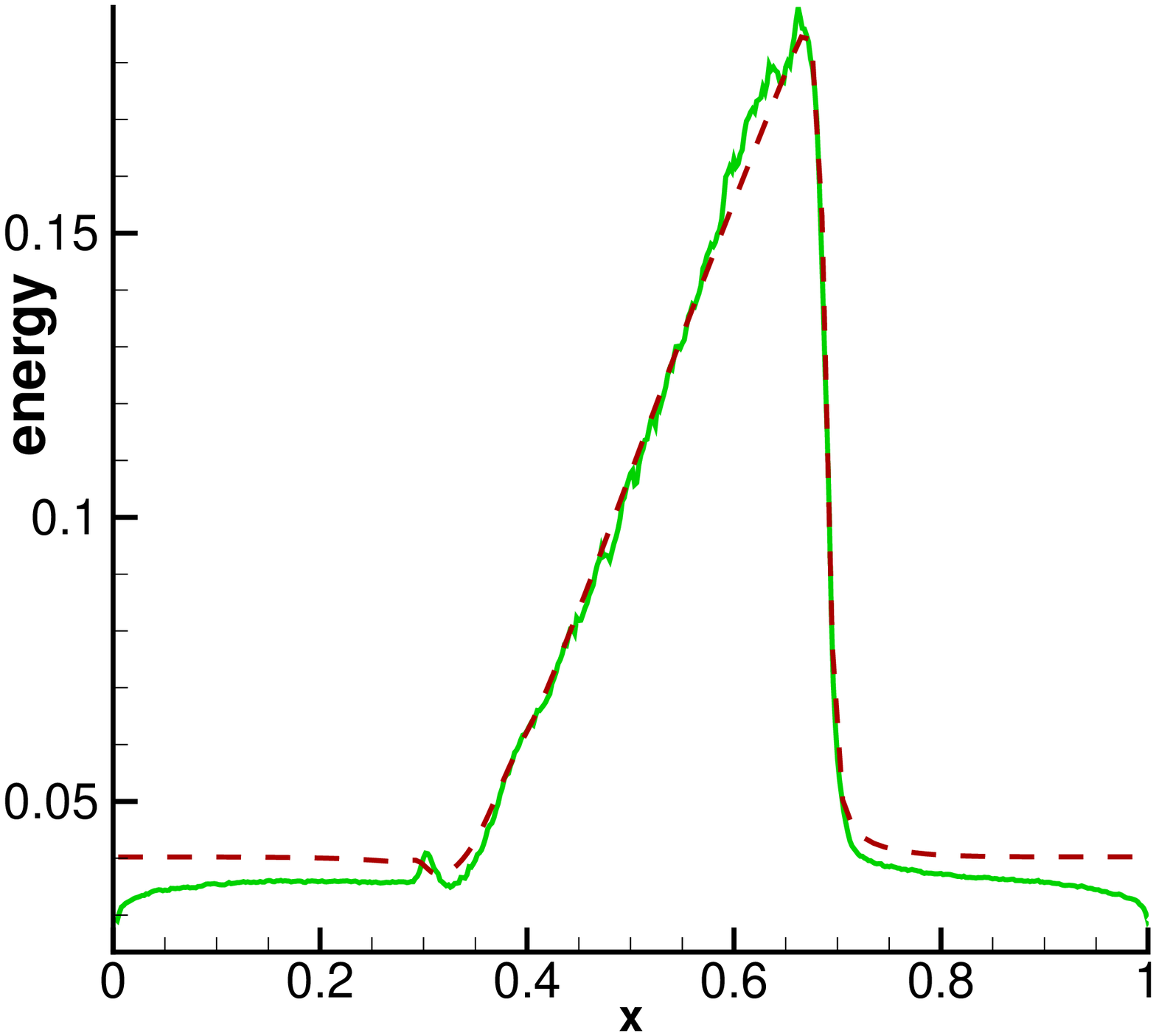}
\includegraphics[width=2.93in,angle=0]{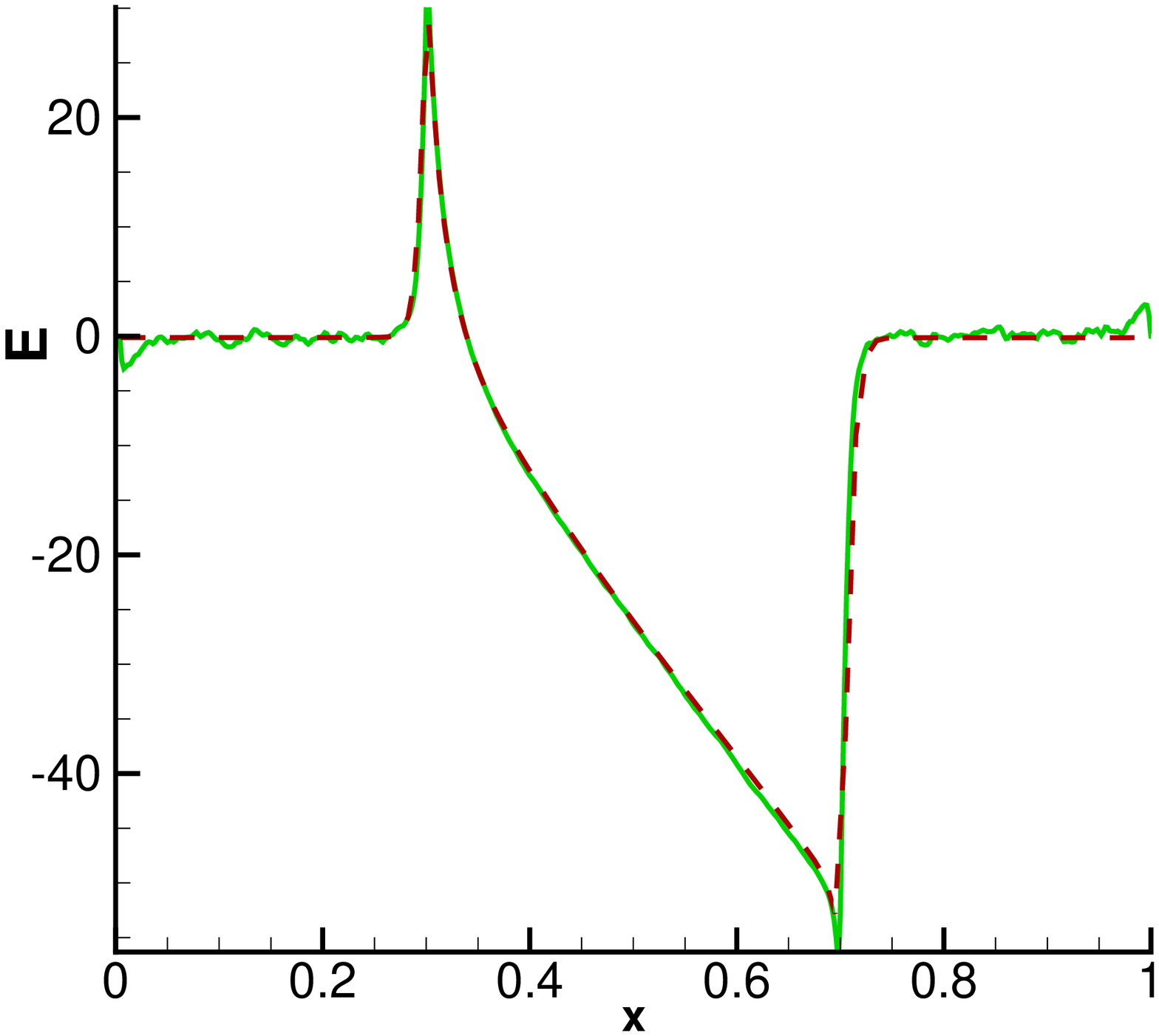}\\
\includegraphics[width=2.93in,angle=0]{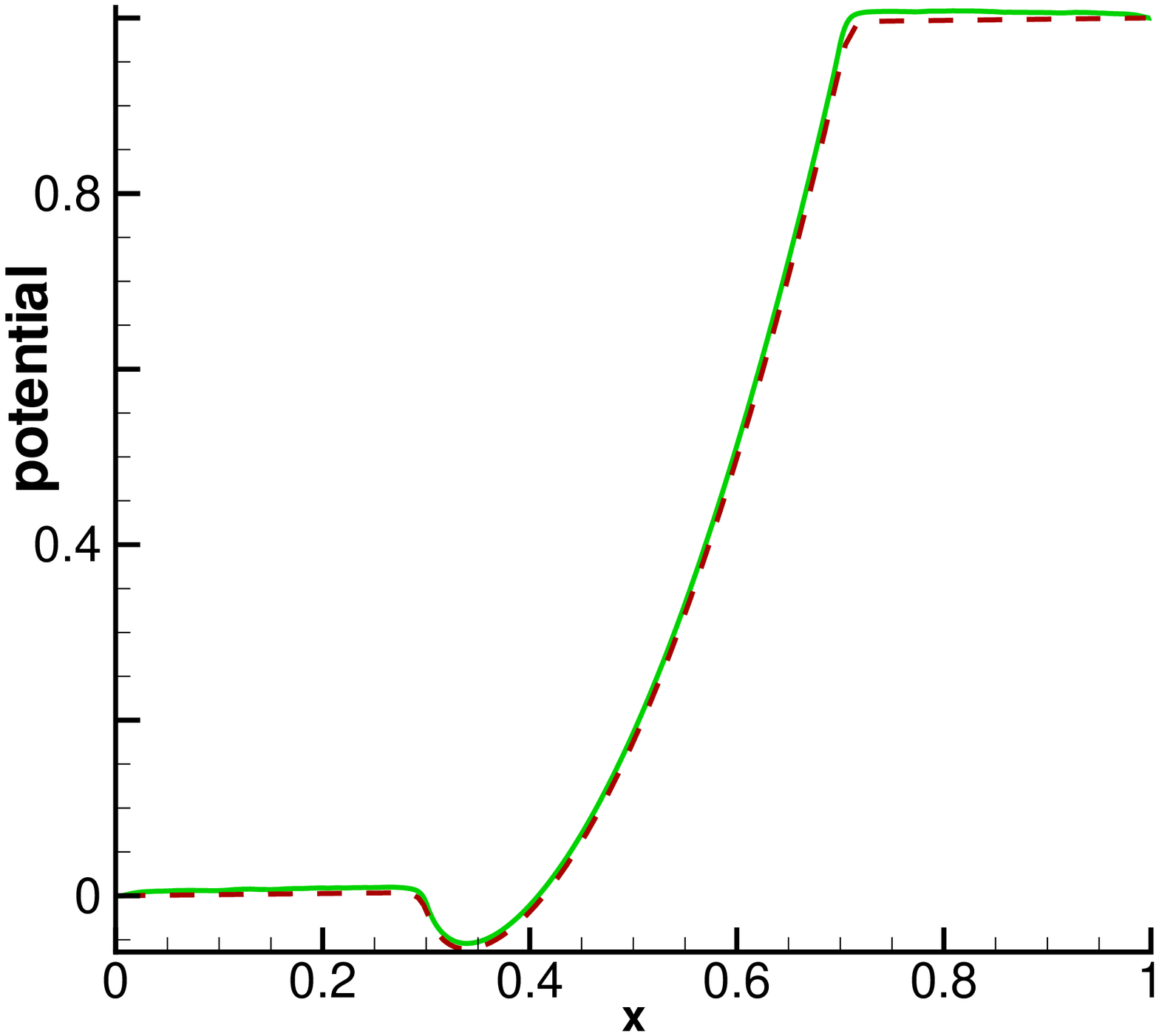}
\includegraphics[width=2.93in,angle=0]{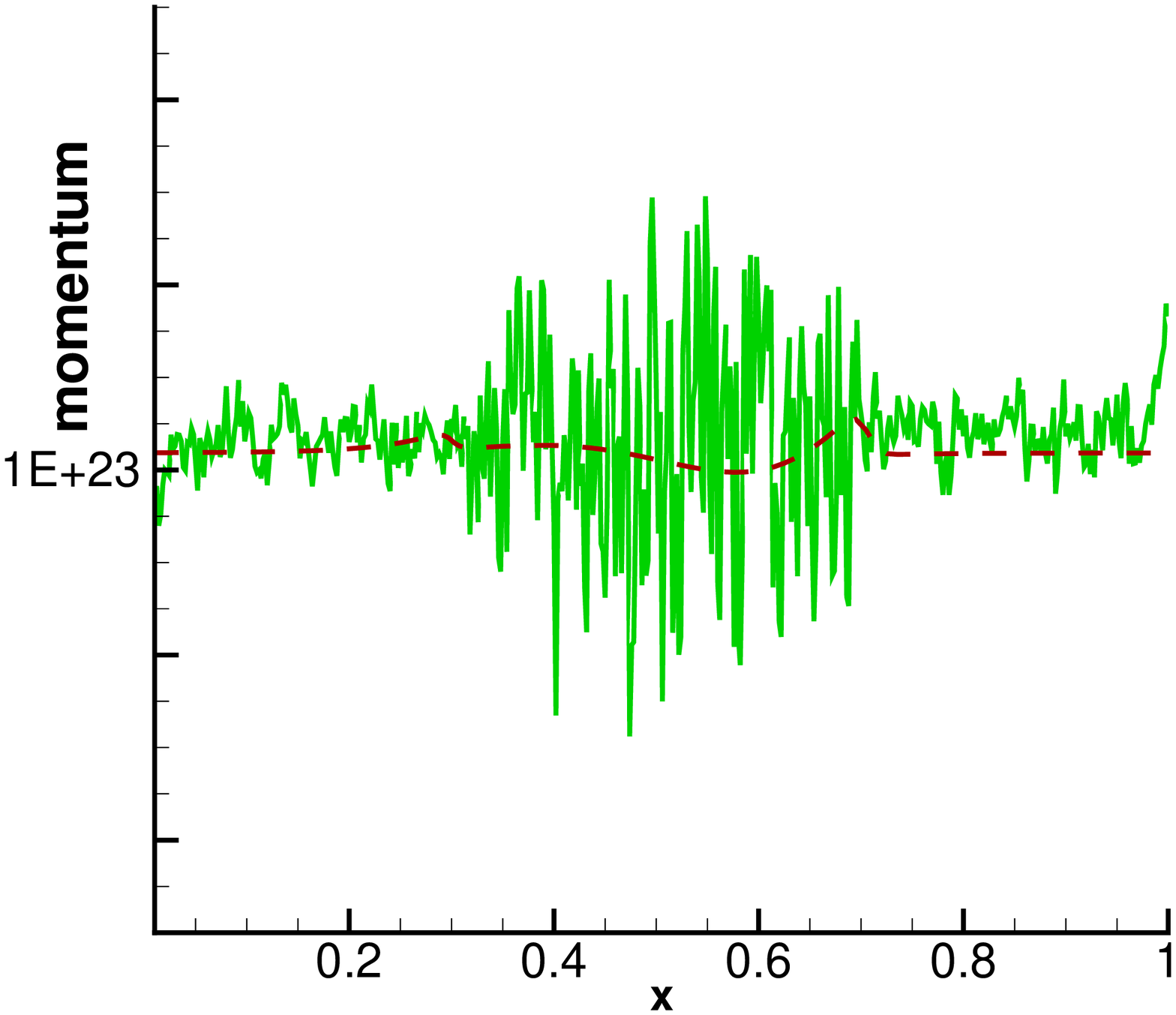}
\caption{Comparison of macroscopic quantities using DG (dashed line)
and DSMC (solid line) for $400$nm channel at $t=5.0$,
$V_{\mbox{bias}}=1.0$. Top left: density in ${cm}^{-3}$; top right:
mean velocity in $cm/s$; middle left: energy in $eV$; middle right:
electric field in $kV/cm$; bottom left: potential in $V$; bottom
right: momentum in ${cm}^{-2} \, s^{-1}$. Solution has reached
steady state.} \label{400dsmc}
\end{figure}

\begin{figure}[htb]
\centering
\includegraphics[width=2.93in,angle=0]{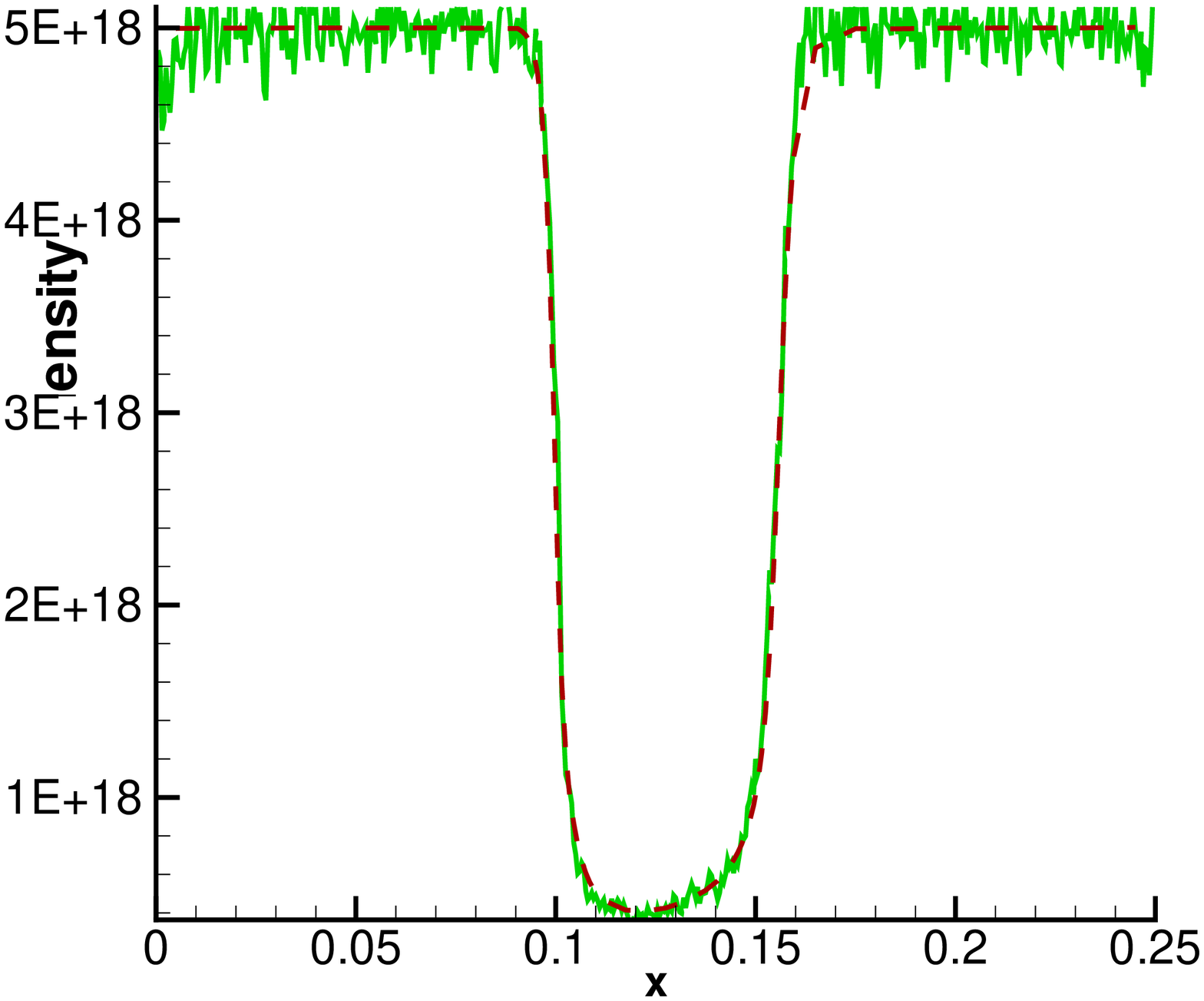}
\includegraphics[width=2.93in,angle=0]{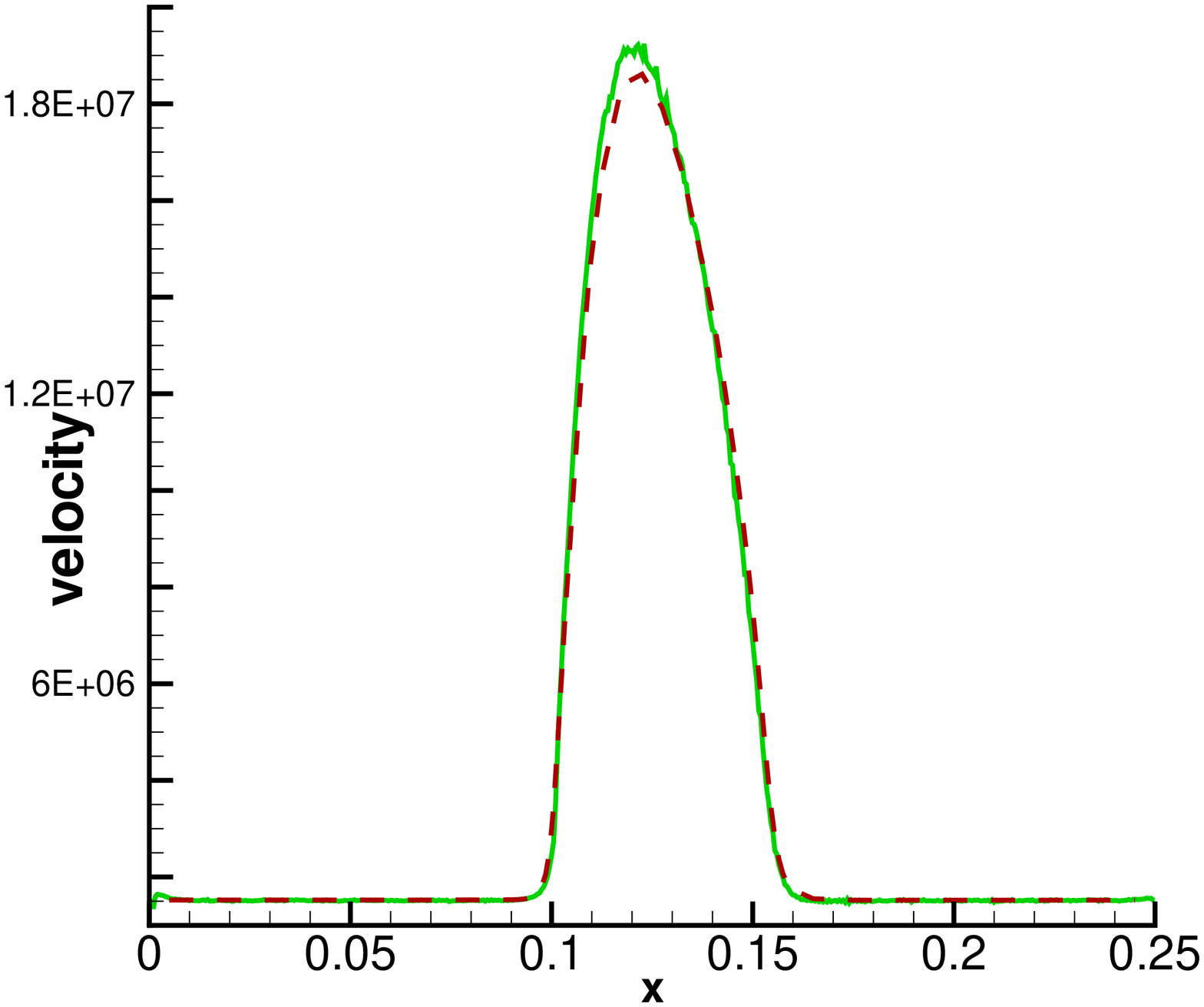}\\
\includegraphics[width=2.93in,angle=0]{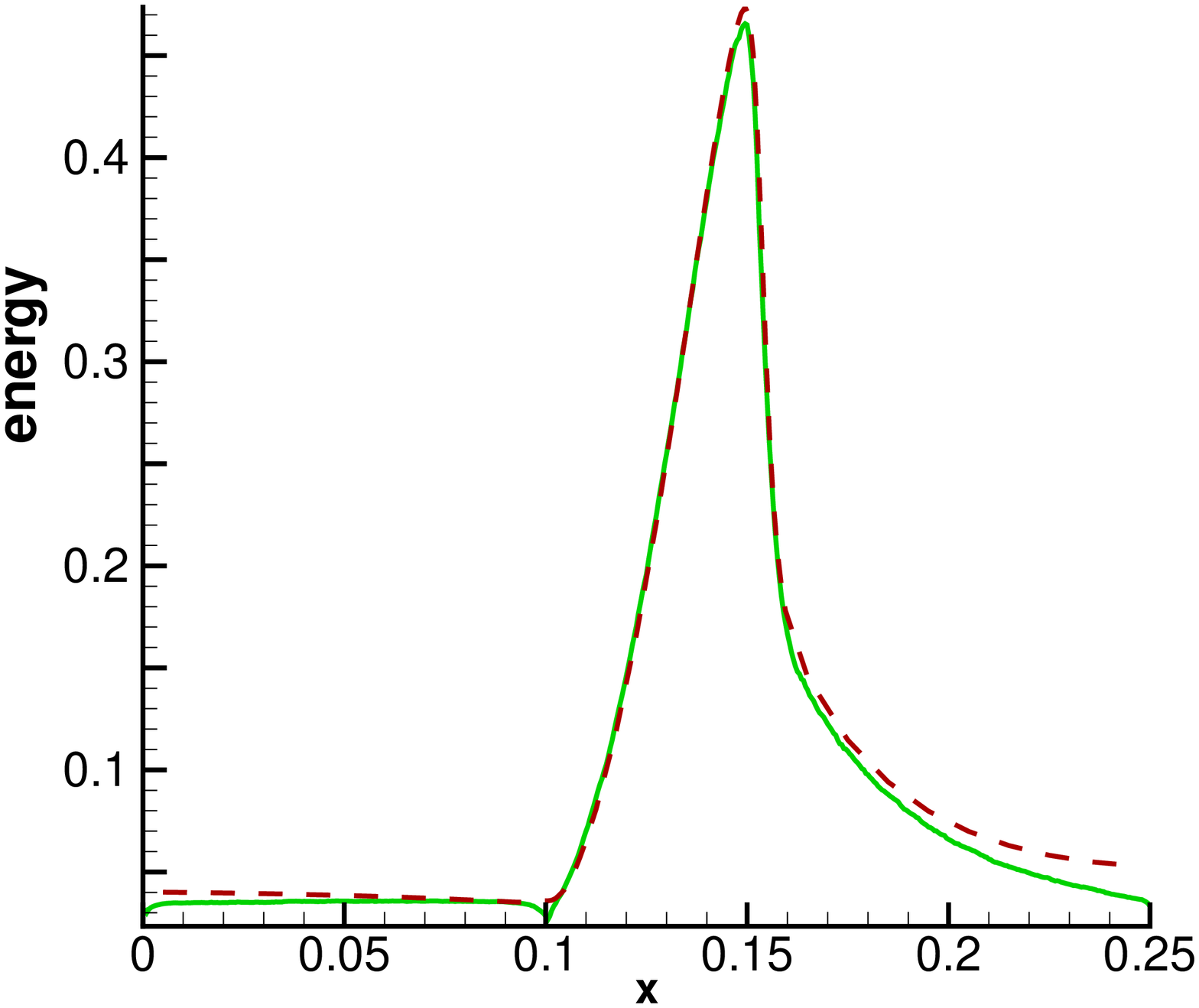}
\includegraphics[width=2.93in,angle=0]{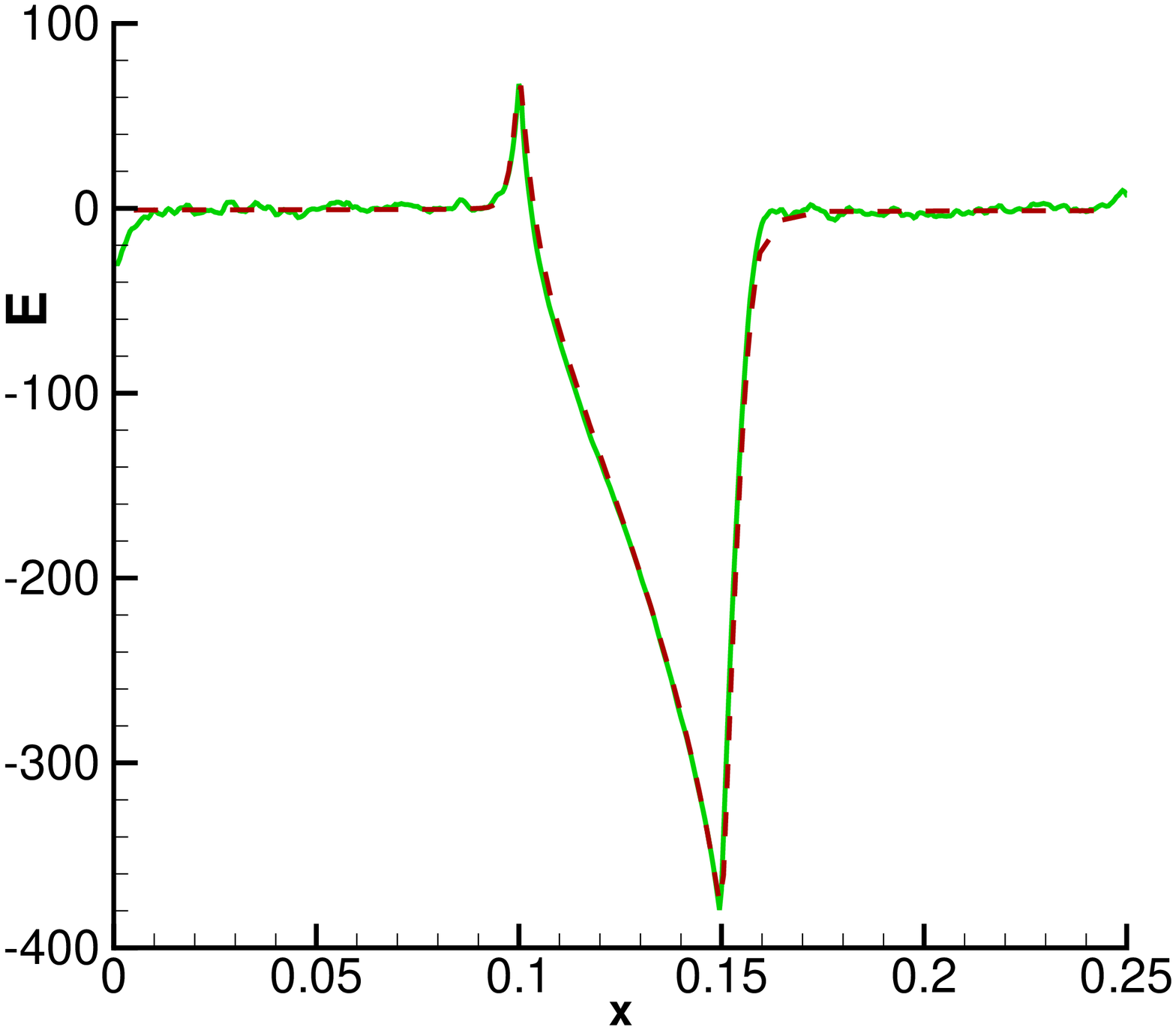}\\
\includegraphics[width=2.93in,angle=0]{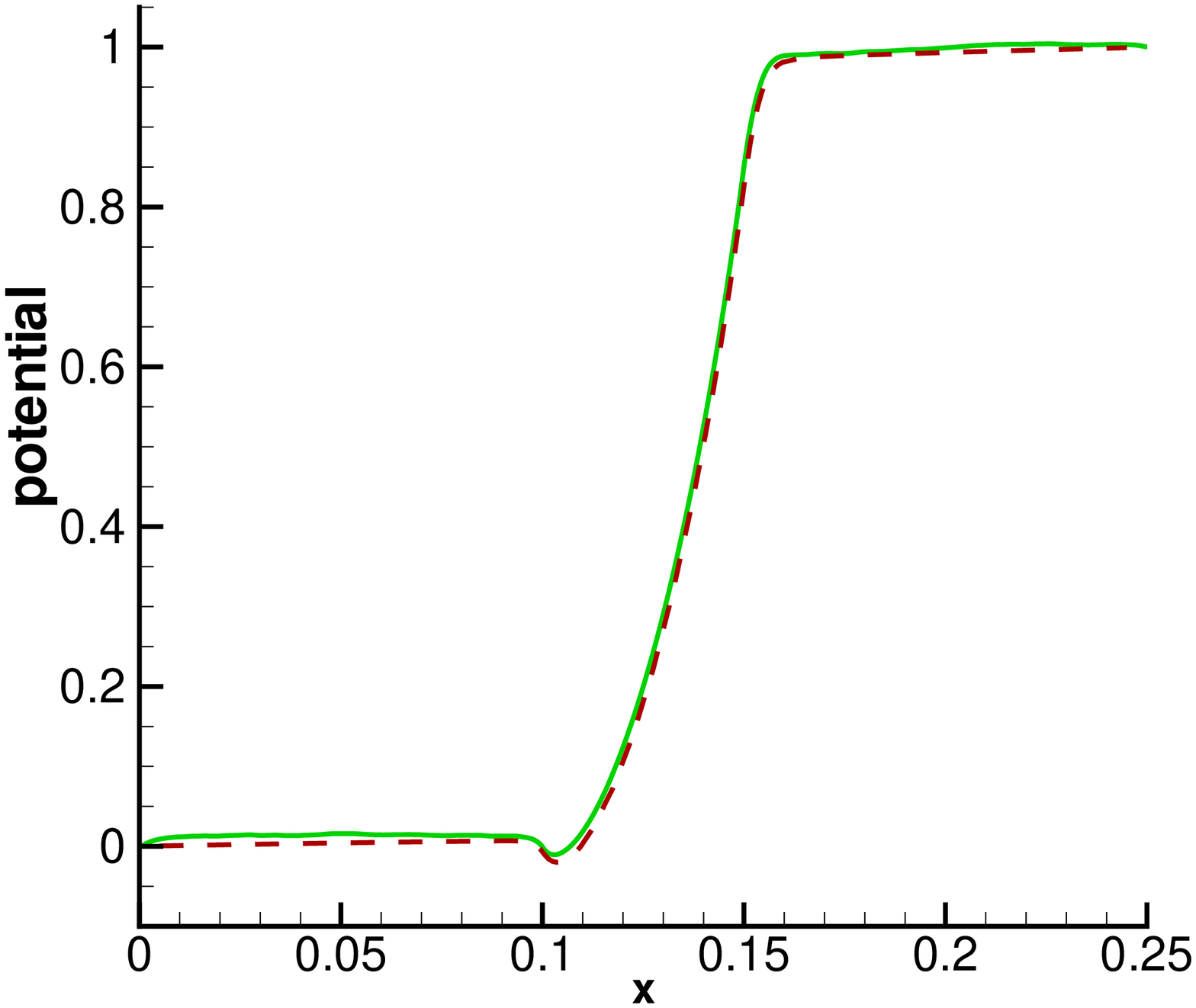}
\includegraphics[width=2.93in,angle=0]{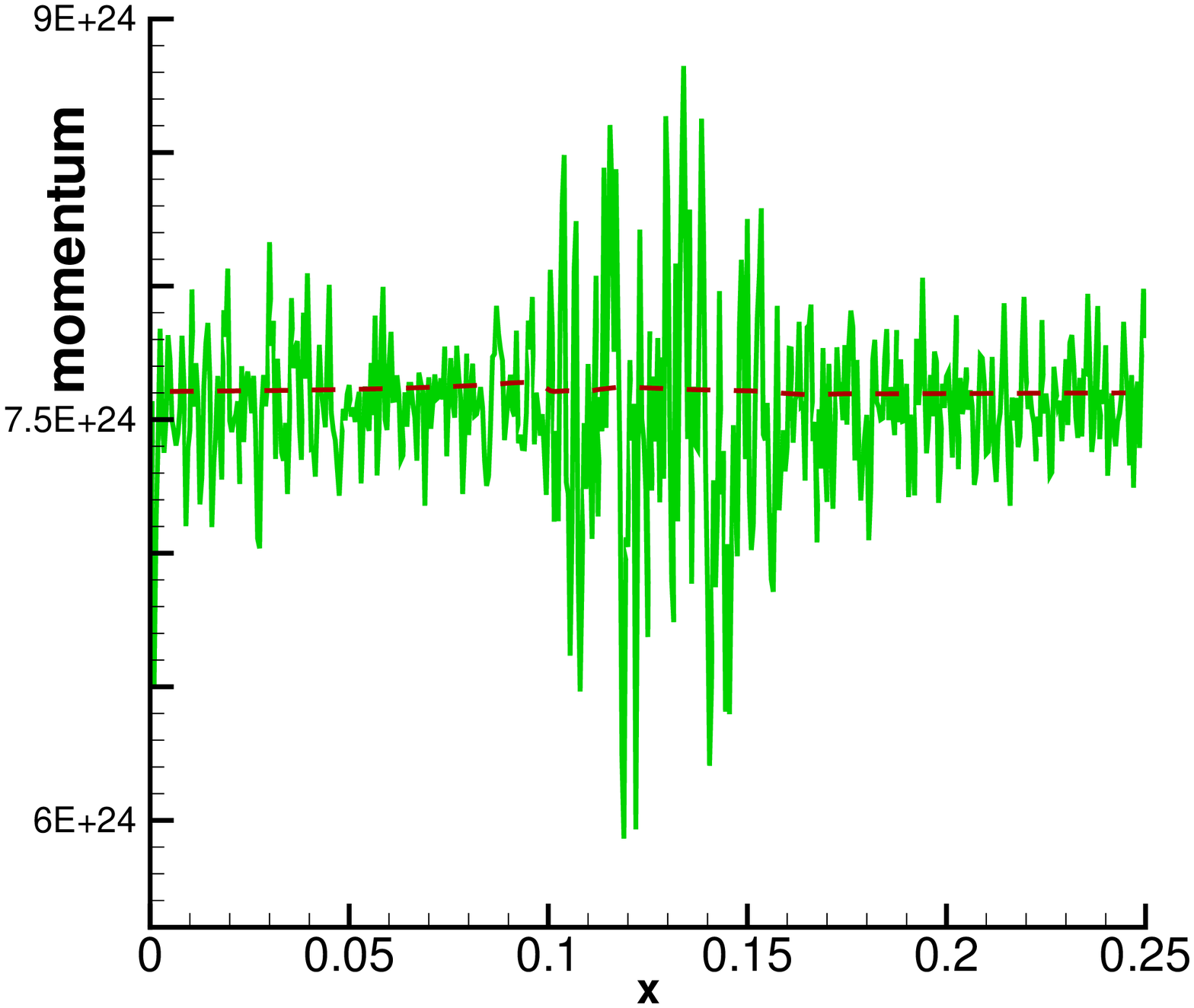}
\caption{Comparison of macroscopic quantities using DG (dashed line)
and DSMC (solid line) for $50$nm channel at $t=3.0$,
$V_{\mbox{bias}}=1.0$. Top left: density in ${cm}^{-3}$; top right:
mean velocity in $cm/s$; middle left: energy in $eV$; middle right:
electric field in $kV/cm$; bottom left: potential in $V$; bottom
right: momentum in ${cm}^{-2} \, s^{-1}$. Solution has reached
steady state.} \label{50dsmc}
\end{figure}

We also compare the results from DG-BTE solver with those obtained
from DSMC simulations, see Figures \ref{400dsmc}, \ref{50dsmc}. The
two simulations show good agreement except for energy plots near the
boundaries. The modeling of the contact boundaries is not simple,
since it requires to know the distribution function of entering
particles. The best way to solve this problem is the inclusion of a
transport kinetic equation for the dynamic of the electron at the
metal junctions; of course this is not realistic due to the
complexity of this new kinetic equation, where the importance of
electron-electron interaction requires a nonlinear collisional
operator, similar to the classical one of the Boltzmann equation for
a rarefied perfect gas. Then, the simplest reasonable rule consists
in assuming that the distribution function near, but outside, the
device is proportional to a Maxwellian (or shifted Maxwellian)
equilibrium distribution function, or to the distribution function
near, but inside, the device boundaries. When there are strong
electric fields also near the boundaries, the first choice is not
reasonable, since, as we show in this paper, the distribution
function is very far away to a Maxwellian distribution function.
Therefore, the second choice is better than the first. We remark
that both choices are simple but only low level approximation of the
true physical phenomena; so many  criticisms are known in the
literature. We assume, as usual, that charge neutrality holds at the
contact; so, the particle density near the contact boundaries
coincides with the doping density. This law is used in all of the
DSMC, WENO and DG simulations. Nevertheless, since we must
approximate this constraint in different way, i.e. at molecular
level for DSMC, introducing suitable ghost points for WENO scheme or
giving appropriate values of $\Phi$ at boundaries in DG simulations,
we cannot have  a unique exact boundary condition in the
computational experiments. Now, it is obvious that this difference
in the boundary treatment has an influence for the solutions at the
stationary regime.
%
%
%
%
\clearpage
%

\clearpage
\section{DG-BTE solver for 2D double gate MOSFET simulation}

In this section, we consider a 2D double gate MOSFET device. In
order to solver the 2D Boltzmann-Poisson system, we choose to
implement a simple rectangular grid and let
$$
 \Omega_{ijkmn} = \left[ x_{i - \ot} , \, x_{i + \ot} \right] \times
                  \left[ y_{j - \ot} , \, y_{j + \ot} \right] \times
                  \left[ w_{k - \ot} , \, w_{k + \ot} \right] \times
                  \left[ \mu_{m - \ot} , \, \mu_{m + \ot} \right] \times
                  \left[ \ph_{n - \ot} , \, \ph_{n + \ot} \right]
$$
where $i=1, \ldots N_x$, $j=1, \ldots N_y$, $k=1, \ldots N_w$, $m=1,
\ldots N_\mu$, $n=1, \ldots N_\ph$, and
$$
 x_{i \pm \ot} = x_{i} \pm \frac{\Delta x_{i}}{2} \, , \quad
 y_{j \pm \ot} = y_{j} \pm \frac{\Delta y_{j}}{2}\, ,  \quad
 w_{k \pm \ot} = w_{k} \pm \frac{\Delta w_{k}}{2}\,
$$
$$
 \mu_{m \pm \ot} = \mu_{m} \pm \frac{\Delta \mu_{m}}{2}\, , \quad
 \ph_{n \pm \ot} = \ph_{n} \pm \frac{\Delta \ph_{n}}{2}.
$$

The approximation space is defined as
\begin{equation}
V_h^\kpol=\{ v : v|_{\Omega_{ijkmn}} \in P^\kpol(\Omega_{ijkmn})\}.
\end{equation}
Here, $P^\kpol(\Omega_{ijkmn})$ is the set of all polynomials of
degree at most $\kpol$ on $\Omega_{ijkmn}$. The DG formulation for
the Boltzmann equation (\ref{eqPhi}) would be: to find $\Phi_h \in
V_h^\kpol$, such that
\begin{eqnarray}
\label{dgb}
 &&\int_{\Omega_{ijkmn}} (\Phi_h)_t  \,v_h  \,d \Omega
- \int_{\Omega_{ijkmn}} g_1 \Phi_h  \,(v_h)_x  \,d \Omega
- \int_{\Omega_{ijkmn}} g_2 \Phi_h  \,(v_h)_y  \,d \Omega \nonumber \\
&& \mbox{} - \int_{\Omega_{ijkmn}} g_3 \Phi_h  \,(v_h)_w  \,d \Omega
- \int_{\Omega_{ijkmn}} g_4 \Phi_h  \,(v_h)_\mu  \,d \Omega -
\int_{\Omega_{ijkmn}} g_5 \Phi_h  \,(v_h)_\ph  \,d \Omega
\\
&& \mbox{} + F_x^+ - F_x^- +F_y^+-F_y^- +F_w^+-F_w^-
+F_\mu^+-F_\mu^- +F_\ph^+-F_\ph^- =\int_{\Omega_{ijkmn}} C(\Phi_h)
\,v_h \,d \Omega . \nonumber
\end{eqnarray}
for any test function $v_h \in V_h^\kpol$. In (\ref{dgb}),
$$
F_x^+=\iy \iw \imu \iphi g_1 \, \check{\Phi} \, v_h^-(x_{i+ \ot}, y,
w, \mu, \ph)dy \, dw \, d\mu \, d\ph  ,
$$
$$
F_x^-=\iy \iw \imu \iphi g_1 \, \check{\Phi} \, v_h^+ (x_{i- \ot},
y, w, \mu, \ph)dy \, dw \, d\mu \, d\ph  ,
$$

$$
F_y^+=\ix \iw \imu \iphi g_2 \, \bar{\Phi} \, v_h^- (x, y_{j+ \ot},
w, \mu, \ph)dx \, dw \, d\mu \, d\ph ,
$$
$$
F_y^-=\ix \iw \imu \iphi g_2 \, \bar{\Phi} \, v_h^+ (x, y_{j- \ot},
w, \mu, \ph)dx \, dw \, d\mu \, d\ph ,
$$

$$
F_w^+=\ix \iy \imu \iphi \widehat{ g_3 \, \Phi} \, v_h^- (x, y,
w_{k+ \ot}, \mu, \ph)dx \, dy \, d\mu \, d\ph ,
$$
$$
F_w^-=\ix \iy \imu \iphi \widehat{ g_3 \, \Phi} \, v_h^+  (x, y,
w_{k- \ot},\mu, \ph)dx \, dy \, d\mu \, d\ph ,
$$

$$
F_\mu^+=\ix \iy \iw \iphi \widetilde{g_4 \, \Phi} \, v_h^- (x, y, w,
\mu_{m+ \ot},  \ph)dx \, dy \, dw \, d\ph ,
$$
$$
F_\mu^-=\ix \iy \iw \iphi \widetilde{g_4 \, \Phi} \, v_h^+ ( x, y,
w,\mu_{m- \ot},  \ph)dx \, dy \, dw \, d\ph ,
$$

$$
F_\ph^+=\ix \iy \iw \imu g_5 \, \dot{\Phi} \, v_h^- (x, y, w, \mu,
\ph_{n+ \ot})dx \, dy \, dw \, d\mu ,
$$
$$
F_\ph^-=\ix \iy \iw \imu g_5 \, \dot{\Phi} \, v_h^+ (x, y, w, \mu,
\ph_{n- \ot})dx \, dy \, dw \, d\mu .
$$
where the upwind numerical fluxes $\check{\Phi}, \bar{\Phi},
\widehat{ g_3 \, \Phi} , \widetilde{g_4 \, \Phi}, \dot{\Phi}$  are
defined in the following way,
\begin{itemize}
\item The sign of $g_1$ only depends on $\mu$,
if $\mu_m >0$, then $ \check{\Phi} =\Phi^-$; otherwise, $
\check{\Phi} =\Phi^+.$

\item The sign of $g_2$ only depends on $\cos \ph $,
if $\cos \ph_n  >0$, then $ \bar{\Phi} =\Phi^- $; otherwise, $
\bar{\Phi} =\Phi^+ .$ Note that in our simulation, $N_\ph$ is always
even.

\item For $\widehat{ g_3 \, \Phi}$, we   let $$\widehat{ g_3 \, \Phi}
=- 2 c_k \frac{\sqrt{w(1+\ak w)}}{1+2 \ak w} \left[  \mu \,
E_x(t,x,y) \hat{\Phi} +
 \sqrt{1-\mu^2} \cos\ph \, E_y(t,x,y) \tilde{\Phi} \right] ,
$$
If $\mu_m E_x(t, x_i, y_j)<0$, then $\hat{\Phi}=\Phi^-$; otherwise,
$ \hat{\Phi} =\Phi^+.$

If $(\cos \ph_n ) E_y(t, x_i, y_j)<0$, then $\tilde{\Phi}=\Phi^-$;
otherwise, $ \tilde{\Phi} =\Phi^+.$

\item For $\widetilde{ g_4 \, \Phi}$, we  let $$\widetilde{ g_4 \, \Phi}
=-  c_k \frac{\sqrt{1-\mu^2}}{\sqrt{w(1+\ak w)}} \left[
\sqrt{1-\mu^2} \, E_x(t,x,y) \hat{\Phi} -
 \mu \cos\ph \, E_y(t,x,y) \tilde{\Phi} \right] ,
$$
If $ E_x(t, x_i, y_j)<0$, then $\hat{\Phi}=\Phi^-$; otherwise, $
\hat{\Phi} =\Phi^+$.

If $\mu_m \cos(\ph_n) E_y(t, x_i, y_j)>0$, then
$\tilde{\Phi}=\Phi^-$; otherwise, $ \tilde{\Phi} =\Phi^+$.

\item The sign of $g_5$ only depends on $E_y(t, x, y)$,
if $E_y(t, x_i, y_j)>0$, then $ \dot{\Phi} =\Phi^- $; otherwise, $
\dot{\Phi} =\Phi^+ .$
\end{itemize}

\begin{figure}[htb]
\centering
\includegraphics[width=0.95\linewidth]{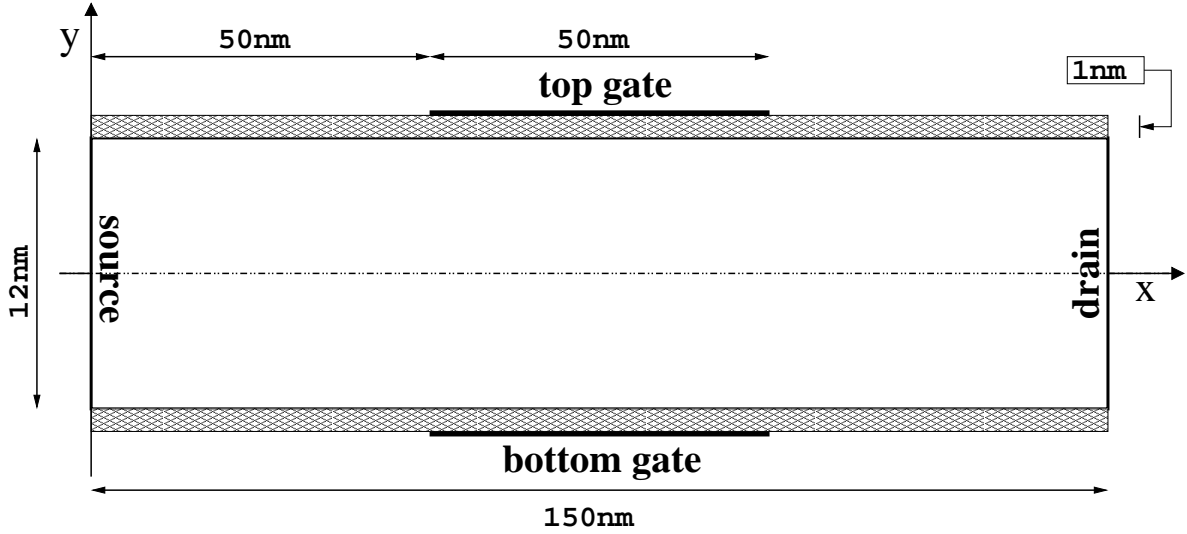}
\caption{Schematic representation of a 2D double gate MOSFET device}
\label{mosfet}
\end{figure}

The schematic plot of the double gate MOSFET device is given in
Figure \ref{mosfet}. The shadowed region denotes the oxide-silicon
region, whereas the rest is the silicon region. Since the problem is
symmetric about the x-axis, we will only need to compute for $y>0$.
At the source and drain contacts, we implement the same boundary
condition as proposed in \cite{cgms06} to realize neutral charges. A
buffer layer of ghost points of $i=0$ and $i=N_x+1$ is used to make
$$\Phi(i=0)=\Phi(i=1)\frac{N_D(i=1)}{\rho(i=1)},$$ and
$$\Phi(i=N_x+1)=\Phi(i=N_x)\frac{N_D(i=N_x)}{\rho(i=N_x)}.$$

At the top and bottom of the computational  domain (the silicon
region), we   impose the classical elastic specular boundary
reflection.

In the $(w,\mu,\ph)$-space, no boundary condition is necessary, the
reason is similar as in 1D,
\begin{itemize}
\item at $w=0$, $g_3=0$. At $w=w_{\mbox{max}}$, $\Phi$ is machine zero;
\item at $\mu=\pm 1$, $g_4=0$;
\item at $\ph=0, \pi$, $g_5=0$,
\end{itemize}
so at the $w, \mu, \ph$ boundary, the numerical flux vanishes, hence
no ghost point is necessary.

For the Poisson equation, $\nV=0.52354$ at source, $\nV=1.5235$ at
drain and $\nV=1.06$ at gate. For the rest of  boundaries, we impose
homogeneous Neumann boundary condition, i.e., $\frac{\partial
\nV}{\partial n}=0$. The relative dielectric constant in the
oxide-silicon region is $\epsilon_r=3.9$, in the silicon region is
$\epsilon_r=11.7$.

The Poisson equation (\ref{pois}) is solved by the LDG method. It
involves rewriting the equation into the following form,
\begin{equation}
\label{pois2} \left\{\begin{array} {l}
        \displaystyle    q= \frac{\partial \nV}{\partial x} , \qquad  s=\frac{\partial \nV}{\partial y}  \\
    \displaystyle     \frac{\partial}{\partial x} \left( \epsilon_{r} q \right)
 + \frac{\partial}{\partial y} \left( \epsilon_{r} s \right)
 = R(t,x,y)
         \end{array}
   \right.
\end{equation}
where $ R(t,x,y)=c_{p} \left[ \rho(t,x,y) - \mathcal{N}_{D}(x,y)
\right]$ is a known function that can be computed at each time step
once $\Phi$ is solved from (\ref{dgb}), and the coefficient
$\epsilon_r$ depends on $x, y$. The Poisson system is only  on the
$(x,y)$ domain. Hence, we use the  grid $I_{ij}=\left[ x_{i - \ot} ,
\, x_{i + \ot} \right] \times
                  \left[ y_{j - \ot} , \, y_{j + \ot} \right] $, with $i=1,\ldots, N_x$, $j=1, \ldots, N_y+M_y$, that includes the oxide-silicon region and is consistent with the five-dimensional rectangular grid for the Boltzmann equation in the silicon region. The approximation space  is defined as
\begin{equation}
W_h^\kpol=\{ v : v|_{I_{ij}} \in P^\kpol(I_{ij})\}.
\end{equation}
Here $P^\kpol(I_{ij})$ denotes the set of all polynomials of degree
at most $\kpol$ on $I_{ij}$.  The LDG scheme for (\ref{pois2}) is:
to find $q_h, s_h, \nV_h \in V_h^\kpol$, such that
\begin{eqnarray}
\label{ldgpois} && \mbox{} \hspace{-22pt} \int_{I_{i,j}} q_h v_h
dxdy + \int_{I_{i,j}} \nV_h (v_h)_x dxdy -\int_{y_{j - \ot}}^{y_{j +
\ot}} \hat{\nV}_h v_h^-( x_{i + \ot}, y) dy
+\int_{y_{j - \ot}}^{y_{j + \ot}} \hat{\nV}_h v_h^+( x_{i - \ot}, y) dy =0, \nonumber \\
&& \mbox{} \hspace{-22pt} \int_{I_{i,j}} \! s_h w_h dxdy +
\int_{I_{i,j}} \! \nV_h (w_h)_y dxdy -\int_{x_{i - \ot}}^{x_{i +
\ot}} \tilde{\nV}_h w_h^-(x, y_{j + \ot}) dx
+\int_{x_{i - \ot}}^{x_{i + \ot}} \tilde{\nV}_h w_h^+(x, y_{j - \ot}) dx =0,  \nonumber \\
&& \mbox{} -\int_{I_{i,j}} \epsilon_{r} q_h (p_h)_x dxdy +\int_{y_{j
- \ot}}^{y_{j + \ot}} \widehat{ \epsilon_{r} q}_h p_h^-( x_{i +
\ot}, y) dy -\int_{y_{j - \ot}}^{y_{j + \ot}} \widehat{\epsilon_{r}
q}_h p_h^+( x_{i - \ot}, y) dy
\nonumber \\
&& \mbox{} -\int_{I_{i,j}} \epsilon_{r} s_h (p_h)_y dxdy +
\int_{x_{i - \ot}}^{x_{i + \ot}} \widetilde{ \epsilon_{r}  s}_h
p_h^-(x, y_{j + \ot}) dx -\int_{x_{i - \ot}}^{x_{i + \ot}}
\widetilde{ \epsilon_{r}  s}_h p_h^+(x, y_{j - \ot}) dx
\nonumber \\
& & \mbox{} = \int_{I_{i,j}} R(t,x,y) p_h dxdy
\end{eqnarray}
hold true for any $v_h, w_h, p_h \in W_h^\kpol$. In the above
formulation, we choose the flux as follows,
 in the $x$-direction, we use $\hat{\nV}_h=\nV^-_h$, $\widehat{ \epsilon_{r} q}_h=\epsilon_{r} q_h^+ -[\nV_h]$. In the $y$-direction, we use $\tilde{\nV}_h=\nV^-_h$, $\widetilde{ \epsilon_{r}  s}_h = \epsilon_{r} s_h^+ -[\nV_h]$.
 Near the drain, we are given Dirichlet boundary condition, so we need to flip
the flux in $x-$direction: let $\hat{\nV}_h ( x_{i + \ot},
y)=\nV^+_h ( x_{i + \ot}, y)$ and $\widehat{ \epsilon_{r} q}_h (
x_{i + \ot}, y)=\epsilon_{r} q_h^-( x_{i + \ot}, y) -[\nV_h]( x_{i +
\ot}, y),$ if the point $( x_{i + \ot}, y)$ is at the drain. For the
gate, we need to flip the flux in $y-$direction: let
$\tilde{\nV}_h(x,y_{j + \ot})=\nV^+_h(x,y_{j + \ot})$ and
$\widetilde{ \epsilon_{r}  s}_h(x,y_{j + \ot}) = \epsilon_{r}
s_h^-(x,y_{j + \ot}) -[\nV_h](x,y_{j + \ot})$, if  the point
$(x,y_{j + \ot})$ is at the gate.  For the bottom, we need to use
the Neumann condition, and flip the flux in y-direction, i.e.,
$\tilde{\nV}_h=\nV^+_h$, $\widetilde{ \epsilon_{r}  s}_h =
\epsilon_{r} s_h^-$.  This scheme described above will enforce the
continuity of $\nV$ and $\epsilon_r \frac{\partial \nV}{\partial n}$
across the interface of silicon and oxide-silicon interface.  The
solution of (\ref{ldgpois}) gives us approximations to both the
potential $\nV_h$ and  the electric field $(E_x)_h=-c_v q_h$,
$(E_y)_h=-c_v s_h$.

To summarize, start with an initial condition for $\Phi_h$, the
DG-LDG algorithm for the 2D double gate MOSFET advances from $t^n$
to $t^{n+1}$ in the following steps:
\begin{description}
\item[Step 1] Compute $\rho_h(t,x,y)= \Iwmp \Phi_h (t,x,y,w,\mu,\ph)$.
\item[Step 2] Use $\rho_h(t,x,y)$ to solve from (\ref{ldgpois}) the electric field $(E_x)_h$ and $(E_y)_h$, and compute $g_i$, $i=1, \ldots, 5$.
\item[Step 3] Solve (\ref{dgb}) and get a method of line ODE for $\Phi_h$.
\item[Step 4] Evolve this ODE by proper time stepping from $t^n$ to $t^{n+1}$, if partial time step is necessary, then repeat Step 1 to 3 as needed.
\end{description}

\begin{figure}[ht]
\centering
\includegraphics[width=3in,angle=0]{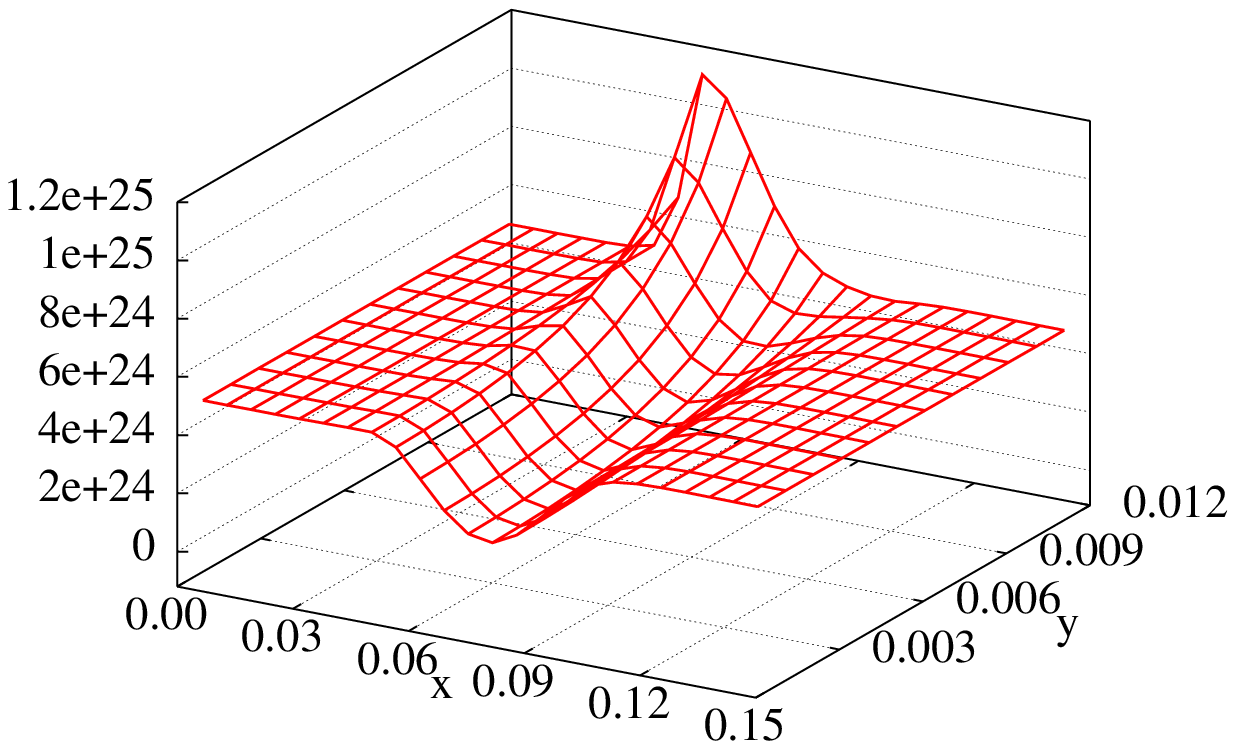}
\includegraphics[width=3in,angle=0]{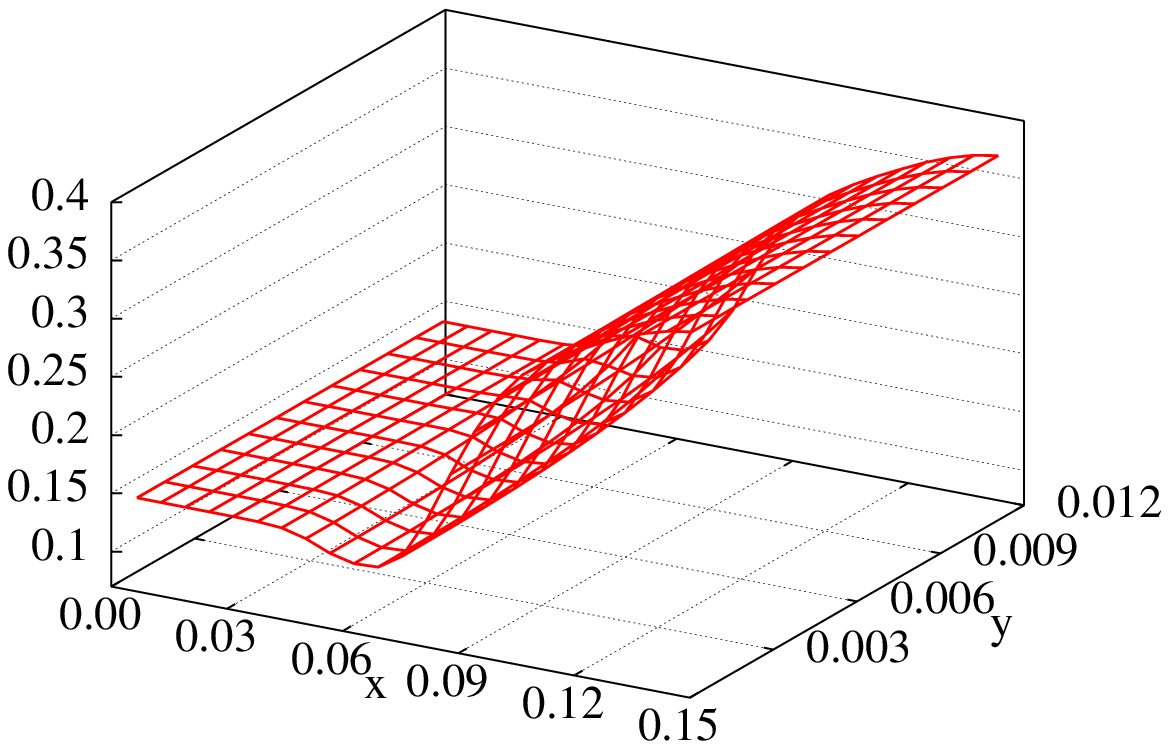}\\
\includegraphics[width=3in,angle=0]{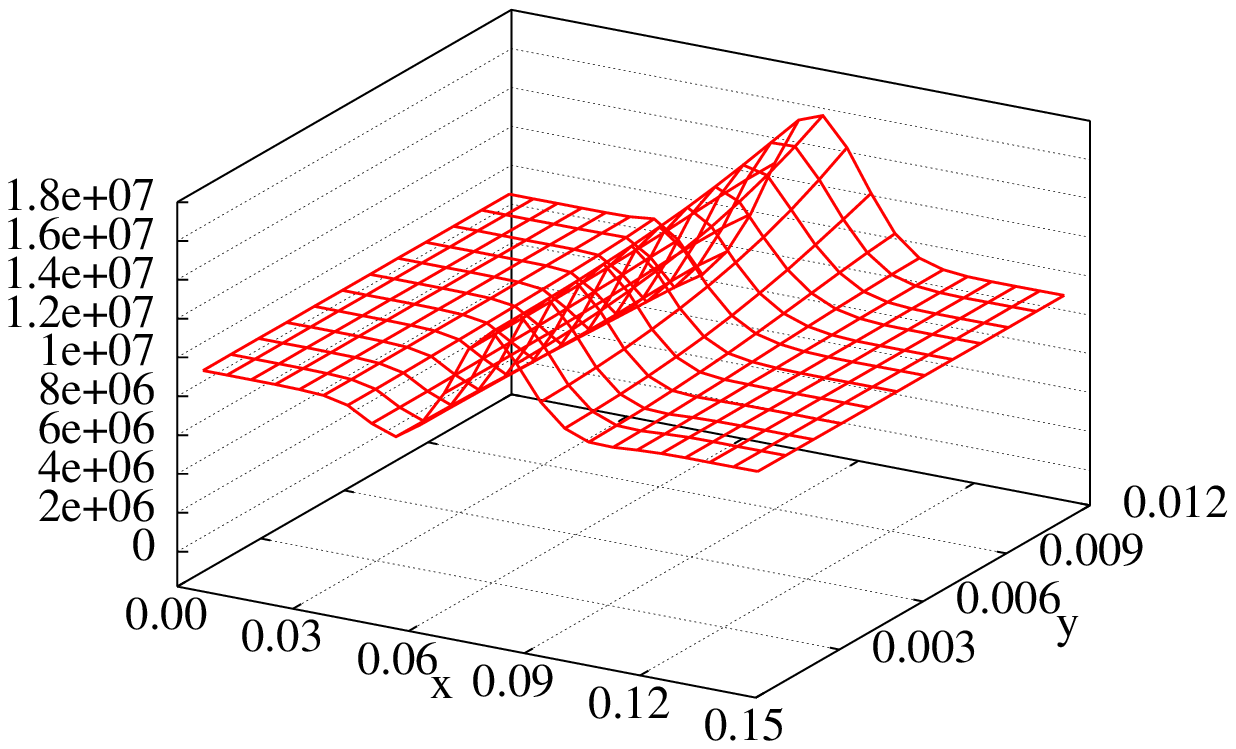}
\includegraphics[width=3in,angle=0]{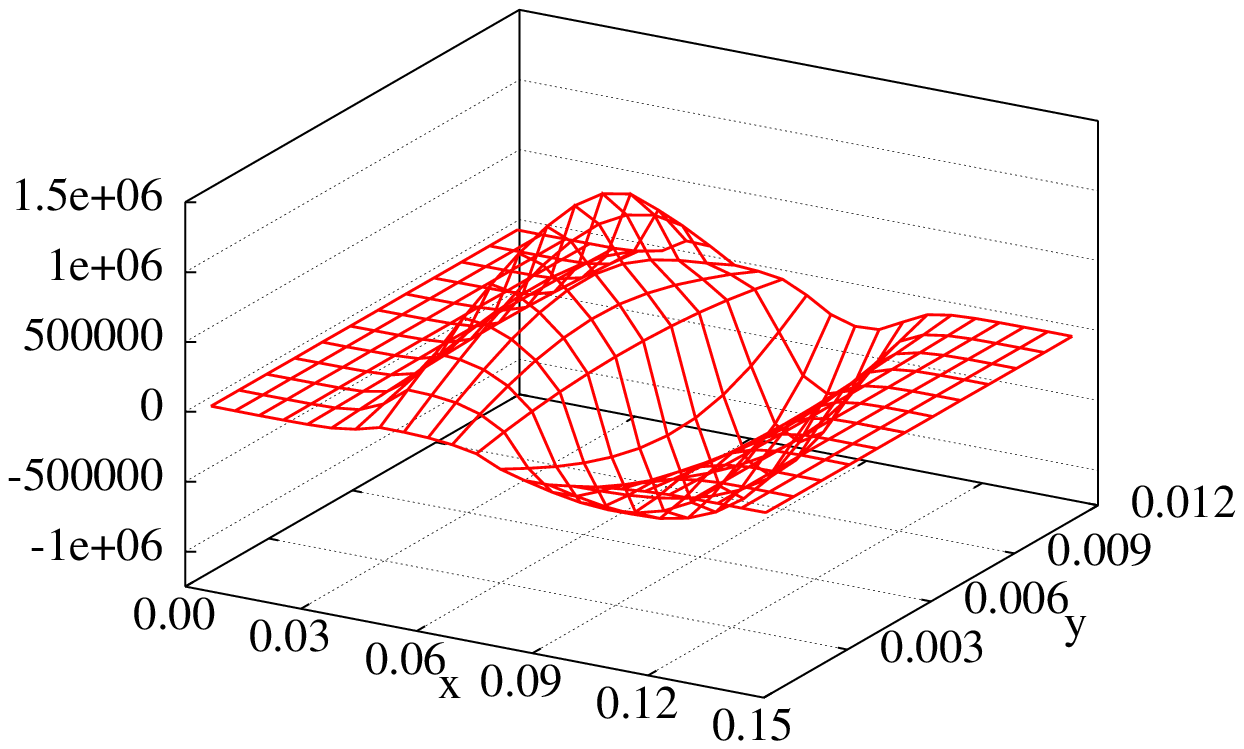}
\caption{Macroscopic quantities of double gate MOSFET device at
$t=0.5$. Top left: density in ${cm}^{-3}$; top right: energy in
$eV$;  bottom left: x-component of velocity in $cm/s$; bottom right:
y-component of velocity in $cm/s$. Solution reached steady state.}
\label{transhyd}
\end{figure}

\begin{figure}[htb]
\centering
\includegraphics[width=3in,angle=0]{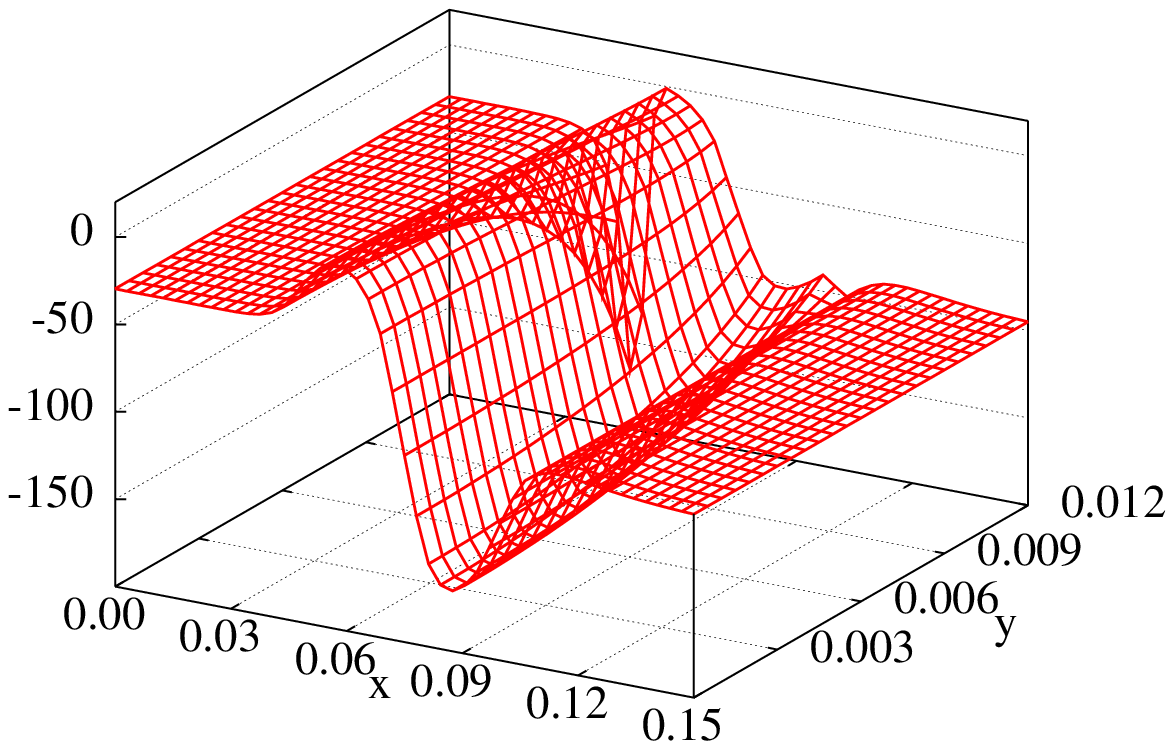}
\includegraphics[width=3in,angle=0]{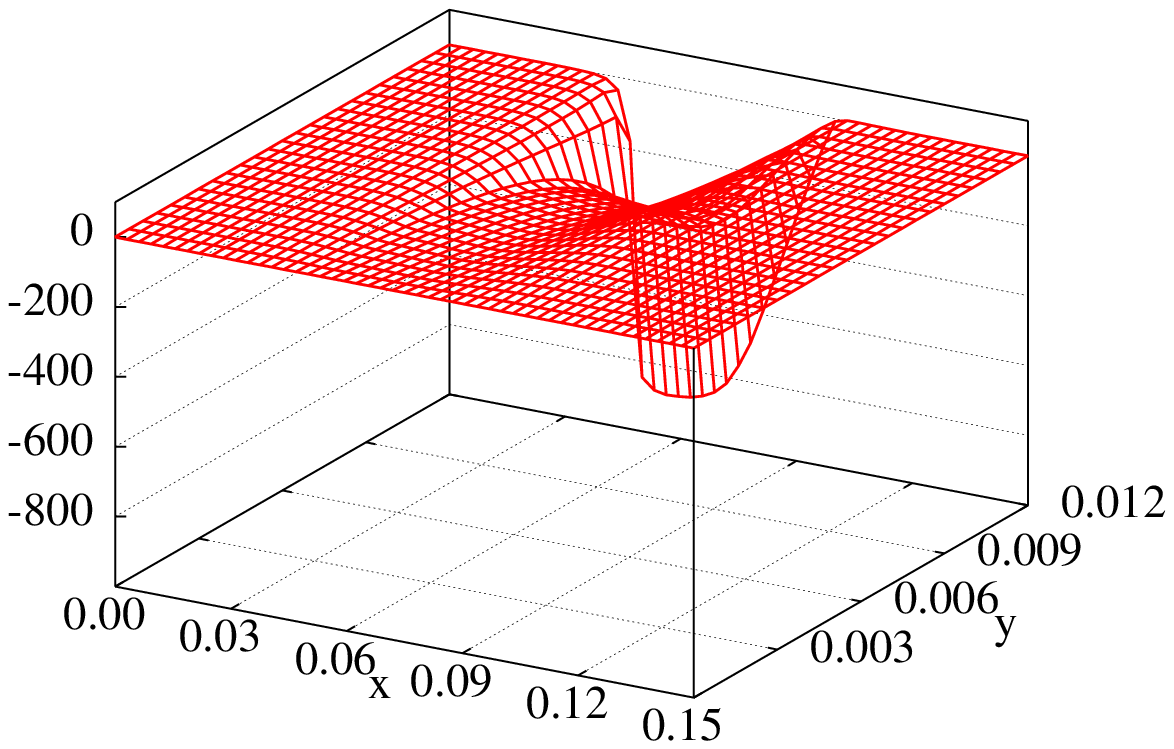}\\
\includegraphics[width=3in,angle=0]{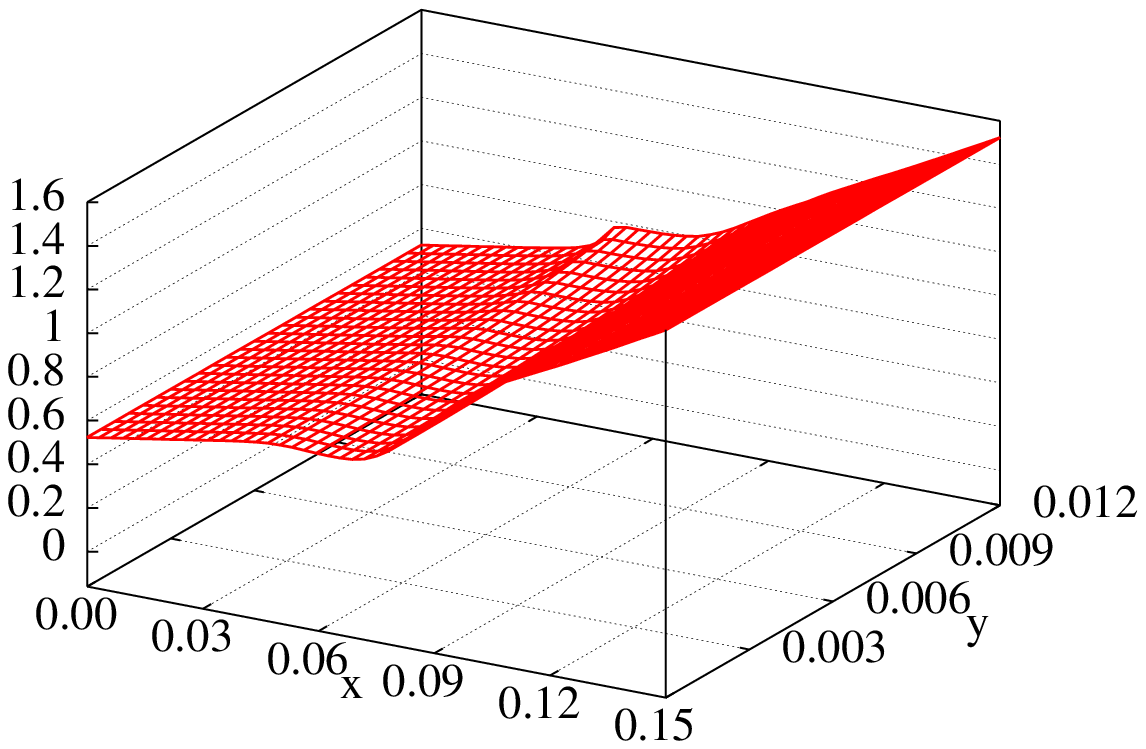}
\caption{Macroscopic quantities of double gate MOSFET device at
$t=0.5$. Top left: x-component of electric field in $kV/cm$; top
right: y-component of electric field in $kV/cm$;  bottom: electric
potential in $V$. Solution has reached steady state.}
\label{transelec}
\end{figure}

\begin{figure}[htb]
\centering
\includegraphics[width=2.8in,angle=0]{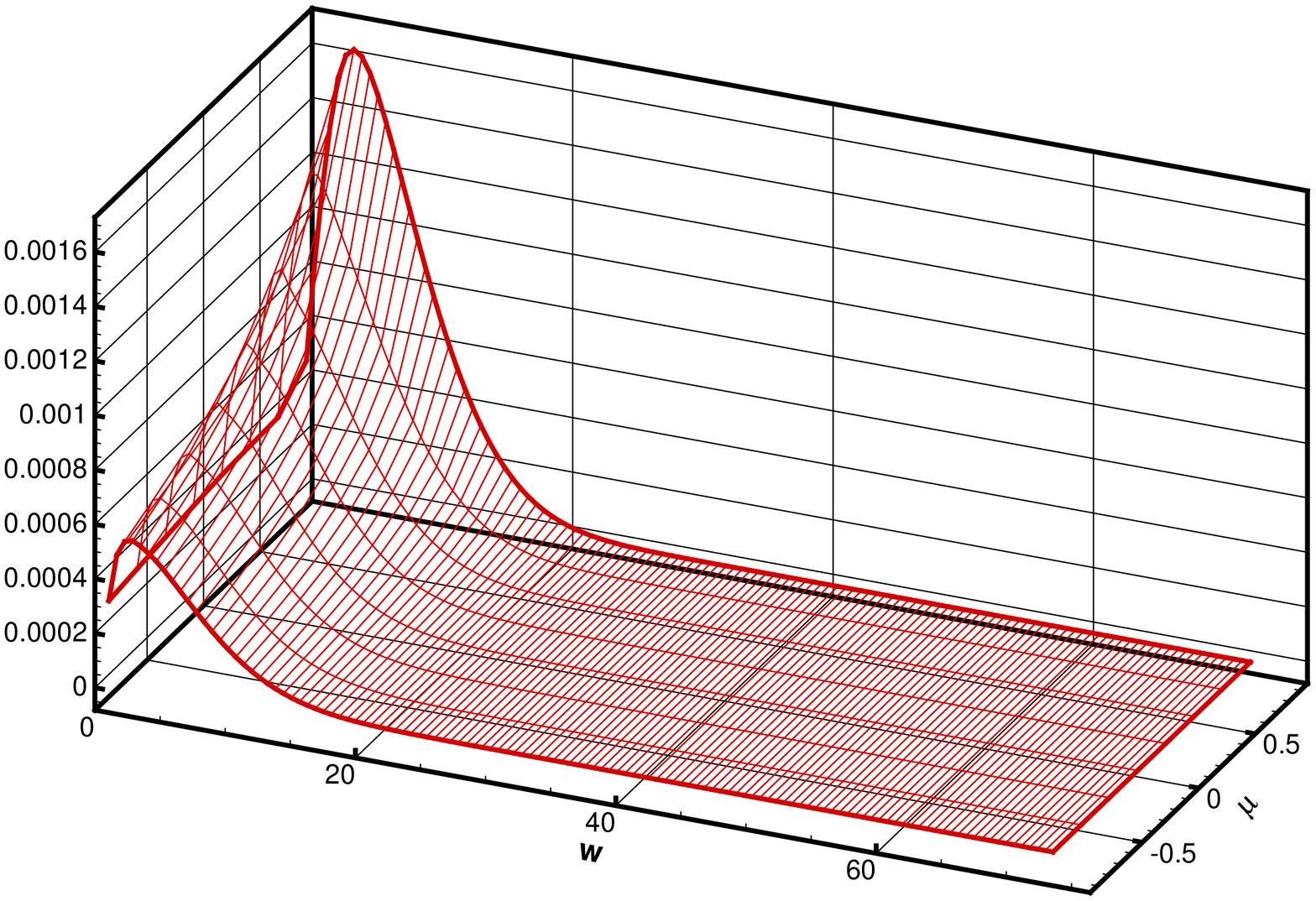}
\includegraphics[width=2.8in,angle=0]{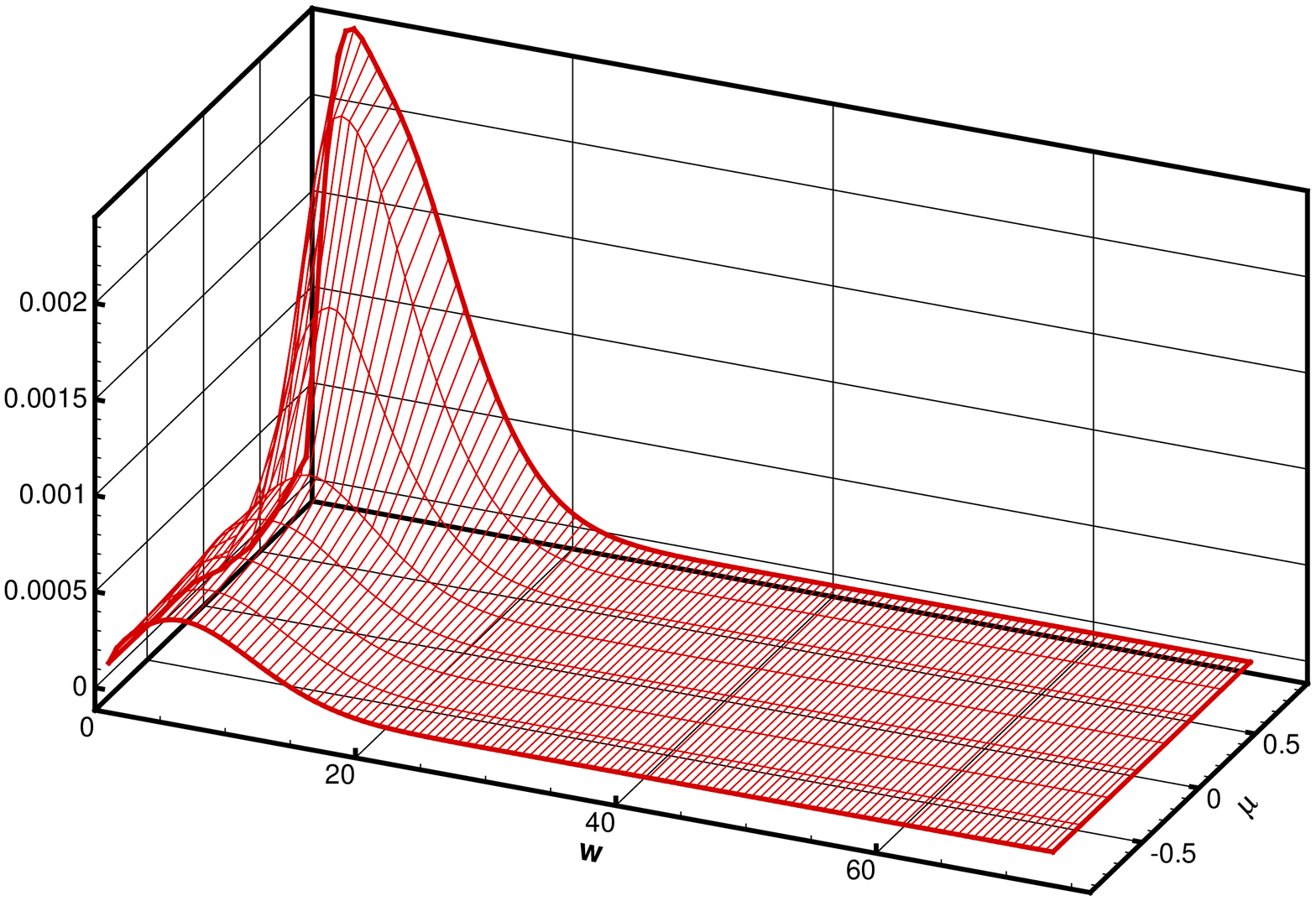}\\
\includegraphics[width=2.8in,angle=0]{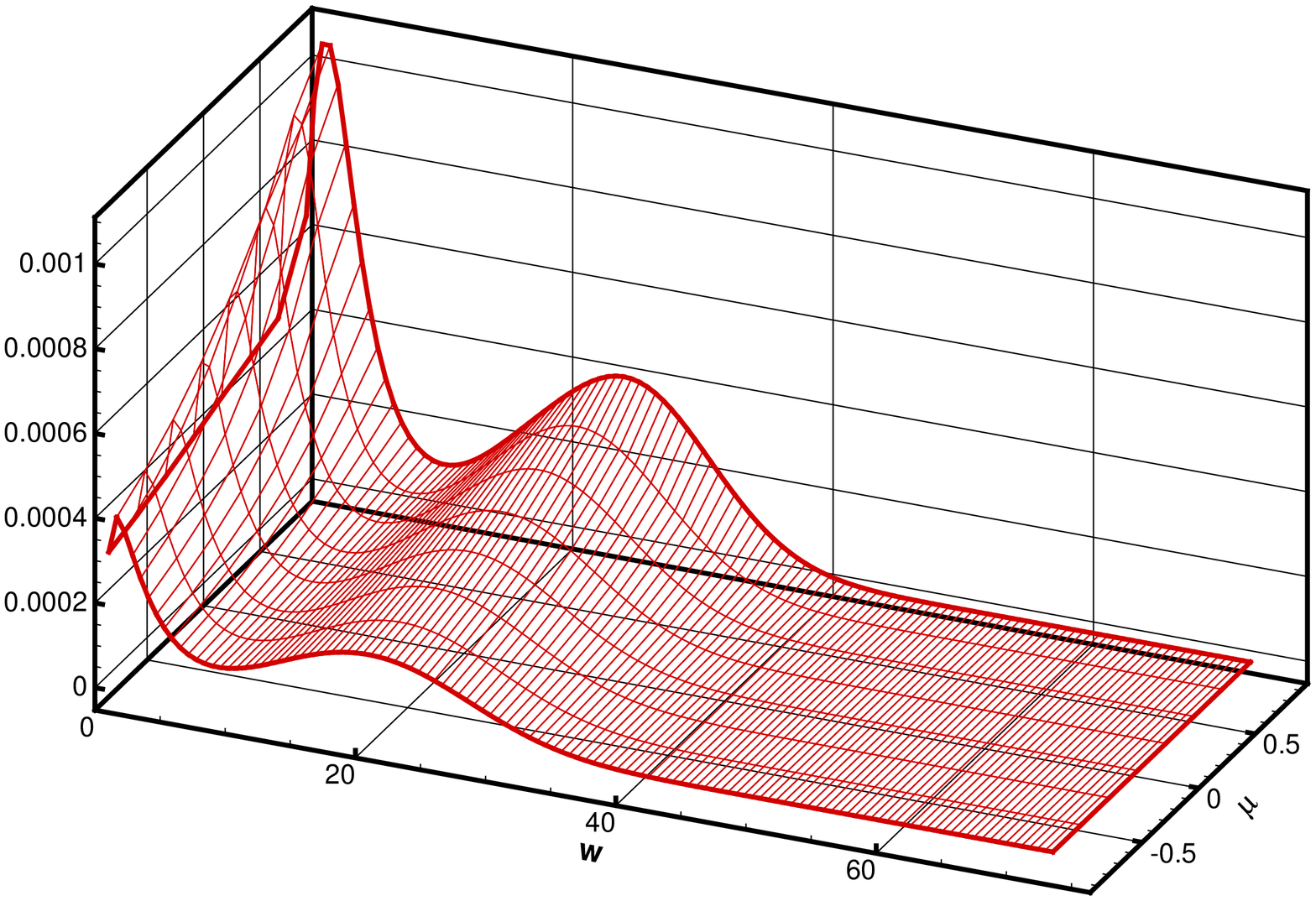}
\includegraphics[width=2.8in,angle=0]{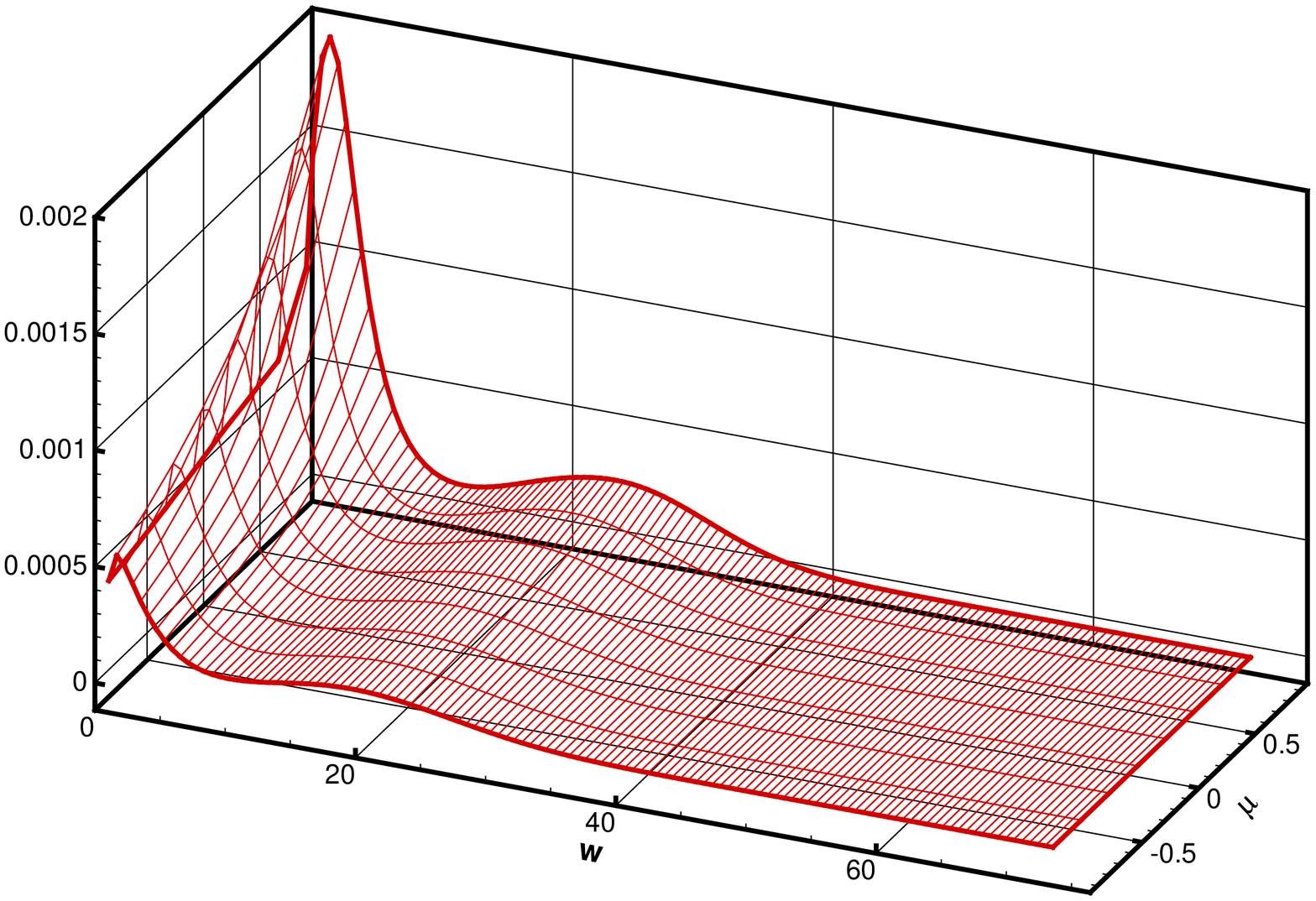}\\
\includegraphics[width=2.8in,angle=0]{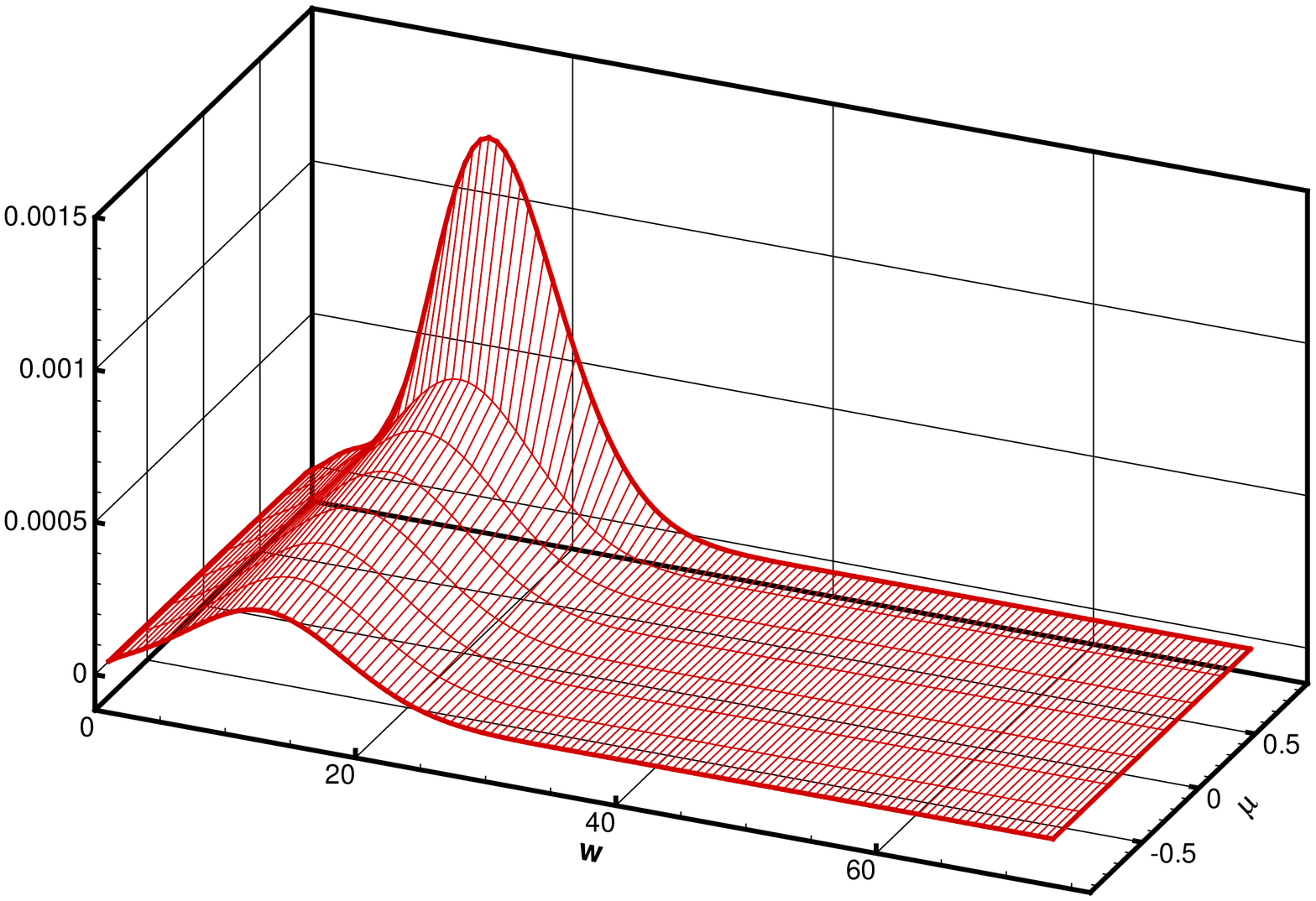}
\includegraphics[width=2.8in,angle=0]{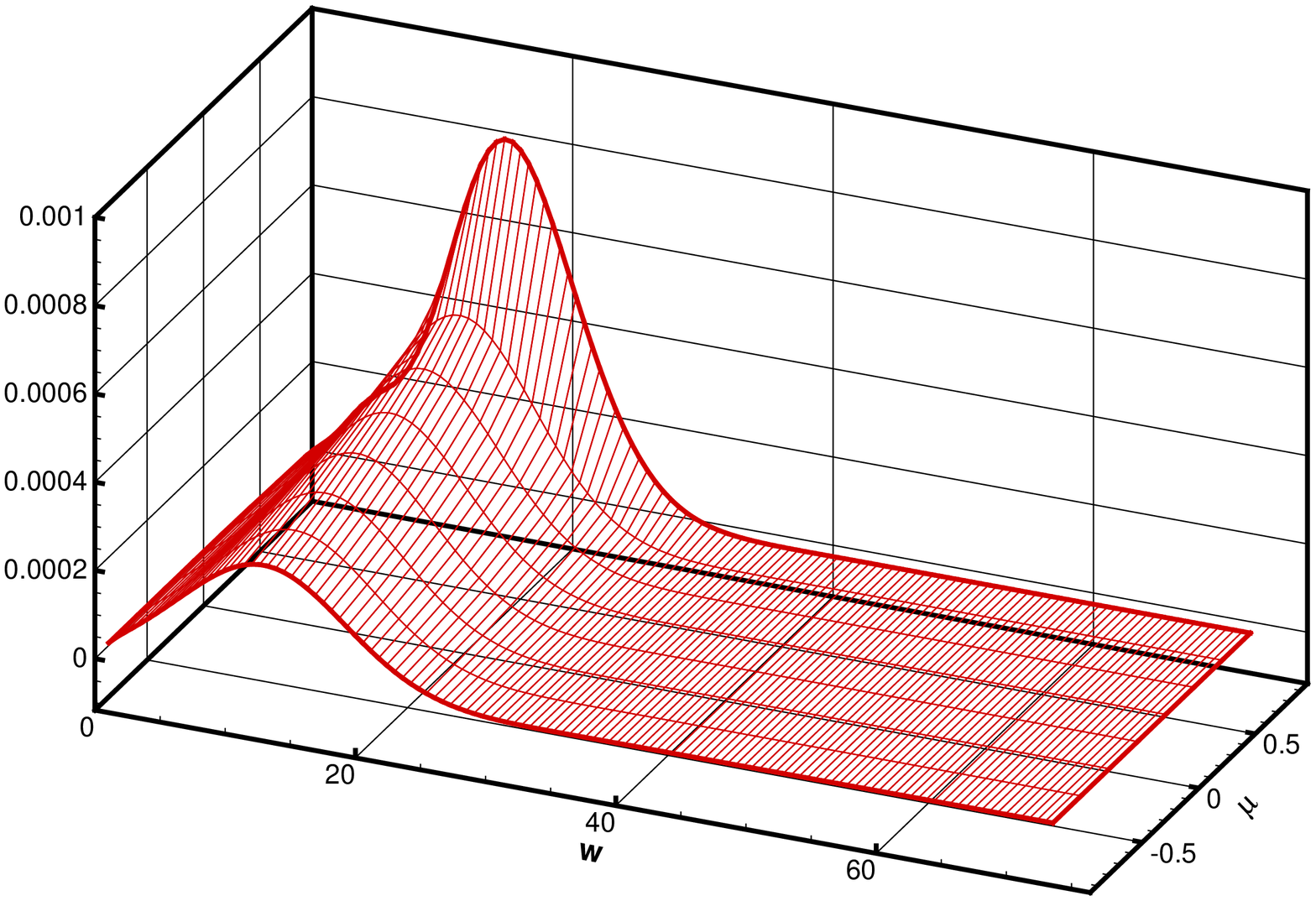}
\caption{PDF of double gate MOSFET device at $t=0.5$. Top left: at
(0.025, 0.12); top right: at at (0.075, 0.12);  middle left: at
(0.125, 0.12) ;  middle right: at (0.1375, 0.06);  bottom left: at
(0.09375, 0.10) ;  bottom right: at (0.09375, 0.)   . Solution
reached steady state.} \label{transpdf}
\end{figure}

\begin{figure}[htb]
\centering
\includegraphics[width=4in,angle=0]{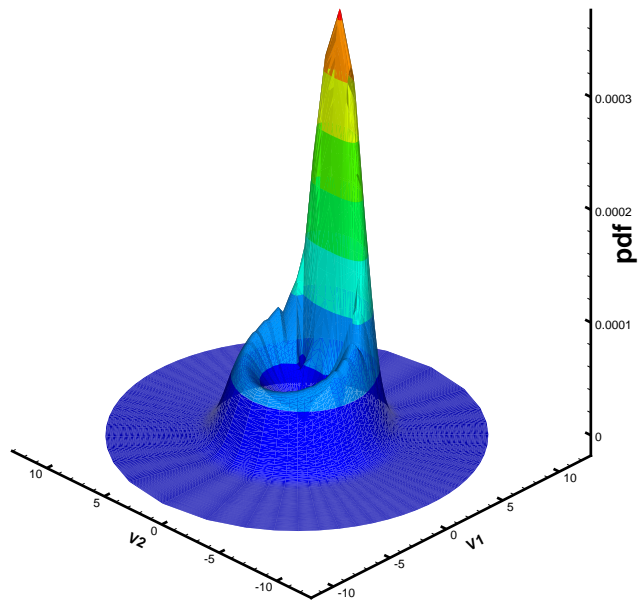}
\caption{PDF for 2D double gate MOSFET at $t=0.5$, $(x,y)=(0.9375,
0.10)$.} \label{2dcart}
\end{figure}

All numerical results  are obtained with a piecewise linear
approximation space and first order Euler time stepping.  Apparently
the collision term makes the Euler forward time stepping stable. We
use a $24 \times 14$ grid in space, $120$ points in $w$, $8$ points
in $\mu$ and  $6$ points in $\ph$. In Figures \ref{transhyd} and
\ref{transelec}, we show the results of the macroscopic quantities.
We also show the \emph{pdf} at six different locations in the device
in Figure \ref{transpdf}. These \emph{pdf}'s have been computed by
averaging the values of $\Phi_h$ over $\ph$. In Figure \ref{2dcart},
we present the cartesian plot for \emph{pdf} at $(x,y)=(0.125,
0.12)$, where a very non-equilibrium \emph{pdf} is observed.
%
%

\clearpage
\section{Conclusions and final remarks}

We have developed a DG scheme for  BTEs of type (\ref{bte}), which
takes into account optical-phonon interactions that become dominant
under strong energetic conditions.   We used the coordinate
transformation proposed in \cite{MP, cgms03} and changed the
collision into an integral-difference operator by using energy band
as one of the variables. The Poisson equation is treated by LDG   on
a  mesh that is consistent with the mesh of the DG-BTE scheme. The
results are compared  to those obtained from a high order WENO
scheme simulation. By a local refinement in mesh,  we were able to
capture the subtle kinetic effects including very non-equilibrium
distributions without a great increase of memory allocation and CPU
time. The advantage of the DG scheme lies in its potential for
implementation on unstructured meshes and  for full {\it
hp-}adaptivity. The simple communication pattern of the DG method
also makes it a good candidate for the domain decomposition method
for the coupled kinetic and macroscopic models.

%

%

%
\clearpage
\appendix
\setcounter{equation}{0}
\section{Appendix}

In this appendix, we collect some technical details for the
implementation of the 2D DG-BTE solver. The discussion for 1D solver
is similar and omitted here.

\subsection{The basis of the finite dimensional function space.}

In  every cell $\Omega_{ijkmn}$, we use piecewise linear polynomials
 and assume
\begin{eqnarray}
&&  \Phi_h (t,x,y,w,\mu,\ph) =
  T_{ijkmn}(t) +
  X_{ijkmn}(t) \, \frac{2(x - x_{i})}{\Delta x_{i}} +
  Y_{ijkmn}(t) \, \frac{2(y - y_{j})}{\Delta y_{j}}
\nonumber \\
& & \mbox{}  +
  W_{ijkmn}(t) \, \frac{2(w - w_{k})}{\Delta w_{k}} +
  M_{ijkmn}(t) \, \frac{2(\mu - \mu_{m})}{\Delta \mu_{m}} +
  P_{ijkmn}(t) \, \frac{2(\ph - \ph_{n})}{\Delta \ph_{n}}
  \, .
\end{eqnarray}
It will be useful to note that
\begin{eqnarray*}
&& \Phi_h (t,x,y,w,\mu,\ph) = \sum_{k=1}^{N_{w}} \left[
  T_{ijkmn}(t) +
  X_{ijkmn}(t) \, \frac{2(x - x_{i})}{\Delta x_{i}} +
  Y_{ijkmn}(t) \, \frac{2(y - y_{j})}{\Delta y_{j}} \right.
  \quad \quad \mbox{ } \\
& & \left. \mbox{} \: +
  W_{ijkmn}(t) \, \frac{2(w - w_{k})}{\Delta w_{k}} +
  M_{ijkmn}(t) \, \frac{2(\mu - \mu_{m})}{\Delta \mu_{m}} +
  P_{ijkmn}(t) \, \frac{2(\ph - \ph_{n})}{\Delta \ph_{n}}
  \right] \chi_{k}(w) \, ,
\end{eqnarray*}
for every $ \dm (x,y,w,\mu,\ph) \in \bigcup_{k=1}^{N_{w}}
\Omega_{ijkmn}$. Here, $\chi_{k}(w)$ is the characteristic function
in the interval $ \dm \left[ w_{k - \ot} , \, w_{k + \ot} \right] $.

\bigskip
\subsection{Treatment of the collision operator}

The gain term of the collisional operator is
\begin{eqnarray}
&&   G(\Phi_h)(t,x,y,w) = s(w) \left\{ c_{0} \Ipm \Phi_h(t,x,y,w,\mu
',\ph')
   \right. \nonumber \\
&& \left. +  \Ipm [ c_{+} \Phi_h(t,x,y,w + \qe,\mu ',\ph')
   + c_{-} \Phi_h(t,x,y,w - \qe,\mu ',\ph') ] \right\} .
\end{eqnarray}
Now, we define
$$
 (\overline{v_h})_{m n}(x,y,w) := \int_{\ph_{n - \ot}}^{\ph_{n + \ot}} d \ph
\int_{\mu_{m - \ot}}^{\mu_{m + \ot}} d \mu \: v_h(x,y,w,\mu, \ph) \,
,
$$
and, for $\sigma = - \qe$, $0$, $\qe$, we have
\begin{eqnarray*}
&& \int_{\Omega_{ijkmn}} v_h(x,y,w,\mu, \ph) \left[ s(w) \Ipm
\Phi(t,x,y,w + \sigma,\mu ',\ph') \right]
dx \, dy \, dw \, d\mu \, d\ph   \\[7pt]
&& \mbox{} =
 \Ipm \hspace{-4pt}
\int_{x_{i - \ot}}^{x_{i + \ot}} dx \hspace{-1pt} \int_{y_{j -
\ot}}^{y_{j + \ot}} dy \hspace{-1pt} \int_{w_{k - \ot}}^{w_{k +
\ot}} dw \: s(w) \, \Phi(t,x,y,w + \sigma,\mu ',\ph') \,
(\overline{v_h})_{m n}(x,y,w)
\end{eqnarray*}
\begin{center}

$ \dm
 \mbox{} =
\sum_{m'=1}^{N_{\mu}} \sum_{n'=1}^{N_{\ph}} \int_{\Omega_{ijkm'n'}}
s(w) \, \Phi_h(t,x,y,w + \sigma,\mu ',\ph') \, (\overline{v_h})_{m
n}(x,y,w) \: dx \, dy \, dw \, d\mu' \, d\ph' \, . $

\end{center}

Now we discuss the following integral for different test function
$v_h$,
$$
I=\dm \int_{\Omega_{ijkmn}} v_h(x,y,w,\mu, \ph) \left[ s(w) \Ipm
\Phi_h(t,x,y,w + \sigma,\mu ',\ph') \right] dx \, dy \, dw \, d\mu
\, d\ph \, .
$$
%

%
\bigskip
\begin{itemize}
\item For $v_h(x,y,w,\mu, \ph) = 1$,
\begin{eqnarray*}
&&I= \sum_{m'=1}^{N_{\mu}} \sum_{n'=1}^{N_{\ph}} \sum_{k'=1}^{N_{w}}
 \Delta \mu_{m'} \, \Delta \ph_{n'}
\left[ T_{ijk'm'n'}(t) \int_{w_{k - \ot}}^{w_{k + \ot}} s(w) \,
\chi_{k'}(w + \sigma) \, dw
\right. \\
&& \left. \mbox{} + W_{ijk'm'n'}(t)
 \int_{w_{k - \ot}}^{w_{k + \ot}} s(w)
 \, \frac{2(w + \sigma - w_{k'})}{\Delta w_{k'}} \,  \chi_{k'}(w + \sigma)  \, dw
\right]
 \Delta x_{i} \, \Delta y_{j} \, \Delta \mu_{m} \, \Delta \ph_{n}  \, .
\\[10pt]
\end{eqnarray*}
\item For $
v_h(x,y,w,\mu, \ph) = \frac{2(x - x_{i})}{\Delta x_{i}} $,
\begin{eqnarray*}
&& I= \frac{1}{3}
 \Delta x_{i} \, \Delta y_{j} \, \Delta \mu_{m} \, \Delta \ph_{n}
\sum_{m'=1}^{N_{\mu}} \sum_{n'=1}^{N_{\ph}} \sum_{k'=1}^{N_{w}}
 \Delta \mu_{m'} \, \Delta \ph_{n'} \, X_{ijk'm'n'}(t)
\\ && \mbox{} \qquad \times
\int_{w_{k - \ot}}^{w_{k + \ot}} s(w) \,  \chi_{k'}(w + \sigma) \,
dw \, .
\\[10pt]
\end{eqnarray*}

\item For $
v_h(x,y,w,\mu, \ph) = \frac{2(y - y_{j})}{\Delta y_{j}}$,
\begin{eqnarray*}
&& I= \frac{1}{3}
 \Delta x_{i} \, \Delta y_{j} \, \Delta \mu_{m} \, \Delta \ph_{n}
\sum_{m'=1}^{N_{\mu}} \sum_{n'=1}^{N_{\ph}} \sum_{k'=1}^{N_{w}}
 \Delta \mu_{m'} \, \Delta \ph_{n'} \, Y_{ijk'm'n'}(t)
\\ && \mbox{} \qquad \times
\int_{w_{k - \ot}}^{w_{k + \ot}} s(w) \,  \chi_{k'}(w + \sigma) \,
dw \, .
\\[10pt]
\end{eqnarray*}

\item For $
v_h(x,y,w,\mu, \ph) = \frac{2(w - w_{k})}{\Delta w_{k}}$,
\begin{eqnarray*}
&& I= \sum_{m'=1}^{N_{\mu}} \sum_{n'=1}^{N_{\ph}}
\sum_{k'=1}^{N_{w}}
 \Delta \mu_{m'} \, \Delta \ph_{n'} \left[ T_{ijk'm'n'}(t)
 \int_{w_{k - \ot}}^{w_{k + \ot}} s(w) \,
 \frac{2(w - w_{k})}{\Delta w_{k}} \chi_{k'}(w + \sigma) \, dw
\right. \\
&& \left.  \mbox{} \hspace{-24pt} +  W_{ijk'm'n'}(t)
 \int_{w_{k - \ot}}^{w_{k + \ot}} s(w) \,
 \frac{4 (w + \sigma - w_{k'}) (w - w_{k})}{\Delta w_{k'} \, \Delta w_{k}} \,
 \chi_{k'}(w + \sigma)  \, dw
\right]   \Delta x_{i} \, \Delta y_{j} \, \Delta \mu_{m} \, \Delta
\ph_{n} \, .
\\[10pt]
\end{eqnarray*}

\item For $v_h(x,y,w,\mu, \ph) = \frac{2(\mu - \mu_{m})}{\Delta \mu_{m}}$, $ I= 0$.
\item For $v_h(x,y,w,\mu, \ph) =  \frac{2(\ph - \ph_{n})}{\Delta \ph_{n}}$, $I = 0$.
\end{itemize}
The lost term in the collision operator is
\begin{equation}
2 \pi [c_0 s(w) +  c_+ s(w - \qe) + c_- s(w + \qe)]
   \Phi(t,x,y,w,\mu,\ph)  \, .
\end{equation}
Let
$$
\nu(w) = 2 \pi [c_0 s(w) +  c_+ s(w - \qe) + c_- s(w + \qe)] \, ,
$$
then we need to evaluate numerically,
$$I^\prime=\dm
\int_{\Omega_{ijkmn}} \nu(w) \, \Phi_h(t,x,y,w,\mu,\ph) \,
v_h(x,y,w,\mu, \ph) \: dx \, dy \, dw \, d\mu \, d\ph \, .
$$

\begin{itemize}%
\item For $v_h(x,y,w,\mu, \ph) = 1 $,
\begin{eqnarray*}
&& I^\prime= \Delta x_{i} \, \Delta y_{j} \, \Delta \mu_{m} \,
\Delta \ph_{n}
\\
&& \mbox{} \quad \times \left[ T_{ijkmn}(t) \int_{w_{k - \ot}}^{w_{k
+ \ot}} \nu(w) \, dw + W_{ijkmn}(t) \int_{w_{k - \ot}}^{w_{k + \ot}}
\nu(w) \,
 \frac{2(w - w_{k})}{\Delta w_{k}}  \, dw
\right] .
\\[10pt]
\end{eqnarray*}

\item For $v_h(x,y,w,\mu, \ph) = \frac{2(x - x_{i})}{\Delta x_{i}}$,
\begin{eqnarray*}
&& I^\prime= \frac{1}{3}\Delta x_{i} \, \Delta y_{j} \, \Delta
\mu_{m} \, \Delta \ph_{n} \, X_{ijkmn}(t) \int_{w_{k - \ot}}^{w_{k +
\ot}} \nu(w) \, dw \, .
\\[10pt]
\end{eqnarray*}

\item For $v_h(x,y,w,\mu, \ph) = \frac{2(y - y_{j})}{\Delta y_{j}}$,
\begin{eqnarray*}
&& I^\prime= \frac{1}{3}\Delta x_{i} \, \Delta y_{j} \, \Delta
\mu_{m} \, \Delta \ph_{n} \, Y_{ijkmn}(t) \int_{w_{k - \ot}}^{w_{k +
\ot}} \nu(w) \, dw \, .
\\[10pt]
\end{eqnarray*}

\item For $
v_h(x,y,w,\mu, \ph) = \frac{2(w - w_{k})}{\Delta w_{k}}$,
\begin{eqnarray*}
&& I^\prime= \Delta x_{i} \, \Delta y_{j} \, \Delta \mu_{m} \,
\Delta \ph_{n} \left[ T_{ijkmn}(t) \int_{w_{k - \ot}}^{w_{k + \ot}}
\, \nu(w) \frac{2(w - w_{k})}{\Delta w_{k}} \, dw \right.
\\
&& \mbox{} \quad  + \left. W_{ijkmn}(t) \int_{w_{k - \ot}}^{w_{k +
\ot}} \nu(w) \,\frac{4 (w - w_{k})^{2}}{\left( \Delta w_{k}
\right)^{2}} \,  dw \right] .
\\[10pt]
\end{eqnarray*}

\item For $
v_h(x,y,w,\mu, \ph) = \frac{2(\mu - \mu_{m})}{\Delta \mu_{m}}$,
\begin{eqnarray*}
&& I^\prime= \frac{1}{3}\Delta x_{i} \, \Delta y_{j} \, \Delta
\mu_{m} \, \Delta \ph_{n} \, M_{ijkmn}(t) \int_{w_{k - \ot}}^{w_{k +
\ot}} \nu(w) \, dw \, .
\\[10pt]
\end{eqnarray*}
%
%
\item For $
v_h(x,y,w,\mu, \ph) =  \frac{2(\ph - \ph_{n})}{\Delta \ph_{n}}$,
\begin{eqnarray*}
&& I^\prime= \frac{1}{3}\Delta x_{i} \, \Delta y_{j} \, \Delta
\mu_{m} \, \Delta \ph_{n} \, P_{ijkmn}(t) \int_{w_{k - \ot}}^{w_{k +
\ot}} \nu(w) \, dw \, .
\end{eqnarray*}
\end{itemize}

\subsection{Integrals related to the collisional operator}

We need to evaluate (some numerically) the following integrals
\begin{eqnarray*}
&& \int_{w_{k - \ot}}^{w_{k + \ot}} s(w) \,  \chi_{k'}(w + \sigma)
\, dw
\\
&&
 \int_{w_{k - \ot}}^{w_{k + \ot}} s(w)
 \, \frac{2(w + \sigma - w_{k'})}{\Delta w_{k'}} \,  \chi_{k'}(w + \sigma)  \, dw
\\
&&
 \int_{w_{k - \ot}}^{w_{k + \ot}} s(w)
 \frac{2 (w - w_{k})}{\Delta w_{k}} \,
 \chi_{k'}(w + \sigma)  \, dw
\\
&&
 \int_{w_{k - \ot}}^{w_{k + \ot}} s(w)
 \frac{4 (w + \sigma - w_{k'}) (w - w_{k})}{\Delta w_{k'} \, \Delta w_{k}} \,
 \chi_{k'}(w + \sigma)  \, dw
\\
&& \int_{w_{k - \ot}}^{w_{k + \ot}} \nu(w) \, dw
\\
&& \int_{w_{k - \ot}}^{w_{k + \ot}} \nu(w) \,
 \frac{2(w - w_{k})}{\Delta w_{k}}  \, dw
\\
&& \int_{w_{k - \ot}}^{w_{k + \ot}} \nu(w) \left[ \frac{2(w -
w_{k})}{\Delta w_{k}} \right]^{2} dw \, .
\end{eqnarray*}
If we  evaluate these integrals by means of numerical quadrature
formulas, then it is appropriate to eliminate the singularity of the
function $s(w)$ at $w=0$ by change of variables.
\begin{eqnarray*}
\int_{a}^{b} s(w) \: dw & = & \int_{a}^{b} \sqrt{w(1+\ak w)}(1+2\ak
w) \: dw
\\ & = &
\int_{\sqrt{a}}^{\sqrt{b}} \sqrt{1+\ak r^{2}} \, (1+2\ak r^{2}) \, 2
\, r^{2} \, dr \, , \qquad ( w = r^{2} )
\\
\int_{a}^{b} s(w-\qe) \: dw & = & \int_{a-\qe}^{b-\qe} s(\hat{w}) \:
d\hat{w}
\\ & = &
\int_{\sqrt{a-\qe}}^{\sqrt{b-\qe}} \sqrt{1+\ak r^{2}} \, (1+2\ak
r^{2}) \, 2 \, r^{2} \, dr \, , \qquad ( w = r^{2} + \qe )
\\
\int_{a}^{b} s(w+\qe) \: dw & = & \int_{a+\qe}^{b+\qe} s(\bar{w}) \:
d\bar{w}
\\ & = &
\int_{\sqrt{a+\qe}}^{\sqrt{b+\qe}} \sqrt{1+\ak r^{2}} \, (1+2\ak
r^{2}) \, 2 \, r^{2} \, dr \, . \qquad ( w = r^{2} - \qe )
\end{eqnarray*}
\subsection{Integrals related to the free streaming operator}

We recall that
\begin{eqnarray*}
g_1 \argf & = & c_x \frac{\mu \sqrt{w(1+\ak w)}}{1+2 \ak w} \, ,
\\
g_2 \argf & = & c_x \frac{\sqrt{1-\mu^2} \sqrt{w(1+\ak w)}
\cos\ph}{1+2 \ak w} \, ,
\\
g_3 \argf & = & \mbox{} - 2 c_k \frac{\sqrt{w(1+\ak w)}}{1+2 \ak w}
\left[ \mu \, E_x(t,x,y) + \sqrt{1-\mu^2} \cos\ph \, E_y(t,x,y)
\right] \, ,
\\
g_4 \argf & = & \mbox{} - c_k \frac{\sqrt{1-\mu^2}}{\sqrt{w(1+\ak
w)}}
 \left[ \sqrt{1-\mu^2} \, E_x(t,x,y) - \mu \cos\ph \, E_y(t,x,y) \right] \, ,
\\
g_5 \argf & = & c_k \frac{\sin\ph}{\sqrt{w(1+\ak w)} \sqrt{1-\mu^2}}
\, E_y(t,x,y) \, .
\end{eqnarray*}
Now, we define:
\begin{equation*}
s_{1}(w) = \frac{\sqrt{w(1+\ak w)}}{1+2 \ak w} \, , \quad s_{2}(w) =
\frac{1}{\sqrt{w(1+\ak w)}} \, .
\end{equation*}
We need to evaluate the integrals.
\begin{eqnarray*}
\int_{a}^{b} s_{1}(w) \: dw & = & \int_{a}^{b} \frac{\sqrt{w(1+\ak
w)}}{1+2\ak w} \: dw
 = \int_{\sqrt{a}}^{\sqrt{b}} \frac{\sqrt{1+\ak r^{2}}}{1+2\ak r^{2}} \, 2 \, r^{2} \, dr \, ,
\\
\int_{a}^{b} s_{2}(w) \: dw & = & \int_{a}^{b}
\frac{1}{\sqrt{w(1+\ak w)}} \: dw = \int_{\sqrt{a}}^{\sqrt{b}}
\frac{2}{\sqrt{1+\ak r^{2}}} \, dr \, ,
\end{eqnarray*}
\begin{eqnarray*}
\int_{a}^{b} s_{1}(w) \, \frac{2(w - w_{k})}{\Delta w_{k}} \: dw & =
& \int_{a}^{b} \frac{\sqrt{w(1+\ak w)}}{1+2\ak w}  \, \frac{2(w -
w_{k})}{\Delta w_{k}} \: dw
\nonumber \\
& = & \int_{\sqrt{a}}^{\sqrt{b}} \frac{\sqrt{1+\ak r^{2}}}{1+2\ak
r^{2}}
 \, \frac{2(r^{2} - w_{k})}{\Delta w_{k}}\, 2 \, r^{2} \, dr \, ,
\\
\int_{a}^{b} s_{2}(w) \, \frac{2(w - w_{k})}{\Delta w_{k}} \: dw & =
& \int_{a}^{b} \frac{1}{\sqrt{w(1+\ak w)}}  \, \frac{2(w -
w_{k})}{\Delta w_{k}} \: dw
\nonumber \\
& = & \int_{\sqrt{a}}^{\sqrt{b}} \frac{2}{\sqrt{1+\ak r^{2}}} \,
\frac{2(r^{2} - w_{k})}{\Delta w_{k}}\,  dr \,
\end{eqnarray*}
\begin{eqnarray*}
\int_{a}^{b} s_{1}(w) \left[  \frac{2(w - w_{k})}{\Delta w_{k}}
\right] ^{2} dw & = & \int_{a}^{b} \frac{\sqrt{w(1+\ak w)}}{1+2\ak
w}  \left[  \frac{2(w - w_{k})}{\Delta w_{k}} \right] ^{2} dw
\nonumber \\
& = & \int_{\sqrt{a}}^{\sqrt{b}} \frac{\sqrt{1+\ak r^{2}}}{1+2\ak
r^{2}}
 \left[  \frac{2(r^{2} - w_{k})}{\Delta w_{k}} \right] ^{2} 2 \, r^{2} \, dr \, ,
\\
\int_{a}^{b} s_{2}(w) \left[  \frac{2(w - w_{k})}{\Delta w_{k}}
\right] ^{2} dw & = & \int_{a}^{b} \frac{1}{\sqrt{w(1+\ak w)}}
\left[  \frac{2(w - w_{k})}{\Delta w_{k}} \right] ^{2} dw
\nonumber \\
& = & \int_{\sqrt{a}}^{\sqrt{b}} \frac{2}{\sqrt{1+\ak r^{2}}} \left[
\frac{2(r^{2} - w_{k})}{\Delta w_{k}} \right] ^{2}  dr \, .
\end{eqnarray*}
\subsection{Initial condition}

Since $\dm \emph{measure} (\Omega_{ijkmn}) = \Delta x_{i} \, \Delta
y_{j} \, \Delta w_{k} \, \Delta \mu_{m} \, \Delta \ph_{n}$, we have
\begin{eqnarray*}
&& \int_{\Omega_{ijkmn}}
 \Phi_h(0,x,y,w,\mu, \ph)  \, dx \, dy \, dw \, d\mu \, d\ph
 = \emph{measure} (\Omega_{ijkmn}) \, T_{ijkmn}(0) \, ,
\\
&& \int_{\Omega_{ijkmn}}
 \Phi_h(0,x,y,w,\mu, \ph)  \, \frac{2(x - x_{i})}{\Delta x_{i}} \, dx \, dy \, dw \, d\mu \, d\ph
 = \frac{1}{3} \emph{measure} (\Omega_{ijkmn}) \, X_{ijkmn}(0) \, ,
\\
&& \int_{\Omega_{ijkmn}}
 \Phi_h(0,x,y,w,\mu, \ph)  \, \frac{2(y - y_{j})}{\Delta y_{j}} \, dx \, dy \, dw \, d\mu \, d\ph
 = \frac{1}{3} \emph{measure} (\Omega_{ijkmn}) \, Y_{ijkmn}(0) \, ,
\\
&& \int_{\Omega_{ijkmn}}
 \Phi_h(0,x,y,w,\mu, \ph)  \, \frac{2(w - w_{k})}{\Delta w_{k}} \, dx \, dy \, dw \, d\mu \, d\ph
 = \frac{1}{3} \emph{measure} (\Omega_{ijkmn}) \, W_{ijkmn}(0) \, ,
\\
&& \int_{\Omega_{ijkmn}}
 \Phi_h(0,x,y,w,\mu, \ph)  \,\frac{2(\mu - \mu_{m})}{\Delta \mu_{m}} \, dx \, dy \, dw \, d\mu \, d\ph
 = \frac{1}{3} \emph{measure} (\Omega_{ijkmn}) \, M_{ijkmn}(0) \, ,
\\
&& \int_{\Omega_{ijkmn}}
 \Phi_h(0,x,y,w,\mu, \ph)  \,\frac{2(\ph - \ph_{n})}{\Delta \ph_{n}} \, dx \, dy \, dw \, d\mu \, d\ph
 = \frac{1}{3} \emph{measure} (\Omega_{ijkmn}) \, P_{ijkmn}(0) \, .
\end{eqnarray*}
If, for each $(x,y)$ and $(\mu, \ph)$,
$$\dm \Phi_h(0,x,y,w,\mu, \ph) = \left\lbrace
\begin{array}{ll}
F(x,y) \, s(w) \, \mbox{e}^{-w} & \mbox{for } w < w_{N_{w} + \ot} \\[5pt]
0 & \mbox{otherwise}
\end{array}
\right.
$$
and
$$
 F(x,y) = F_{ij} \mbox{ (constant)} \qquad \forall (x,y) \in
\left[ x_{i - \ot} , \, x_{i + \ot} \right] \times \left[ y_{j -
\ot} , \, y_{j + \ot} \right] ,
$$
then it is reasonable to assume
\begin{eqnarray}
& & T_{ijkmn}(0) = \frac{F_{ij}}{\Delta w_{k}}
    \int_{w_{k - \ot}}^{w_{k + \ot}} s(w) \, \mbox{e}^{-w} \, dw \, ,
\\[10pt]
& & X_{ijkmn}(0) = Y_{ijkmn}(0) = 0 \, ,
\\[10pt]
& & W_{ijkmn}(0)  = 3 \, \frac{F_{ij}}{\Delta w_{k}}
    \int_{w_{k - \ot}}^{w_{k + \ot}} s(w) \, \mbox{e}^{-w}
    \, \frac{2(w - w_{k})}{\Delta w_{k}} \, dw \, ,
\\[10pt]
& & M_{ijkmn}(0) = P_{ijkmn}(0) = 0 \, .
\end{eqnarray}
Recalling the definition (\ref{dens}) of the dimensionless charge
density, we have
$$
  \rho(0,x,y) = 2 \, \pi \, F(x,y)  \sum_{k=1}^{N_{w}}
\int_{w_{k - \ot}}^{w_{k + \ot}} s(w) \,  \mbox{e}^{-w} \, dw \, ,
$$
which gives a relationship between $F(x,y)$ and the initial charge
density $\rho$.

\subsection{Hydrodynamical variables}

If $\dm (x, y) \in \left[ x_{i - \ot}, x_{i + \ot} \right] \times
\left[ y_{j - \ot}, y_{j + \ot} \right]$, then
\begin{eqnarray}
\rho_h(t,x,y) & =& \sum_{k=1}^{N_{w}} \sum_{m=1}^{N_{\mu}}
\sum_{n=1}^{N_{\ph}}
  \left[ T_{ijkmn}(t) +
  X_{ijkmn}(t) \, \frac{2(x - x_{i})}{\Delta x_{i}} \right.
\nonumber \\
& & \left. \hspace{50pt} \mbox{} +
  Y_{ijkmn}(t) \, \frac{2(y - y_{j})}{\Delta y_{j}} \right]
 \Delta w_{k} \Delta \mu_{m} \, \Delta \ph_{n} \, .
\end{eqnarray}
We define
\begin{eqnarray*}
\hat{T}_{ij}(t) & = & \sum_{k=1}^{N_{w}} \sum_{m=1}^{N_{\mu}}
\sum_{n=1}^{N_{\ph}} T_{ijkmn}(t) \, \Delta w_{k} \Delta \mu_{m} \,
\Delta \ph_{n}
\\
\hat{X}_{ij}(t) & = & \sum_{k=1}^{N_{w}} \sum_{m=1}^{N_{\mu}}
\sum_{n=1}^{N_{\ph}} X_{ijkmn}(t) \, \Delta w_{k} \Delta \mu_{m} \,
\Delta \ph_{n}
\\
\hat{Y}_{ij}(t) & = & \sum_{k=1}^{N_{w}} \sum_{m=1}^{N_{\mu}}
\sum_{n=1}^{N_{\ph}} Y_{ijkmn}(t) \, \Delta w_{k} \Delta \mu_{m} \,
\Delta \ph_{n} \, .
\end{eqnarray*}
Therefore, for every $\dm (x, y) \in \left[ x_{i - \ot}, x_{i + \ot}
\right] \times \left[ y_{j - \ot}, y_{j + \ot} \right]$,
\begin{equation}
\rho_h(t,x,y) =
 \hat{T}_{ij}(t) + \hat{X}_{ij}(t)  \, \frac{2(x - x_{i})}{\Delta x_{i}} +
\hat{Y}_{ij}(t)  \, \frac{2(y - y_{j})}{\Delta y_{j}} \, .
\end{equation}
For every $\dm (x, y) \in \left[ x_{i - \ot}, x_{i + \ot} \right]
\times \left[ y_{j - \ot}, y_{j + \ot} \right]$, the approximate
momentum in $x$-direction is
\begin{eqnarray}
&& \sum_{k=1}^{N_{w}} \sum_{m=1}^{N_{\mu}} \sum_{n=1}^{N_{\ph}}
\Delta \ph_{n} \left[
 g_{1, k m} \, T_{ijkmn}(t) + g_{1 w, k m} \, W_{ijkmn}(t)  + g_{1 \mu, k m} \,M_{ijkmn}(t)
\right] \nonumber
\\
&& \mbox{ } + \left[ \sum_{k=1}^{N_{w}} \sum_{m=1}^{N_{\mu}}
\sum_{n=1}^{N_{\ph}}  \Delta \ph_{n} \,
  g_{1, k m} \, X_{ijkmn}(t) \right]  \frac{2(x - x_{i})}{\Delta x_{i}} \nonumber
\\
&& \mbox{ } + \left[ \sum_{k=1}^{N_{w}} \sum_{m=1}^{N_{\mu}}
\sum_{n=1}^{N_{\ph}}  \Delta \ph_{n} \,
  g_{1, k m} \, Y_{ijkmn}(t) \right]  \frac{2(y - y_{j})}{\Delta y_{j}} \, ,
\end{eqnarray}
and the approximate momentum in $y$-direction is
\begin{eqnarray}
 &&
\sum_{k=1}^{N_{w}} \sum_{m=1}^{N_{\mu}} \sum_{n=1}^{N_{\ph}} \left[
 g_{2, k m n} \, T_{ijkmn}(t) + g_{2 w, k m n} \, W_{ijkmn}(t)  + g_{2 \mu, k m n} \,M_{ijkmn}(t)
 + g_{2 \ph, k m n} \, P_{ijkmn}(t)
\right] \nonumber
\\
&& \mbox{ } + \left[ \sum_{k=1}^{N_{w}} \sum_{m=1}^{N_{\mu}}
\sum_{n=1}^{N_{\ph}}
  g_{2, k m n} \, X_{ijkmn}(t) \right]  \frac{2(x - x_{i})}{\Delta x_{i}} \nonumber
\\
&& \mbox{ } + \left[ \sum_{k=1}^{N_{w}} \sum_{m=1}^{N_{\mu}}
\sum_{n=1}^{N_{\ph}}
  g_{2, k m n} \, Y_{ijkmn}(t) \right]  \frac{2(y - y_{j})}{\Delta y_{j}} \, ,
\end{eqnarray}
where
\begin{eqnarray*}
&& g_{1,k m} =  \iw \imu g_{1}(w,\mu) \, dw \, d\mu
\\
&& g_{1 w,k m} =  \iw \imu g_{1}(w,\mu)  \, \frac{2(w -
w_{k})}{\Delta w_{k}} \, dw \, d\mu
\\
&& g_{1 \mu,k m} =  \iw \imu g_{1}(w,\mu) \, \frac{2(\mu -
\mu_{m})}{\Delta \mu_{m}} \, dw \, d\mu
\\
&& g_{2,k m n} = \iw \imu \iphi g_{2}(w,\mu, \ph) \, dw \, d\mu  \,
d\ph
\\
&& g_{2 w,k m n} = \iw \imu \iphi g_{2}(w,\mu, \ph)  \, \frac{2(w -
w_{k})}{\Delta w_{k}} \, dw \, d\mu  \, d\ph
\\
&& g_{2 \mu,k m n} =  \iw \imu \iphi g_{2}(w,\mu, \ph) \,
\frac{2(\mu - \mu_{m})}{\Delta \mu_{m}} \, dw \, d\mu  \, d\ph
\\
&& g_{2 \ph,k m n} = \iw \imu \iphi g_{2}(w,\mu, \ph) \, \frac{2(\ph
- \ph_{n})}{\Delta \ph_{n}} \, dw \, d\mu  \, d\ph \, .
\end{eqnarray*}
Analogously, the energy multiplied by the charge density
$\rho_h(t,x,y)$ is
\begin{eqnarray}
&& \sum_{k=1}^{N_{w}} \sum_{m=1}^{N_{\mu}} \sum_{n=1}^{N_{\ph}}
\Delta \ph_{n} \, \Delta \mu_{m} \left[
 w_{k} \, \Delta w_{k}  \, T_{ijkmn}(t) + \frac{(\Delta w_{k})^{2}}{6} \,  W_{ijkmn}(t)
\right] \nonumber
\\
&& \mbox{ } + \left[ \sum_{k=1}^{N_{w}} \sum_{m=1}^{N_{\mu}}
\sum_{n=1}^{N_{\ph}}  \Delta \ph_{n} \, \Delta \mu_{m}
 \, w_{k} \, \Delta w_{k}  \, X_{ijkmn}(t) \right]  \frac{2(x - x_{i})}{\Delta x_{i}} \nonumber
\\
&& \mbox{ } + \left[ \sum_{k=1}^{N_{w}} \sum_{m=1}^{N_{\mu}}
\sum_{n=1}^{N_{\ph}}  \Delta \ph_{n} \, \Delta \mu_{m}
 \,  w_{k} \, \Delta w_{k}  \, Y_{ijkmn}(t) \right]  \frac{2(y - y_{j})}{\Delta y_{j}} \, ,
\end{eqnarray}
for every $\dm (x, y) \in \left[ x_{i - \ot}, x_{i + \ot} \right]
\times \left[ y_{j - \ot}, y_{j + \ot} \right]$.

\end{document}